\documentclass[final,1p,times]{elsarticle}

\usepackage{amssymb}
\usepackage{amsthm}
\usepackage{amsmath}
\usepackage[scr=rsfs]{mathalpha}
\usepackage{lineno}
\usepackage{algorithm}
\usepackage{algpseudocode}
\usepackage{hyperref}
\hypersetup{
    colorlinks=true,
    linkcolor=blue,
    filecolor=magenta,      
    urlcolor=cyan,
    pdftitle={JCP-tumor},
    pdfpagemode=FullScreen,
    }

\journal{Journal of Computational Physics}

\begin{document}

\begin{frontmatter}



\title{Nonlinear simulation of vascular tumor growth with chemotaxis and the control of necrosis}


\author[inst1]{Min-Jhe Lu}
\affiliation[inst1]{organization={Department of Applied Mathematics, Illinois Institute of Technology},
            city={Chicago},
            postcode={60616}, 
            state={Illinois},
            country={United States}}
            
\author[inst3]{Wenrui Hao}

\affiliation[inst3]{organization={Department of Mathematics, Pennsylvania State University},
            city={University Park},
            postcode={16802}, 
            state={Pennsylvania},
            country={United States}}

\author[inst1]{Chun Liu}
\author[inst2]{John Lowengrub}

\affiliation[inst2]{organization={Departments of Mathematics and Biomedical Engineering, Center for Complex Biological Systems, Chao Family Comprehensive Cancer Center, University of California at Irvine},
            city={Irvine},
            postcode={92617}, 
            state={California},
            country={United States}}
            
\author[inst1]{Shuwang Li\corref{cor1}}
\ead{sli15@iit.edu}
\cortext[cor1]{Corresponding author} 

\begin{abstract}
In this paper, we develop a sharp interface tumor growth model in two dimensions to study the effect of both the intratumoral structure using a controlled necrotic core and the extratumoral nutrient supply from vasculature on tumor morphology. We first show that our model extends the benchmark results in the literature using linear stability analysis. Then we solve this generalized model numerically using a spectrally accurate boundary integral method in an evolving annular domain, not only with a Robin boundary condition on the outer boundary for the nutrient field which models tumor vasculature, but also with a static boundary condition on the inner boundary for pressure field which models the control of tumor \textit{necrosis}. The discretized linear systems for both pressure and nutrient fields are shown to be well-conditioned through tracing GMRES iteration numbers. Our nonlinear simulations reveal the stabilizing effects of \textit{angiogenesis} and the destabilizing ones of \textit{chemotaxis} and \textit{necrosis} in the development of tumor morphological instabilities if the necrotic core is fixed in a circular shape. When the necrotic core is controlled in a non-circular shape, the stabilizing effects of \textit{proliferation} and the destabilizing ones of \textit{apoptosis} are observed. Finally, the values of the nutrient concentration with its fluxes and the pressure level with its normal derivatives, which are solved accurately at the boundaries, help us to characterize the corresponding tumor morphology and the level of the biophysical quantities on interfaces required in keeping various shapes of the necrotic region of the tumor. Interestingly, we notice that when the necrotic region is fixed in a 3-fold non-circular shape, even if the initial shape of the tumor is circular, the tumor will evolve into a shape corresponding to the 3-fold symmetry of the shape of the fixed necrotic region.
\end{abstract}



\begin{keyword}
Tumor growth with chemotaxis
\sep angiogenesis 
\sep control of necrosis 
\sep boundary integral method 
\sep Robin boundary condition 
\sep annular domain.
\end{keyword}
\end{frontmatter}
\section{Introduction}

The growth of a solid tumor is characterized by several increasingly aggressive stages of development. In the first stage, \textit{carcinogenesis}, genetic mutations result in the occurrence of abnormal cell proliferation and cell apoptosis. In the second stage, \textit{avascular growth}, the tumor receives nutrients (\textit{e.g.}, oxygen) by diffusion through the surrounding tissues. As the size of tumor increases, \textit{necrosis} occurs in the center of tumor since the interior tumor cells start to die due to the lack of nutrient supply. The necrotic core is thus formed, and tumor cells within this region will secrete Tumor Angiogenesis Factor (TAF) in order to grasp more nutrient supply. The control of the development of necrotic region is thus significant. 

After the development of a tumor-induced neovasculature from  \textit{angiogenesis}, the tumor arrives at the stage of \textit{vascular growth} and receives nutrients from the vasculature. Along with this development, the heterogeneous nutrient distributions could induce diffusional instability through nonuniform rates of cell proliferation, apoptosis, and migration. Tumor morphological instability, in turn, is capable of bringing more available nutrients to the tumor by increasing its surface-to-volume ratio.  In particular, regions where instabilities first occur tend to grow at a faster rate than the rest of the tumor tissue (\textit{e.g.}, differential growth) that further enhances the instabilities and leads to complicated tumor morphologies, which was shown to increase the invasive behavior of tumor, \textit{e.g.}, in \cite{cristini2005morphologic}. In order to study the evolution of tumor morphologies for clinical purposes, it is therefore important to build a mathematical model not only tracking the tumor interface accurately but also incorporating the \textit{angiogenesis} process under the control of the necrotic region.

The mathematical modeling and the nonlinear simulation of the process of tumor growth has been studied since the mid-1960s (see, for example, the reviews
\cite{araujo2004history,fasano2006mathematical,roose2007mathematical,bellomo2008selected,lowengrub2009nonlinear,byrne2010dissecting,byrne2012mathematical,KimOthmer2015,alfonso2017biology,Yankeelov2018} and books \cite{cristini2010multiscale,cristini2017introduction}). 
Following the early biomechanical models of avascular tumor growth proposed by Greenspan \cite{greenspan1976growth}, 
the bifurcation analysis 
%
\cite{friedman2001existence,friedman2001symmetry,friedman2006bifurcation,Hu20071,Hu20072,friedman2008stability,zhao2020symmetry,zhao2020impact},
numerical simulations, and computational modeling \cite{cristini2003nonlinear,cristini2009,fritz2019local,mcdougall2002mathematical,mcdougall2006mathematical,foo2011stochastic,hillen2013tumor,pham2018nonlinear} have contributed significantly to the tumor modeling area. 
Recently, tumor growth models with a necrotic core have also been developed and analyzed via the bifurcation theory \cite{cui2001analysis,hao2012bifurcation,kohlmann2012necrotic,hao2012continuation,wu2019bifurcation,zhuang2018analysis,song2021stationary}. 

In this paper, based on the tumor model with a complex far-field geometry \cite{lu2020complex}, we extend the tumor microenvironment with  a heterogeneous distribution of vasculature to include  \textit{angiogenesis, necrosis} and \textit{chemotaxis}. More specifically, we develop a novel boundary integral formulation featuring both the Robin boundary condition on tumor boundary for the nutrient field, signifying the \textit{angiogenesis} effect, and a static boundary condition on the necrotic boundary for pressure field, signifying the control of \textit{necrosis}. Such boundary integral formulation comes from a quasi-steady assumption of the nutrient field and is justified through an estimation of the taxis time scale versus the diffusion time scale. The main goals of this paper are to analyze the linear stability of the tumor model with \textit{angiogenesis, necrosis} and \textit{chemotaxis}, and numerically simulate the fully nonlinear dynamics, using a novel formulation with boundary integral method in two dimensions. 

The work presented in this paper is unique in the following aspects. First, we consider the boundary of the necrotic core fixed with a threshold nutrient level (also fixed), which allows us to study the nutrient distribution which induces morphological instabilities under the scenario that the necrotic core is controlled. We remarked that the role of necrosis in destabilizing the tumor morphology is also investigated in \cite{macklin2007nonlinear}. Second, from the numerical perspective, we develop a new boundary integral method (BIM), which naturally incorporates the Robin boundary condition which models \textit{angiogenesis} without approximation errors introduced by spatial meshes at tumor boundary. The integral equations uniquely determine the nutrient concentration and the normal derivative of the pressure at the tumor boundary, and also determine the pressure level and the nutrient flux across the necrotic boundary.  
Such a sharp interface model that solves integral equations with spectral accuracy enables us to accurately track the evolution of the tumor and all the biophysical quantities required on both interfaces. Those quantities on the interfaces are important for us to understand the mechanism of the control of tumor \textit{necrosis}. Lastly, we investigate the significance of the factors which contributes to \textit{angiogenesis, necrosis} and \textit{chemotaxis} on the fully nonlinear tumor dynamics.

The paper is organized as follows. In Section \ref{sec:model}, we formulate the sharp interface model. In Section \ref{sec:Non-dimensionalization}, we nondimensionalize the resulting systems. In Section \ref{sec:bimreformulation}, we develop the BIM formulation and summarize our numerical method. In Section \ref{sec:LA}, we analyze the linear stability of the system. In Section \ref{sec:Result}, we present our simulation results of the fully nonlinear system including a numerical convergence study, a comparison with linear analysis results, and parameter studies under various effects. The conclusion is presented in Section \ref{sec:conclusions}. In Appendices A, B, C, and D, we give a complete derivation of our linear stability analysis, the details of our numerical method including layer potential evaluations for boundary integrals, the small-scale decomposition to remove the stiffness from the high-order derivatives in high curvature region on the interface and the semi-implicit time-stepping scheme to evolve the tumor interface.

\section{Mathematical model}
\label{sec:model}
\begin{figure}
 \centering
  \label{fig:domain}\includegraphics[width=0.8\textwidth]{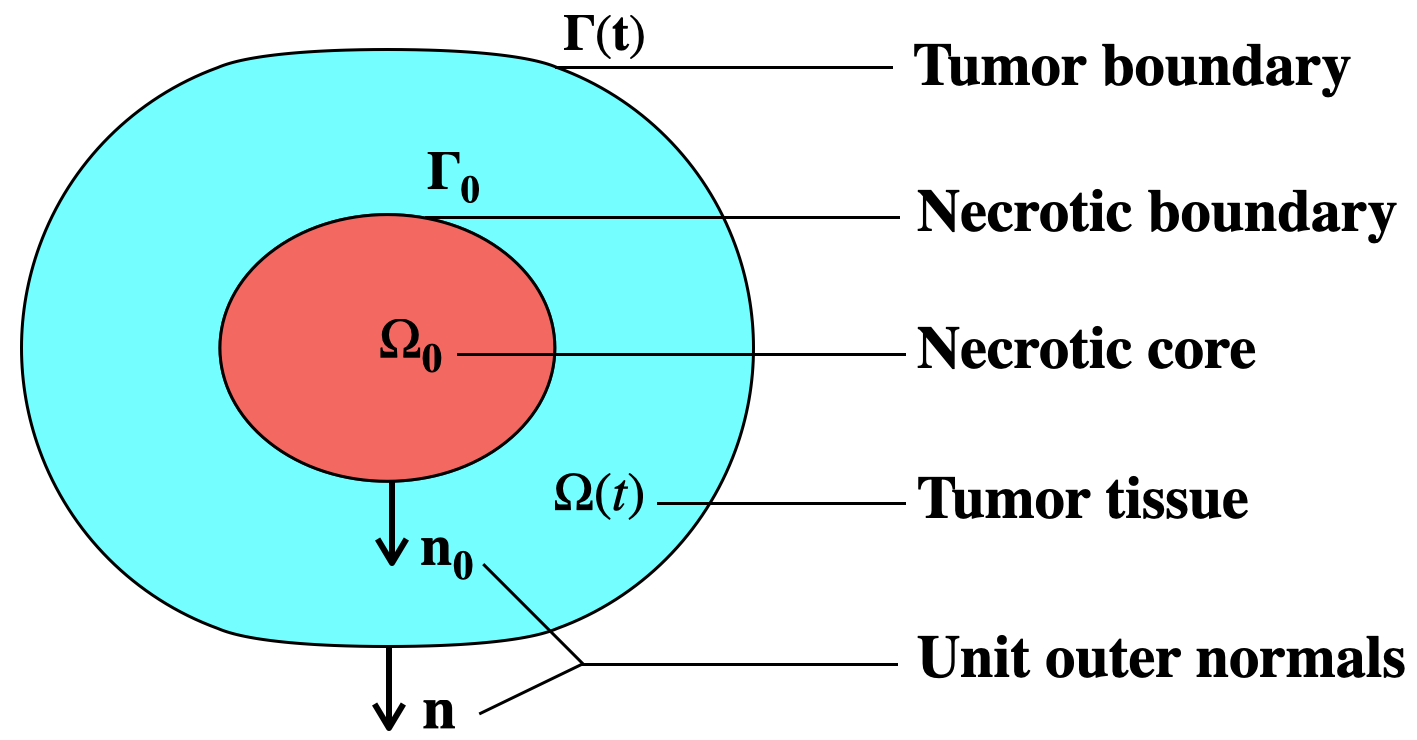}
  \caption{Illustration of the computation domain.}
\end{figure}

\paragraph{\textbf{Computation domain}}
As illustrated in Fig. \ref{fig:domain}, let $\Omega_0$ be the necrotic core, $\Omega(t)$ be the tumor tissue, $\Gamma_0$ be the controlled necrotic boundary and $\Gamma(t)$ be the tumor boundary. 

\paragraph{\textbf{Nutrient field}}
The nutrient field in $\Omega(t)$ is governed by:
\begin{equation}\label{nut rde}
\sigma_t=D\Delta \sigma-\lambda\sigma\quad \text{in }\Omega(t),
\end{equation}
where $D,\lambda$ are the diffusion constant and uptake rate, respectively.
We assume Dirichlet boundary condition on the necrotic boundary: 
\begin{equation}
\sigma=\sigma^N\quad\text{on }\Gamma_0,
\end{equation}
where $\sigma^N$ is the constant nutrient level at the necrotic boundary.
We assume Robin boundary condition on the tumor boundary:
\begin{equation}
\frac{\partial\sigma}{\partial\mathbf{n}}+\beta(\sigma-\overline\sigma)=0\quad\text{on }\Gamma(t),
\end{equation}
where $\mathbf{n}$ is the outward normal, $\bar{\sigma}$ is the constant nutrient level outside the tumor, $\beta$ is the rate of nutrient supply to the tumor, which reflects the extent of angiogenesis. We remark that we can also use the Robin boundary condition to replace the Dirichlet boundary condition on the necrotic boundary, since as $\beta$ is large enough, the Robin boundary condition will converge to the Dirichlet boundary condition, which is used for tumor boundary in \cite{cristini2003nonlinear}.

\paragraph{\textbf{Pressure field}}
To introduce chemotaxis, the directed cell migration up gradients of nutrients, we use the Chemo-Darcy's law:
\begin{equation}\label{chemodarcy}
\mathbf{u}=-\mu \nabla p + \chi_{\sigma}\nabla \sigma\quad\text{in } \Omega(t),
\end{equation}
where $\mathbf{u}$ is the tumor cell velocity, $\mu$ is the cell mobility and $\chi_\sigma$ is the chemotaxis coefficient.
The mass conservation:
\begin{equation}\label{massconserv}
\nabla  \cdot \mathbf{u}=\lambda_M\frac{\sigma}{\sigma_\infty}-\lambda_A\quad\text{in }\Omega(t),
\end{equation}
where $\lambda_M$, $\lambda_A$ are the rates of mitosis (cell birth) and apoptosis (cell death), respectively.
Applying Eq. $\eqref{chemodarcy}$ into Eq.  $\eqref{massconserv}$, we have
\begin{equation}
-\mu\Delta p=\left(\frac{\lambda_M}{\sigma^\infty}-\chi_\sigma\right)\sigma-\lambda_A\quad\text{in }\Omega(t).
\end{equation}
We assume a static boundary condition on the necrotic boundary:
\begin{equation}\label{fix nec}
0=-\mu \frac{\partial p}{\partial \mathbf{n}_0}+\chi_{\sigma}\frac{\partial \sigma}{\partial \mathbf{n}_{0}}\quad\text{on }\Gamma_0,
\end{equation}
which corresponds to our assumption that the necrotic boundary is fixed.

The Laplace-Young condition is assumed on the tumor boundary:
\begin{equation}\label{eq6}
p=\gamma \left.\kappa\right|_{\Gamma(t)}\quad\text{on }\Gamma(t),
\end{equation}
where $\gamma$ is the constant representing cell-cell adhesions and $\left.\kappa\right|_{\Gamma(t)}$ is the mean curvature of the curve $\Gamma(t)$.

\paragraph{\textbf{Equation of motion}}
The equation of motion for the interface $\Gamma(t)$ is given by:
\begin{equation}\label{eq8}
V\equiv\mathbf u
\cdot 
\mathbf{n}=-\mu \left.\frac{\partial p}{\partial \mathbf{n}}\right|_{\Gamma(t)} 
+ \chi_\sigma \left.\frac{\partial \sigma}{\partial \mathbf{n}}\right|_{\Gamma(t)}\quad\text{on }\Gamma(t).
\end{equation}

\section{Non-dimensionalization}
\label{sec:Non-dimensionalization}
We introduce the diffusion length $L$, the intrinsic taxis time scale $\lambda_\chi^{-1}$, and the characteristic pressure $p_s$ by:
\begin{equation}\label{eq9}
L=\sqrt{\frac{D}{\lambda}}, \quad
\lambda_\chi=\frac{\overline{\chi_{\sigma}}\sigma^{\infty}}{L^2}, \quad 
p_s=\frac{\lambda_{\chi}L^2}{\mu},
\end{equation}
where $\overline{\chi_\sigma}$ is a characteristic taxis coefficient. The length scale $L$ and the time scale ${\lambda_\chi}^{-1}$ are used to non-dimensionalize the space and time variables by $\mathbf{x}=L\widetilde{\mathbf{x}}$, $t={\lambda_\chi}^{-1}\widetilde{t}$.
Define
\begin{equation}\label{eq10}
\widetilde{\sigma}=\frac{\sigma}{\sigma^{\infty}} , \quad
\underline{\sigma}=\frac{\sigma^N}{\sigma^{\infty}} , \quad
\widetilde{p}=\frac{p}{p_s}, \quad
\widetilde{\chi_{\sigma}}=\frac{\chi_{\sigma}}{\overline{\chi_{\sigma}}}, \quad
\widetilde{\beta}=L\beta.
\end{equation}
Since taxis occurs more slowly than diffusion (\textit{e.g.} minutes vs hours), we assume  $\lambda_{\chi} \ll \lambda$, which leads to a quasi-steady reaction-diffusion equation for the nutrient field. We remark that by the term ``taxis'', we mean taxis of tumor cells up nutrient gradients, as embodied in Eq. \eqref{chemodarcy}. Then Eq. \eqref{nut rde} becomes $\varepsilon\widetilde\sigma_{\widetilde t}=\widetilde\Delta\widetilde\sigma-\widetilde \sigma$, where $\varepsilon=\frac{\lambda_\chi}{\lambda}\approx\frac{T_\text{diffusion}}{T_\text{taxis}}$. For the nutrient diffusion time scale $T_\text{diffusion}$, typically it can be assumed to occur in the scale of minutes, say 1 minute (see p.226 in \cite{friedman2006cancer}). For the taxis time scale $T_\text{taxis}$, we can estimate it through dividing the diameter  of the diffusion-limited tumor spheroid by the speed of migration of tumor cell up chemical gradients.  
For the tumor diameter, as summarized in \cite{grimes2014method}: ``oxygen diffusion limits are typically 100–200 $\mu m$'', and here we take the average 150 $\mu m$. For the speed of tumor migration, as summarized in \cite{roussos2011chemotaxis}: ``Some carcinoma cells with an amoeboid morphology can move at high speeds inside the tumours ($\sim 4 \mu m\text{ min}^{-1} $) ... At the other end of the range of modes of motility, ... mesenchymal migration ... (0.1--1 $\mu m\text{ min}^{-1} $)'', and here we take the average of the two types $\sim 2\mu m \text{ min}^{-1}$. Therefore, the taxis time scale can be estimated as $\frac{150 \mu m}{2\mu m \text{ min}^{-1}}=1.25$ hour.
Hence we have $\displaystyle\varepsilon\approx\frac{T_\text{diffusion}}{T_\text{taxis}}\approx\frac{1\  minute}{1\ hour} \ll 1$.
The dimensionless system is thus given by:
\paragraph{\textbf{Nutrient field}}
We have governing equations for the nutrient field:

\begin{equation}
\left\{\begin{aligned}
\widetilde{\Delta} \tilde{\sigma} &=\widetilde{\sigma} & & \text { in } \Omega(t), \\
\tilde{\sigma} &=\underline{\sigma} & & \text { on } \Gamma_{0}, \\
\frac{\widetilde{\partial} \widetilde{\sigma}}{\widetilde{\partial} \widetilde{\mathbf{n}}} &=\widetilde{\beta}(1-\widetilde{\sigma}) && \text { on } \Gamma(t),
\end{aligned}\right.
\end{equation}
where $\widetilde\beta$ (angiogenesis factor) represents the extent of angiogenesis.
\paragraph{\textbf{Pressure field}}
\begin{itemize}
    \item Non-dimensional Chemo-Darcy's law.
    \begin{equation}
\mathbf{\widetilde{u}}=-\widetilde{\nabla} \widetilde{p} + \widetilde{\chi_{\sigma}} \widetilde{\nabla} \widetilde{\sigma}\quad\text{in }\Omega(t),
\end{equation}
where $\widetilde{\chi_\sigma}$ (chemotaxis constants) represents taxis effect.
\item Conservation of tumor mass.
\begin{equation}
\widetilde{\nabla} \cdot \widetilde{\mathbf{u}}=\mathcal{P}\left(\widetilde{\sigma}- \mathcal{A}\right)\quad\text{in }\Omega(t),
\end{equation}
where $\displaystyle \mathcal{P}=\frac{\lambda_M}{\lambda_\chi}$ (proliferation rate) represents the rate of cell mitosis relative to taxis, $\displaystyle \mathcal{A}=\frac{\lambda_A}{\lambda_M}$ (apoptosis rate) represents apoptosis relative to cell mitosis.
\item Boundary conditions.
\begin{align}
\left.\frac{\widetilde\partial \widetilde p}{\widetilde\partial \widetilde{\mathbf{n}}_0}\right|_{\Gamma_0}&=\widetilde{\chi_{\sigma}}\left.\frac{\widetilde\partial\widetilde\sigma}{\widetilde\partial \widetilde{ \mathbf{n}_0}}\right|_{\Gamma_0}&\text{on }&\Gamma_0,\label{psbc1}\\
\left.\widetilde p\right|_{\Gamma(t)}&=\widetilde{\mathcal{G}}^{-1}\left.\widetilde\kappa\right|_{\Gamma(t)}&\text{on }&\Gamma(t),\label{psbc2}
\end{align}
where $\displaystyle \widetilde{\mathcal{G}}^{-1} 
=\frac{\mu\gamma}{\lambda_\chi L^3}$ represents the relative strength of cell-cell interactions (adhesion).
\end{itemize}
Therefore, we have governing equations for pressure field:
\begin{equation}
\left\{
\begin{array}{ccc}
\begin{aligned}
-\widetilde \Delta \widetilde p&=\mathcal{P}(\widetilde\sigma-\mathcal{A}) -\widetilde\chi_\sigma\widetilde\sigma&\text{in }&\Omega(t),\\
\left.\frac{\widetilde\partial \widetilde p}{\widetilde\partial \widetilde{\mathbf{n}}_0}\right|_{\Gamma_0}&=\widetilde{\chi_{\sigma}}\left.\frac{\widetilde\partial\widetilde\sigma}{\widetilde\partial \widetilde{ \mathbf{n}_0}}\right|_{\Gamma_0}&\text{on }&\Gamma_0,\\
\left.\widetilde{p}\right|_{\Gamma(t)}&=\widetilde{\mathcal{G}}^{-1}\left.\widetilde\kappa\right|_{\Gamma(t)
}&\text{on }&\Gamma(t).
\end{aligned}
\end{array}
\right.
\end{equation}
\paragraph{\textbf{Equation of motion}}
\begin{equation}
\widetilde{V}=-\left.\frac{\widetilde{\partial}\widetilde{p}}{\widetilde{\partial}\widetilde{\mathbf{n}}}\right|_{\Gamma(t)}
+\widetilde{\chi_{\sigma}}\left.\frac{\widetilde{\partial}\widetilde{\sigma}}{\widetilde{\partial}\widetilde{\mathbf{n}}}\right|_{\Gamma(t)}\quad\text{on }\Gamma(t).
\end{equation}
Equivalently, we have
\begin{equation}
\widetilde{V}=-\left.\frac{\widetilde{\partial}\widetilde{p}}{\widetilde{\partial}\widetilde{\mathbf{n}}}\right|_{\Gamma(t)}
+\widetilde{\chi_{\sigma}}\widetilde\beta(1-\widetilde\sigma)\quad\text{on }\Gamma(t).
\end{equation}

\section{Boundary Integral Method (BIM) reformulation}
\label{sec:bimreformulation}

Recall that we have Poisson equation for the non-dimensional pressure $\widetilde{p}$:
\begin{equation}
-\widetilde \Delta \widetilde p=(\mathcal{P}-\widetilde\chi_\sigma)\widetilde\sigma -\mathcal{P}\mathcal{A}\quad\text{in }\Omega(t).
\end{equation}
Consider an algebraic transformation in $\Omega(t)$:
\begin{equation}\label{eq:modpressure}
\overline{p}=\widetilde{p}+(\mathcal{P}-\widetilde{\chi_{\sigma}})\widetilde\sigma - \mathcal{P}\mathcal{A}\frac{\mathbf{\widetilde{x}}\cdot \mathbf{\widetilde{x}}}{2 d}\quad\text{in } \Omega(t),
\end{equation}
which satisfies:
\begin{equation}
-\widetilde{\Delta} \overline{p}= -\widetilde{\Delta} \widetilde{p} - (\mathcal{P}-\widetilde{\chi_{\sigma}})\widetilde{\sigma}+\mathcal{P}\mathcal{A}=0\quad\text{in }\Omega(t),
\end{equation}
where $d$ is the dimension of $\mathbf{R}^d \supseteq \Omega(t)$, and $\overline{p}$ is the modified pressure.
Dropping all tildes and overbars for brevity, we have the Laplace equation for the modified pressure $p$:
\begin{equation}\label{pressurefield}
\left\{
\begin{array}{ccc}
\begin{aligned}
\Delta p&=0
&\text{in }&\Omega(t),\\
\left.\frac{\partial p}{\partial \mathbf{n}_0}\right|_{\Gamma_0}&=\mathcal{P}\left.\frac{\partial\sigma}{\partial\mathbf{n}_0}\right|_{\Gamma_0} -\mathcal{P}\mathcal{A}\left.\frac{\mathbf{n}_0\cdot\mathbf{x}}{d}\right|_{\Gamma_0}
&\text{on }&\Gamma_0,\\
\left.p\right|_{\Gamma(t)}&={\mathcal{G}}^{-1}\left.\kappa\right|_{\Gamma(t)}
+(\mathcal{P}-{\chi_{\sigma}})\left.\sigma\right|_{\Gamma(t)}
-\mathcal{P}\mathcal{A}\left.\frac{\mathbf{x}\cdot \mathbf{x}}{2 d}\right|_{\Gamma(t)}
&\text{on }&\Gamma(t).
\end{aligned}
\end{array}
\right.
\end{equation}
where $\sigma$ satisfies modified Helmholtz equations in the annular domain $\Omega(t)$:
\begin{equation}\label{nutrientfield}
\left\{
\begin{array}{ccc}
\begin{aligned}
{\Delta}{\sigma}&={\sigma}&\text{in }&\Omega(t),\\
\left.\sigma\right|_{\Gamma_0}&=\underline\sigma &\text{on }&\Gamma_0, \\
\left.\frac{\partial\sigma}{\partial\mathbf{n}}\right|_{\Gamma(t)}&=\beta(1-\left.\sigma\right|_{\Gamma(t)}) &\text{on }&\Gamma(t).
\end{aligned}
\end{array}
\right.
\end{equation}
The equations of motion for $\Gamma(t)$ are determined by:
\begin{equation}
V=-\left.\frac{\partial p}{\partial \mathbf{n}}\right|_{\Gamma(t)}
+\mathcal{P}\left.\frac{\partial \sigma}{\partial \mathbf{n}}\right|_{\Gamma(t)}
-\mathcal{P}\mathcal{A} \left.\frac{\mathbf{n}\cdot \mathbf{x}}{d}\right|_{\Gamma(t)} \quad\text{on }\Gamma(t).
\end{equation}
Equivalently, we have
\begin{equation}\label{normal velocity}
V=-\left.\frac{\partial p}{\partial \mathbf{n}}\right|_{\Gamma(t)}
-\mathcal{P}
\left(\mathcal{A} \left.\frac{\mathbf{n}\cdot \mathbf{x}}{d}\right|_{\Gamma(t)}
-\beta\left(1-\left.\sigma\right|_{\Gamma(t)}\right)\right) \quad\text{on }\Gamma(t),
\end{equation}
From potential theory, the solutions to Eqs. $\eqref{pressurefield}$ and $\eqref{nutrientfield}$ can be represented as boundary integrals with single layer and double layer potentials. We use direct BIM formulations to both $\sigma$ and $p$ in the annular domain $\Omega(t)$:
\subsection{Direct BIM for the modified Helmholtz equation}
Consider the Green's function for modified Helmholtz equations in $\Omega(t)$:
\begin{equation}\label{mHHgreen}
\Delta G_{\mathbf{x}_*} - G_{\mathbf{x}_*} = -\delta_{\mathbf{x}_*}\quad\text{in }\Omega(t),
\end{equation}
where $G_{\mathbf{x}_*}=G_\mathbf{x_*}(\mathbf{x}')$, $\mathbf{x}_*,\mathbf{x}'\in \Omega(t)$ are the source point and field point, respectively, and $\delta_{\mathbf{x}_*}(\mathbf{x}')$ is the Dirac delta function. 
The fundamental solution to Eq. $\eqref{mHHgreen}$ is
\begin{equation}\label{eq34}
G_{\mathbf{x}_*}(\mathbf{x}')=\frac{1}{2 \pi} K_0(r),
\end{equation}
where $K_0$ is a modified Bessel function of the second kind, $r\equiv|\mathbf{x}_*-\mathbf{x}'|$.
Multiplying the $1^{st}$ equation in Eq. $\eqref{nutrientfield}$ by $G_{\mathbf{x}_*}$ and Eq. $\eqref{mHHgreen}$ by $-\sigma$ and summing them up, we obtain
\begin{equation}\label{Greenidsigma}
G_{\mathbf{x}_*} \Delta \sigma - \sigma \Delta G_{\mathbf{x}_*}=\sigma \delta_{\mathbf{x}_*}\quad\text{in }\Omega(t).
\end{equation}
Integrating Eq. $\eqref{Greenidsigma}$ over $\Omega(t)$ and using Green's second identity give
\begin{equation}\label{greenint}
\int_{\Gamma_0\cup\Gamma(t)} 
\left(G_{\mathbf{x}_*} \frac{\partial\sigma'}{\partial\mathbf{n}_{\star}'}
-\sigma' \frac{\partial G_{\mathbf{x}_*}}{\partial \mathbf{n}_\star'}\right) ds'=\sigma(\mathbf{x}_*)\quad \forall \mathbf{x}_*\in\Omega(t),
\end{equation}
where the symbol prime denotes the evaluation on field points, for example, $\mathbf{n}_\star'=\mathbf{n}_\star(\mathbf{x}')$ is the unit exterior (w.r.t. the domain $\Omega(t)$ enclosed by its boundaries) normal on the corresponding boundaries, i.e., $\mathbf{n}_\star'=-\mathbf{n}_0'$ on $\Gamma_0$, $\mathbf{n}_\star'= \mathbf{n}'$ on $\Gamma(t)$ and $\mathbf{n}_0'=\mathbf{n}_0(\mathbf{x}'),\mathbf{n}'=\mathbf{n}(\mathbf{x}')$ are the unit outer (pointing away from origin) normals on $\Gamma_0$ and $\Gamma(t)$, respectively. 
Letting $\mathbf{x}_* \to \mathbf{x}_0 \in \Gamma_0, \mathbf{x}_* \to \mathbf{x} \in \Gamma(t)$  in Eq. $\eqref{greenint}$ and using the $2^{nd}$ and the $3^{rd}$ equation in Eq. $\eqref{nutrientfield}$, we obtain
\begin{align}\label{nutbim1}
\frac{1}{2} \underline{\sigma}+
\int_{\Gamma_0} 
\left(
\underline{\sigma}
\frac{\partial G_{\mathbf{x}_0} }{\partial \mathbf{n}_0'}
- G_{\mathbf{x}_0}
\frac{\partial \sigma'}{\partial \mathbf{n}_0'}
\right) ds'
+\int_{\Gamma(t)} 
\left(
G_{\mathbf{x}_0}
\beta(1-\sigma')
-\frac{\partial G_{\mathbf{x}_0}}{\partial \mathbf{n}'}
\sigma'
\right)
 ds'
&=
\underline{\sigma}
&\forall& \mathbf{x}_0\in\Gamma_0,\\ \label{nutbim2}
\int_{\Gamma_0} 
\left(
\underline\sigma
\frac{\partial G_\mathbf{x} }{\partial \mathbf{n}_0'}
- G_\mathbf{x}
\frac{\partial \sigma'}{\partial \mathbf{n}_0'}
\right) ds'
+\int_{\Gamma(t)} 
\left(
G_\mathbf{x}
\beta(1-\sigma')
-\frac{\partial G_\mathbf{x} }{\partial \mathbf{n}'}
\sigma'
\right)
ds'
+\frac{1}{2}\sigma
&= \sigma
&\forall& \mathbf{x}\in\Gamma(t),
\end{align}
\normalsize
where we used standard jump relations of double layer potentials. 
Before taking the limit, alternatively, if we introduce linear operators by denoting single layer potentials by $\displaystyle\left.\mathcal{S}^{G_{\mathbf{z}}}\right|_{\mathcal C}[\phi]\equiv\int_{\mathcal C}G_{\mathbf{z}}\phi'ds'$ and double layer potentials by $\displaystyle\left.\mathcal{D}^{G_{\mathbf{z}}}_{\mathbf{n}_\star}\right|_{\mathcal C}[\phi]\equiv\int_{\mathcal C}\frac{\partial G_{\mathbf{z}}}{\partial \mathbf{n}_\star'}\phi' ds',$
 where $\mathcal C$ is the integral domain $\Gamma_0,\Gamma(t)$ or $\Gamma_0 \cup \Gamma(t)$ and $\mathbf{z}$ can be $\mathbf{x}_0,\mathbf{x}$ or $\mathbf{x}_*$, then Eq. $\eqref{greenint}$ can be rewritten as
\begin{equation}\label{greenintlayer}
\left.\mathcal{D}^{G_{\mathbf{x}_*}}_{\mathbf{n}_\star}\right|_{\Gamma_0\cup\Gamma(t)}\left[\sigma\right]
-\left.\mathcal{S}^{G_{\mathbf{x}_*}}\right|_{\Gamma_0\cup\Gamma(t)}\left[\frac{\partial \sigma}{\partial \mathbf{n}_\star}\right]
=-\sigma(\mathbf{x}_*)\quad\forall\mathbf{x}_*\in\Omega(t).
\end{equation}
The jump relations imply
\begin{equation}\label{jump}
\left\{
\begin{array}{ccc}
\displaystyle\lim_{h\rightarrow0^+}\left.\mathcal{D}^{G_{\mathbf{x}_0\pm h\mathbf{n}_{0}}}_{\mathbf{n}_0}\right|_{\Gamma_0}[\phi]&=\left(\left.\mathcal{D}^{G_{\mathbf{x}_0}}_{\mathbf{n_0}}\right|_{\Gamma_0}\pm\frac{1}{2}I\right)[\phi],&\forall\mathbf{x}_0\in\Gamma_0,\\
\displaystyle\lim_{h\rightarrow0^+}\left.\mathcal{D}^{G_{\mathbf{x}\pm h\mathbf{n}}}_{\mathbf{n}}\right|_{\Gamma(t)}[\phi]&=\left(\left.\mathcal{D}^{G_{\mathbf{x}}}_{\mathbf{n}}\right|_{\Gamma(t)}\pm\frac{1}{2}I\right)[\phi],&\forall\mathbf{x}\in\Gamma(t).
\end{array}
\right.
\end{equation}
Letting $\mathbf{x}_* \to \mathbf{x}_0 \in \Gamma_0, \mathbf{x}_* \to \mathbf{x} \in \Gamma(t)$ in Eq. $\eqref{greenintlayer}$ and using the $2^{nd}$ and the $3^{rd}$ equation in $\eqref{nutrientfield}$ with Eq. $\eqref{jump}$, we obtain
\begin{equation}\label{bimnut}
\left\{
\begin{array}{ccc}
\begin{aligned}
\left.\mathcal{S}^{G_{\mathbf{x}_0}}\right|_{\Gamma_0}\left[\frac{\partial\sigma}{\partial\mathbf{n}_0}\right]
+\left(\beta\left.\mathcal{S}^{G_{\mathbf{x}_0}}
+\mathcal{D}^{G_{\mathbf{x}_0}}_{\mathbf{n}}\right)\right|_{\Gamma(t)}\left[\sigma\right]
&=\underline\sigma\left(\left.\mathcal{D}^{G_{\mathbf{x}_0}}_{\mathbf{n}_0}-\frac{1}{2}I\right)\right|_{\Gamma_0}\left[1\right]
+\beta\left.\mathcal{S}^{G_{\mathbf{x}_0}}\right|_{\Gamma(t)}\left[1\right],\\
\left.\mathcal{S}^{G_{\mathbf{x}}}\right|_{\Gamma_0}\left[\frac{\partial\sigma}{\partial\mathbf{n}_0}\right]
+\left(\beta\left.\mathcal{S}^{G_{\mathbf{x}}}
+\mathcal{D}^{G_{\mathbf{x}}}_{\mathbf{n}}+\frac{1}{2}I\right)\right|_{\Gamma(t)}\left[\sigma\right]
&=\underline\sigma\left.\mathcal{D}^{G_{\mathbf{x}}}_{\mathbf{n}_0}\right|_{\Gamma_0}\left[1\right]
+\beta\left.\mathcal{S}^{G_{\mathbf{x}}}\right|_{\Gamma(t)}\left[1\right]
,
\end{aligned}
\end{array}
\right.
\end{equation}
\normalsize
which is the reformulation of Eqs. \eqref{nutbim1},  \eqref{nutbim2} into operator form (see \cite{veerapaneni2016integral} for similar treatment of expression).
This system needs to be solved for the unknowns $\left.\frac{\partial \sigma}{\partial \mathbf{n}_0}\right|_{\Gamma_0}$ and $\left.\sigma\right|_{\Gamma(t)}$.

\subsection{Direct BIM for the Laplace equation}
Consider the Green's function for Laplace equation in $\Omega(t)$:
\begin{equation}\label{LaplaceGreen}
\Delta \Phi_{\mathbf{x}_*} = -\delta_{\mathbf{x}_*}\quad\text{in }\Omega(t),
\end{equation}
where $\Phi_{\mathbf{x}_*}=\Phi_{\mathbf{x}_*}(\mathbf{x}')$, $\mathbf{x}_*,\mathbf{x}'\in \Omega(t)$ are the source point and field point respectively, and $\delta_{\mathbf{x}_*}(\mathbf{x}')$ is the Dirac delta function. 
The fundamental solution to Eq. $\eqref{LaplaceGreen}$ is 
\begin{equation}\label{fundamentalLaplace}
\Phi_{\mathbf{x}_*}(\mathbf{x}')=\frac{1}{2 \pi} \ln \frac{1}{r},
\end{equation}
where $r\equiv|{\mathbf{x}_*}-\mathbf{x}'|$.
Multiplying the $1^{st}$ equation in $\eqref{pressurefield}$ by $\Phi_{\mathbf{x}_*}$ and Eq. $\eqref{LaplaceGreen}$ by $-p$ and summing them up, we obtain
\begin{equation}\label{greenidp}
\Phi_{\mathbf{x}_*} \Delta p - p \Delta \Phi_{\mathbf{x}_*}
=p\delta_{\mathbf{x_*}}\quad\text{in }\Omega(t).
\end{equation}
Integrating Eq. $\eqref{greenidp}$ over $\Omega(t)$ and using Green's second identity, we have
\begin{equation}\label{greenintp}
\int_{\Gamma_0\cup\Gamma(t)} 
\left(\Phi_{\mathbf{x}_*} \frac{\partial p'}{\partial\mathbf{n}_\star'}
-p' \frac{\partial \Phi_{\mathbf{x}_*}}{\partial \mathbf{n}_\star'}\right) ds'
=p(\mathbf{x}_*)\quad\forall\mathbf{x}_*\in\Omega(t),
\end{equation}
where the symbol prime denotes the evaluation on field points, for example, $\mathbf{n}_\star'=\mathbf{n}_\star(\mathbf{x}')$ is the unit exterior (w.r.t. the domain $\Omega(t)$ enclosed by its boundaries) normal on the corresponding boundaries, i.e., $\mathbf{n}_\star'=-\mathbf{n}_0'$ on $\Gamma_0$, $\mathbf{n}_\star'= \mathbf{n}'$ on $\Gamma(t)$ and $\mathbf{n}_0'=\mathbf{n}_0(\mathbf{x}'),\mathbf{n}'=\mathbf{n}(\mathbf{x}')$ are the unit outer (pointing away from origin) normals on $\Gamma_0,\Gamma(t)$, respectively. 
Denote single layer potentials by $\displaystyle\left.\mathcal{S}^{\Phi_{\mathbf{z}}}\right|_{\mathcal C}[\phi]\equiv\int_{\mathcal C}\Phi_{\mathbf{z}}\phi' ds'$ and double layer potentials by $\displaystyle\left.\mathcal{D}^{\Phi_{\mathbf{z}}}_{\mathbf{n}_\star}\right|_{{\mathcal C}}[\phi]\equiv\int_{\mathcal C}\frac{\partial \Phi_{\mathbf{z}}}{\partial \mathbf{n}_\star'}\phi' ds'$, where $\mathcal C$ is the integral domain $\Gamma_0,\Gamma(t)$ or $\Gamma_0 \cup \Gamma(t)$ and $\mathbf{z}$ can be $\mathbf{x}_0,\mathbf{x}$ or $\mathbf{x}_*$. 
Then Eq. $\eqref{greenintp}$ can be rewritten as
\begin{equation}\label{greenintlayerp}
\left.\mathcal{S}^{\Phi_{\mathbf{x}_*}}\right|_{\Gamma_0\cup\Gamma(t)}\left[\frac{\partial p}{\partial \mathbf{n}_\star}\right]-
\left.\mathcal{D}^{\Phi_{\mathbf{x}_*}}_{\mathbf{n}_\star}\right|_{\Gamma_0\cup\Gamma(t)}\left[p\right]=p(\mathbf{x}_*)\quad\forall\mathbf{x}_*\in\Omega(t).
\end{equation}
The jump relations imply
\begin{equation}\label{jumpp}
\left\{\begin{array}{ccc}\displaystyle\lim_{h\rightarrow0^+}\left.\mathcal{D}^{\Phi_{\mathbf{x}_0\pm h\mathbf{n}_{0}}}_{\mathbf{n}_0}\right|_{\Gamma_0}[\phi]&=\left(\left.\mathcal{D}^{\Phi_{\mathbf{x}_0}}_{\mathbf{n_0}}\right|_{\Gamma_0}\pm\frac{1}{2}I\right)[\phi]&\forall\mathbf{x}_0\in\Gamma_0,\\\displaystyle\lim_{h\rightarrow0^+}\left.\mathcal{D}^{\Phi_{\mathbf{x}\pm h\mathbf{n}}}_{\mathbf{n}}\right|_{\Gamma(t)}[\phi]&=\left(\left.\mathcal{D}^{\Phi_{\mathbf{x}}}_{\mathbf{n}}\right|_{\Gamma(t)}\pm\frac{1}{2}I\right)[\phi]&\forall\mathbf{x}\in\Gamma(t).\end{array}\right.
\end{equation}
Letting $\mathbf{x}_* \to \mathbf{x}_0 \in \Gamma_0, \mathbf{x}_* \to \mathbf{x} \in \Gamma(t)$ in Eq. $\eqref{greenintlayerp}$ and using the $2^{nd}$ and the $3^{rd}$ equation in $\eqref{pressurefield}$ with Eq. $\eqref{jumpp}$, we obtain
\begin{equation}\label{bimpressure}
\left\{
\begin{array}{ccc}
\begin{aligned}
\left(-\frac{1}{2}I+\left.\mathcal{D}^{\Phi_{\mathbf{x}_0}}_{\mathbf{n}_0}\right)\right|_{\Gamma_0}\left[p\right]
+\left.\mathcal{S}^{\Phi_{\mathbf{x}_0}}\right|_{\Gamma(t)}\left[\frac{\partial p}{\partial\mathbf{n}}\right]
&=
\left.\mathcal{S}^{\Phi_{\mathbf{x}_0}}\right|_{\Gamma_0}\left[\frac{\partial p}{\partial\mathbf{n}_0}\right]
+\left.\mathcal{D}^{\Phi_{\mathbf{x}_0}}_{\mathbf{n}}\right|_{\Gamma(t)}\left[p\right]
,\\
\left.\mathcal{D}^{\Phi_{\mathbf{x}}}_{\mathbf{n}_0}\right|_{\Gamma_0}\left[p\right]
+\left.\mathcal{S}^{\Phi_{\mathbf{x}}}\right|_{\Gamma(t)}\left[\frac{\partial p}{\partial\mathbf{n}}\right]
&=
\left.\mathcal{S}^{\Phi_{\mathbf{x}}}\right|_{\Gamma_0}\left[\frac{\partial p}{\partial\mathbf{n}_0}\right]
+\left(\left.\mathcal{D}^{\Phi_{\mathbf{x}}}_{\mathbf{n}}+\frac{1}{2}I\right)\right|_{\Gamma(t)}\left[p\right]
.
\end{aligned}
\end{array}
\right.
\end{equation}
\normalsize
This system needs to be solved for the unknowns $\left.p\right|_{\Gamma_0}$ and $\left.\frac{\partial p}{\partial \mathbf{n}}\right|_{\Gamma(t)}$.


\subsection{Summary of BIM linear systems}
By rewriting Eqs. \eqref{bimnut} and \eqref{bimpressure}, we summarize this section as the following linear systems:
\begin{itemize}
    \item \textbf{\textit{Nutrient field.}}
        \begin{equation}\label{eq:linsysnut}
        \begin{pmatrix}
            \left.\mathcal{S}^{G_{\mathbf{x}_0}}\right|_{\Gamma_0}
            &\left(\beta\left.\mathcal{S}^{G_{\mathbf{x}_0}}
+\mathcal{D}^{G_{\mathbf{x}_0}}_{\mathbf{n}}\right)\right|_{\Gamma(t)}\\
            \left.\mathcal{S}^{G_{\mathbf{x}}}\right|_{\Gamma_0}
            &\left(\beta\left.\mathcal{S}^{G_{\mathbf{x}}}
+\mathcal{D}^{G_{\mathbf{x}}}_{\mathbf{n}}+\frac{1}{2}I\right)\right|_{\Gamma(t)}
        \end{pmatrix}
        \underbrace{
        \begin{pmatrix}
            \left.\frac{\partial\sigma}{\partial\mathbf{n}_0}\right|_{\Gamma_0}
            \\
            \left.\sigma\right|_{\Gamma(t)}
        \end{pmatrix}}_\text{unknown}
        =
        \begin{pmatrix}
            \underline\sigma\left(\left.\mathcal{D}^{G_{\mathbf{x}_0}}_{\mathbf{n}_0}-\frac{1}{2}I\right)\right|_{\Gamma_0}\left[1\right]
+\beta\left.\mathcal{S}^{G_{\mathbf{x}_0}}\right|_{\Gamma(t)}\left[1\right]
        \\
        \underline\sigma\left.\mathcal{D}^{G_{\mathbf{x}}}_{\mathbf{n}_0}\right|_{\Gamma_0}\left[1\right]
+\beta\left.\mathcal{S}^{G_{\mathbf{x}}}\right|_{\Gamma(t)}\left[1\right]
        \end{pmatrix}
        \end{equation}
    \item \textbf{\textit{Pressure field.}}
        \begin{equation}\label{eq:linsysps}
        \begin{pmatrix}
            \left(-\frac{1}{2}I+\left.\mathcal{D}^{\Phi_{\mathbf{x}_0}}_{\mathbf{n}_0}\right)\right|_{\Gamma_0}
            &\left.\mathcal{S}^{\Phi_{\mathbf{x}_0}}\right|_{\Gamma(t)}\\
            \left.\mathcal{D}^{\Phi_{\mathbf{x}}}_{\mathbf{n}_0}\right|_{\Gamma_0}
            &\left.\mathcal{S}^{\Phi_{\mathbf{x}}}\right|_{\Gamma(t)}
        \end{pmatrix}
        \underbrace{
        \begin{pmatrix}
            \left.p\right|_{\Gamma_0}
            \\
            \left.\frac{\partial p}{\partial\mathbf{n}}\right|_{\Gamma(t)}
        \end{pmatrix}}_\text{unknown}
        =
        \begin{pmatrix}
            \left.\mathcal{S}^{\Phi_{\mathbf{x}_0}}\right|_{\Gamma_0}\left[\frac{\partial p}{\partial\mathbf{n}_0}\right]
+\left.\mathcal{D}^{\Phi_{\mathbf{x}_0}}_{\mathbf{n}}\right|_{\Gamma(t)}\left[p\right]
        \\\left.\mathcal{S}^{\Phi_{\mathbf{x}}}\right|_{\Gamma_0}\left[\frac{\partial p}{\partial\mathbf{n}_0}\right]
+\left(\left.\mathcal{D}^{\Phi_{\mathbf{x}}}_{\mathbf{n}}+\frac{1}{2}I\right)\right|_{\Gamma(t)}\left[p\right]
        \end{pmatrix}
        \end{equation}
\end{itemize}
Note that the right hand side of Eq. \eqref{eq:linsysps} is known from the boundary conditions in Eq. $\eqref{pressurefield}$, which are reformulated from Eqs. \eqref{psbc1} and \eqref{psbc2}.
We then outline our numerical method in the following algorithm:
\begin{algorithm}[H]
  \caption{Numerical algorithm}\label{Algorithm}
  \begin{algorithmic}[1]
    \Procedure{Evolve tumor interface $\Gamma(t)$}{}
      \State Initialization at $t=0$\Comment{Given initial shape $\Gamma(0)$ and obtain equal arclength meshes}
      \For{$t=0$\texttt{ to }$t_\text{final}$}
        \State \texttt{solve BIM linear systems in Eqs. \eqref{eq:linsysnut}, \eqref{eq:linsysps} by GMRES}
        \State \texttt{compute normal velocity $V$ with Eq. \eqref{normal velocity}}
        \State \texttt{compute tangent velocity $T$ with Eq. \eqref{tangent velocity}}
        \State \texttt{update tangent angle $\theta$ via Eq. \eqref{linear propagator method}}
        \State \texttt{update arclength $s_\alpha$ via Eq. \eqref{linear propagator method 2}}
        \State \texttt{update tumor interface $\Gamma(t)$}
      \EndFor
    \EndProcedure
  \end{algorithmic}
\end{algorithm}

For details on how these integrals are discretized and how the systems Eq. \eqref{bimnut} and Eq. \eqref{bimpressure} are solved, we refer the reader to \ref{Appendix:Evaluation of BIM}.
For details on the corresponding discretizations used and the time-stepping method for evolving the interface $\Gamma(t)$ in time, we refer the reader to \ref{Appendix:Evolution of the interface} and \ref{Appendix:Time stepping}.

\section{Linear Analysis}\label{sec:LA}
The purpose of performing linear analysis is twofold: (1) to study the morphologically unstable regime of parameters by analyzing the linear solutions; (2) to validate our numerical method by checking the agreement of the nonlinear simulation results with the linear solutions at early times. In this section, we focus on the first aspect and we will study the second aspect in Section \ref{sec:comparison non vs lin}. 
We present the results of a linear stability analysis (details are provided in \ref{appendix:linear analysis}) of the non-dimensional sharp interface equations \eqref{pressurefield}--\eqref{normal velocity} reformulated in the preceding section. The linear stability of perturbed radially symmetric tumors was previously analyzed in \cite{cristini2003nonlinear}. Here, we extend their results to take into account the nutrient field with Robin boundary condition and the fixed necrotic core $\Gamma_0$ with radius $R_0$.
Consider a $l^{th}$ mode perturbation of a radially symmetric tumor interface $\Gamma$:
\begin{equation}
r(\theta,t)=R(t)+\delta(t) \cos l \theta,
\end{equation} 
where $r$ is the tumor/host interface, $R$ is the radius of the underlying circle, $\delta$ is the dimensionless perturbation size, and $\theta$ is
the polar angle. We first deduce that on the necrotic boundary, the pressure and the nutrient flux are given by:
\small
\begin{align}\notag
\left.p\right|_{\Gamma_0}=&\ \mathcal{P}\left(A_{1}\left(I_{0}(R)+I_{1}\left(R_{0}\right) R_{0} \ln\left(\frac{R_0}{R}\right)\right)+A_{2}\left(K_{0}(R)-K_{1}\left(R_{0}\right) R_{0} \ln \left(\frac{R_0}{R}\right)\right)\right)\\
&-\chi_{\sigma}\left(A_{1} I_{0}(R)+A_{2} K_{0}(R)\right) 
-\frac{\mathcal{P} \mathcal{A}}{2}\left(R_{0}^{2} \ln \left(\frac{R_0}{R}\right)+\frac{R^{2}}{2}\right)+\frac{\mathcal{G}^{-1}}{R}\notag\\
&+\delta e^{i l \theta}\left(\mathcal{P}\left(\left(A_{1} I_{1}(R)-A_{2} K_{1}(R)+B_{1} I_{l}(R)+B_{2} K_{l}(R)-\left(A_{1} I_{1}\left(R_{0}\right)-A_{2} K_{1}\left(R_{0}\right)\right) \frac{R_{0}}{R}\right) \frac{2(R R_0)^{l}}{R^{2l}+R_{0}^{2 l}}\right.\right. \notag\\
&+\left.\frac{R_{0}}{l}\left(B_{1}I_{l-1}\left(R_{0}\right)-B_{2}K_{l-1}\left(R_{0}\right)\right)\frac{R_{0}^{2l}-R^{2 l}}{R^{2 l}+R_{0}^{2 l}}\right)\notag\\
&-\chi_\sigma\left(A_{1} I_{1}(R)-A_{2} K_{1}(R)+B_{1} I_{l}(R)+B_{2} K_{l}(R)-\left(A_{1} I_{1}\left(R_{0}\right)-A_{2} K_{1}\left(R_{0}\right)\right) \frac{R_{0}}{R}\right) \frac{R_{0}^{2l}-R^{2 l}}{R^{2 l}+R_{0}^{2 l}}
\notag\\
&+\left.\left(\frac{\mathcal{P} \mathcal{A} R_{0}^{2}}{2 R}-\frac{\mathcal{P} \mathcal{A} R}{2}+\mathcal{G}^{-1} \frac{l^{2}-1}{R^{2}}\right) \frac{2(R R_0)^{l}}{R^{2 l}+R_{0}^{2 l}}\right),\notag\\
\left.\frac{\partial \sigma}{\partial n}\right|_{\Gamma_0} 
=&A_{1} I_{1}\left( R_0\right)-A_2K_1(R_0)+\delta e^{i l \theta}\left(B_{1}\left(I_{l-1}\left(R_0\right)-\frac{l}{R_0} I_{l}\left( R_0\right)\right)-B_2\left(K_{l-1}\left(R_0\right)+\frac{l}{R_0} K_{l}\left( R_0\right)\right)\right) .
\end{align}
\normalsize

And on the tumor boundary, the nutrient concentration and the pressure flux are given by:
\small
$$
\begin{aligned}
\left.\sigma\right|_{\Gamma}
=&\ A_{1} I_{0}\left(R\right)+A_2K_0(R)+\delta e^{i l \theta}\left(A_{1} I_{1}\left(R\right)-A_2K_1(R)+B_{1} I_{l}(R)+B_2 K_l(R)\right) ,\notag\\
\left.\frac{\partial p}{\partial n}\right|_{\Gamma}=&\left(\mathcal{P}\left(A_{1} I_{1}\left(R_{0}\right)-A_{2} K_{1}\left(R_{0}\right)\right)-\frac{\mathcal{P} \mathcal{A} R_{0}}{2}\right) \frac{R_{0}}{R}\notag\\
&+\delta e^{i l \theta}\left({\mathcal{P}}\left(-{\left(A_{1} I_{1}\left(R_{0}\right)-A_{2} K_{1}\left(R_{0}\right)\right)} \frac{R_{0}}{R^{2}}\right.\right.\notag\\
&+\left(A_{1} I_{1}(R)-A_{2} K_{1}(R)+B_{1} I_{l}(R)+B_{2} K_{l}(R)-\left(A_{1} I_{1}\left(R_{0}\right)-A_{2} K_{1}\left(R_{0}\right)\right) \frac{R_{0}}{R}\right) \frac{R^{2 l}-R_{0}^{2 l}}{R^{2 l}+R_{0}^{2 l}} \frac{l}{R}\notag\\
&+\left.2\left(B_{1} I_{l-1}\left(R_{0}\right)-B_{2} K_{l-1}\left(R_{0}\right)\right) \frac{R^{l} R_{0}^{l}}{R^{2 l}+R_{0}^{2l} } \frac{R_{0}}{R}\right)\notag\\
&-\chi_{\sigma}\left(A_{1} I_{1}(R)-A_{2} K_{1}(R)+B_{1} I_{l}(R)+B_{2} K_{l}(R)-\left(A_{1} I_{1}\left(R_{0}\right)-A_{2} K_{1}\left(R_{0}\right)\right) \frac{R_{0}}{R}\right) \frac{R^{2 l}-R_{0}^{2 l}}{R^{2 l}+R_{0}^{2 l}} \frac{l}{R}\notag\\
&\left.+\frac{\mathcal{P} \mathcal{A}}{2}\left(\frac{R_{0}}{R}\right)^{2}
+\left(\frac{\mathcal{P} \mathcal{A} R_{0}^{2}}{2 R}-\frac{\mathcal{P} \mathcal{A} R}{2}+\mathcal{G}^{-1} \frac{l^{2}-1}{R^{2}}\right) \frac{R^{2 l}-R_{0}^{2 l}}{R^{2 l}+R_{0}^{2 l}} \frac{l}{R}\right).\notag
\end{aligned}
$$
\normalsize
Then we deduce that the evolution equation for the tumor radius $R$ is given by:
\begin{equation}\label{Rdot}
\frac{dR}{dt}=\underbrace{
\mathcal{P} 
\left(A_1 I_1(R)-A_2 K_1(R)
-\frac{R_0 }{R}\left(A_1 I_1(R_0)-A_2 K_1(R_0)\right)\right)}_\text{Proliferation}
-
\underbrace{\frac{\mathcal{P}\mathcal{A}}{2}
\frac{R^2-R_0^2
}{R}}_\text{Apoptosis}
,
\end{equation}

where
\begin{eqnarray}{}
A_1&=&
\frac{\underline{\sigma}\left(K_{1}(R)-\beta K_{0}(R)\right)+\beta K_{0}\left(R_{0}\right)}{K_{0}\left(R_{0}\right)\left(\beta I_{0}(R)+I_{1}(R)\right)+I_{0}\left(R_{0}\right)\left(K_{1}(R)-\beta K_{0}(R)\right)},\label{A1}\\
A_2&=&
\frac{\underline{\sigma}\left(\beta I_{0}(R)+I_{1}(R)\right)-\beta I_{0}\left(R_{0}\right)}{K_{0}\left(R_{0}\right)\left(\beta I_{0}(R)+I_{1}(R)\right)+I_{0}\left(R_{0}\right)\left(K_{1}(R)-\beta K_{0}(R)\right)}\label{A2}
.
\end{eqnarray}

The equation of the shape perturbation $\frac{\delta}{R}$ is given by:
\small
\begin{eqnarray}\label{ddovRdt}
&&\left(\frac{\delta}{R}\right)^{-1}
\frac{d}{dt}{\left(\frac{\delta}{R}\right)}\notag\\
&=&\overbrace{\mathcal{P} \mathcal{A}\left(\left(1-\left(\frac{R_{0}}{R}\right)^{2}\right)\left(1-\frac{2 R_{0}^{2 l}}{R^{2 l}+R_{0}^{2 l}}\right) \frac{l}{2}-\left(\frac{R_{0}}{R}\right)^{2}\right)}^\textbf{Apoptosis}-\overbrace{\mathcal{G}^{-1} \frac{l\left(l^{2}-1\right)}{R^{3}}\left(1-\frac{2 R_{0}^{2 l}}{R^{2 l}+R_{0}^{2 l}}\right)}^\textbf{Cell-cell adhesion} \notag\\ 
&&
-\overbrace{\mathcal{P}{\beta}\left(\frac{1}{R}+A_{1}\left(I_{1}(R)-\frac{I_{0}(R)}{R}\right)-A_{2}\left(K_{1}(R)+\frac{K_{0}(R)}{R}\right)+B_{1} I_{l}(R)+B_{2} K_{l}(R)\right)}^\textbf{Angiogenesis}\notag\\
&&+\overbrace{\chi_{\sigma}\left(A_{1} I_{1}(R)-A_{2} K_{1}(R)+B_{1} I_{l}(R)+B_{2} K_{l}(R)-\left(A_{1} I_{1}\left(R_{0}\right)-A_{2} K_{1}\left(R_{0}\right)\right) \frac{R_{0}}{R}\right) \left(1-\frac{2 R_{0}^{2 l}}{R^{2 l}+R_{0}^{2l}}\right) \frac{l}{R}
}^\textbf{Chemotaxis}\notag\\
&&+\overbrace{\mathcal{P}\left(\left(A_{1} I_{1}\left(R_{0}\right)-A_{2} K_{1}\left(R_{0}\right)\right) \frac{R_{0}}{R^{2}}\left(2+l\left(1-\frac{2 R_{0}^{2 l}}{R^{2 l}+R_{0}^{2l}}\right)\right)-2\left(B_{1} I_{l-1}\left(R_{0}\right)-B_{2} K_{l-1}\left(R_{0}\right)\right) \frac{R^{l} R_{0}^{l}}{R^{2 l}+R_{0}^{2 l}} \frac{R_{0}}{R}\right)
}^\textbf{Proliferation}\notag\\
&&-\overbrace{\mathcal{P}\left(\left(A_{1} I_{1}(R)-A_{2} K_{1}(R)+B_{1} I_{l}(R)+B_{2} K_{l}(R) \right)\left(1-\frac{2 R_{0}^{2 l}}{R^{2 l}+R_{0}^{2l}}\right) \frac{l}{R}\right)
}^\textbf{Proliferation},
\end{eqnarray}
\normalsize
where $I_{l}(R)$ and $K_{l}(R)$ are modified Bessel functions of the first and of the second kind, respectively. A complete derivation and the expressions of $A_i, B_i,i=1,2$, which contain the parameters $R_{0}$, $\underline{\sigma}$ and $\beta$ are given in \ref{appendix:linear analysis}.

In Fig. \ref{fig:compareCristini2003}, we first show that the linear stability analysis results by Cristini et al. on the evolution of both the tumor radius $R$ (Eq. (7) in \cite{cristini2003nonlinear}) and the shape factor $\frac{\delta}{R}$ (Eq. (11) in \cite{cristini2003nonlinear} but with typo, corrected formula can be found as Eq. (4.7a) in \cite{cristini2010multiscale}), which are plotted in dash lines, can be recovered from our results as $R_0\rightarrow 0$ and $\beta\rightarrow\infty$ in Eq. \eqref{Rdot} and Eq. \eqref{ddovRdt}, which is verified numerically by decreasing $R_0$ from $0.1$ (left) to $\epsilon_{machine}\approx 2.2204\times 10^{-16}$ (right) and by increasing $\beta$ from $0.5$ to $100$ as plotted in solid lines. 

\begin{figure}
    \begin{minipage}{1.0\linewidth}
        \centering
        \includegraphics[width=\textwidth]{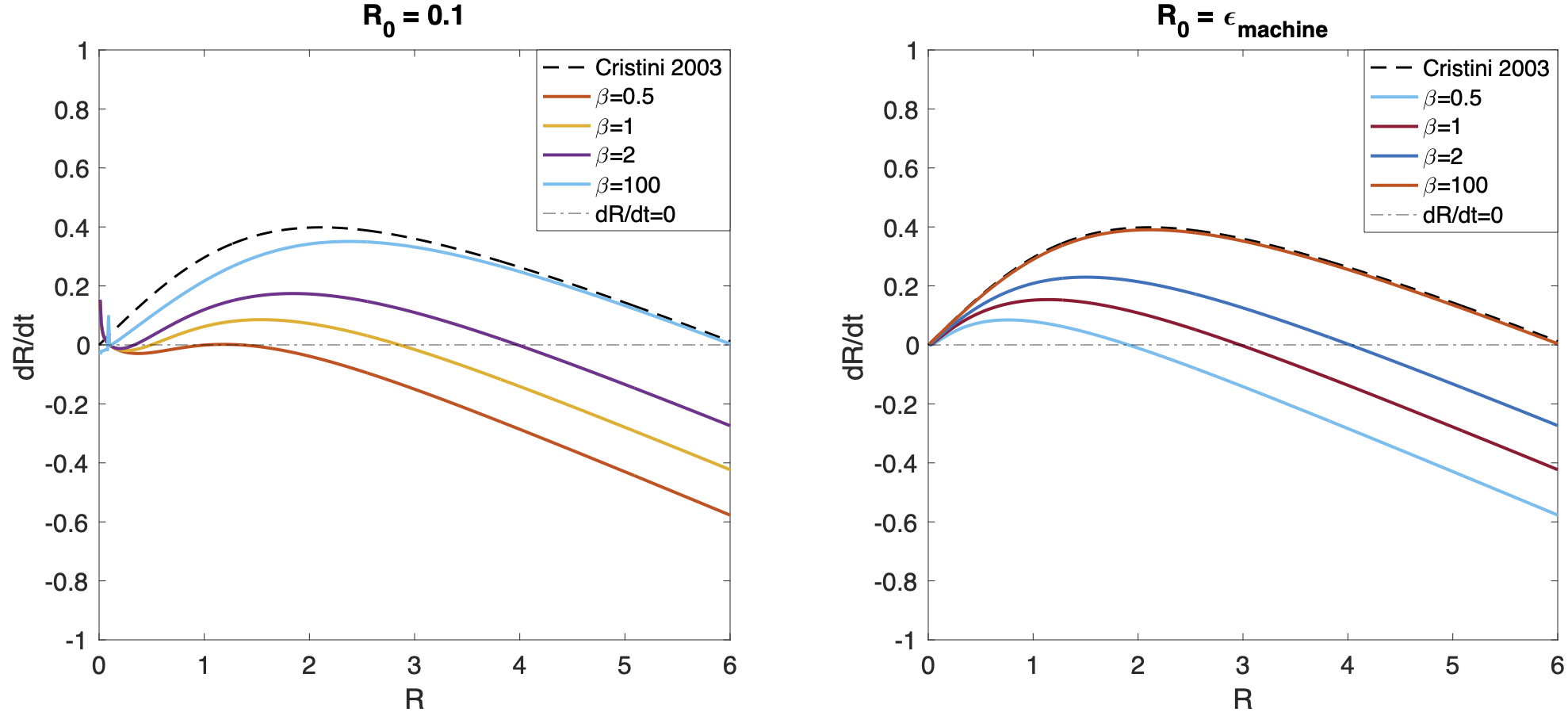}[a]
    \end{minipage}
    \begin{minipage}{1.0\linewidth}
        \centering
        \includegraphics[width=\textwidth]{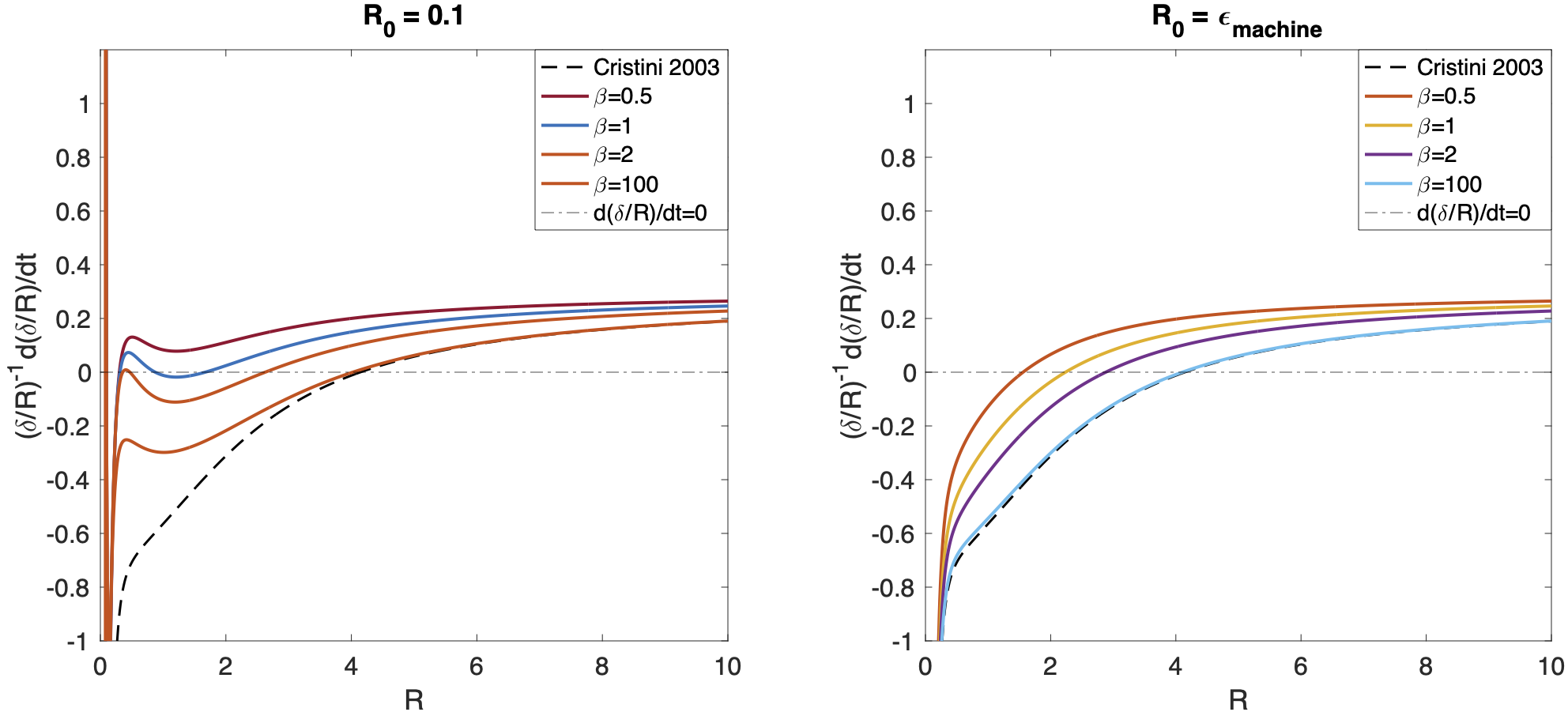}[b]
    \end{minipage}
    \caption{The rate of change of the tumor radius $R$ in [a] and the shape factor $\frac{\delta}{R}$ in [b] for the case with parameters $\mathcal{P}=1,\mathcal{A}=0.3,\chi_\sigma=0,\underline{\sigma}=0,\mathcal{G}^{-1}=0.001$. In both [a] and [b], the necrotic boundary is circular with radius $R_0=0.1$ (left) and $R_0=\epsilon_{machine}\approx 2.2204\times 10^{-16}$ (right). The solid lines are plotted with Eq. \eqref{Rdot} in [a] and Eq. \eqref{ddovRdt} in [b] for $\beta=0.5,1,2,100$. The dash lines are plotted with Eq. (7) in 
    \cite{cristini2003nonlinear} in [a] and Eq. (4.7a) in \cite{cristini2010multiscale} in [b].}
    \label{fig:compareCristini2003}
\end{figure}

We next characterize the stability regime by determining $\mathcal{A}=\mathcal{A}_c$ as a function of the unperturbed radius $R$ such that $\frac{d}{dt}{\left(\frac{\delta}{R}\right)}=0$:
\small
\begin{eqnarray}
\mathcal{A}_{c}&=& \left( \overbrace{ \mathcal{G}^{-1} \frac{l\left(l^{2}-1\right)}{\mathcal{P} R^{3}}\left(1-\frac{2 R_{0}^{2 l}}{R^{2 l}+R_{0}^{2 l}}\right)}^\textbf{Cell-cell adhesion}\right.\notag\\
&&+\overbrace { \beta\left(\frac{1}{R}+A_{1}\left(I_{1}(R)-\frac{I_{0}(R)}{R}\right)-A_{2}\left(K_{1}(R)+\frac{K_{0}(R)}{R}\right)+B_{1} I_{l}(R)+B_{2} K_{l}(R)\right)}^{\textbf{Angiogenesis}}\notag\\
&&-\overbrace{\frac{\chi_{\sigma}}{\mathcal{P}}\left(A_{1} I_{1}(R)-A_{2} K_{1}(R)+B_{1} I_{l}(R)+B_{2} K_{l}(R)-\left(A_{1} I_{1}\left(R_{0}\right)-A_{2} K_{1}\left(R_{0}\right)\right) \frac{R_{0}}{R}\right) \left(1-\frac{2 R_{0}^{2 l}}{R^{2 l}+R_{0}^{2l}}\right) \frac{l}{R}
}^\textbf{Chemotaxis to Proliferation}\notag\\
&&-
\left(A_{1} I_{1}(R_{0})-A_{2} K_{1}(R_{0}\right) \frac{R_{0}}{R^{2}}
\left(2+l\left(1-\frac{2R_{0}^{2l}}{R^{2l}+R_{0}^{2l}}\right)\right)\notag\\
&&+\left(A_{1} I_{1}(R)-A_{2} K_{1}(R)+B_{1} I_{l}(R)+B_{2} K_{l}(R)\right)\left(1-\frac{2 R_{0}^{2 l}}{R^{2 l}+R_{0}^{2 l}}\right) \frac{l}{R}\notag\\
&&\left.+2\left(B_{1} I_{l-1}(R_{0})-B_{2} K_{l-1}(R_{0})\right) 
\frac{R^{l} R_{0}^{l}}{R^{2 l}+R_{0}^{2 l}} \frac{R_{0}}{R}\right)\notag\\
&&/\left(\left(1-\left(\frac{R_{0}}{R}\right)^{2}\right)\left(1-\frac{2 R_{0}^{2 l}}{R^{2 l}+R_{0}^{2 l}}\right) \frac{l}{2}-\left(\frac{R_{0}}{R}\right)^{2}\right).
\nonumber
\label{apoptc}
\end{eqnarray}

\normalsize
Here, $\mathcal{A}_c$ is the critical value of apoptosis that divides regimes of stable growth ($\mathcal{A} < \mathcal{A}_c$, \textit{e.g.}, below the curve) and regimes of unstable growth ($\mathcal{A} > \mathcal{A}_c$, \textit{e.g.}, above the curve) for a given mode $l=2$ and necrotic radius $R_0=0.1$. We focus on the parameters $\mathcal{P}=5,\mathcal{G}^{-1}=0.001,\underline{\sigma}=0.2$.

\begin{figure}
\centering
\includegraphics[width=\textwidth]{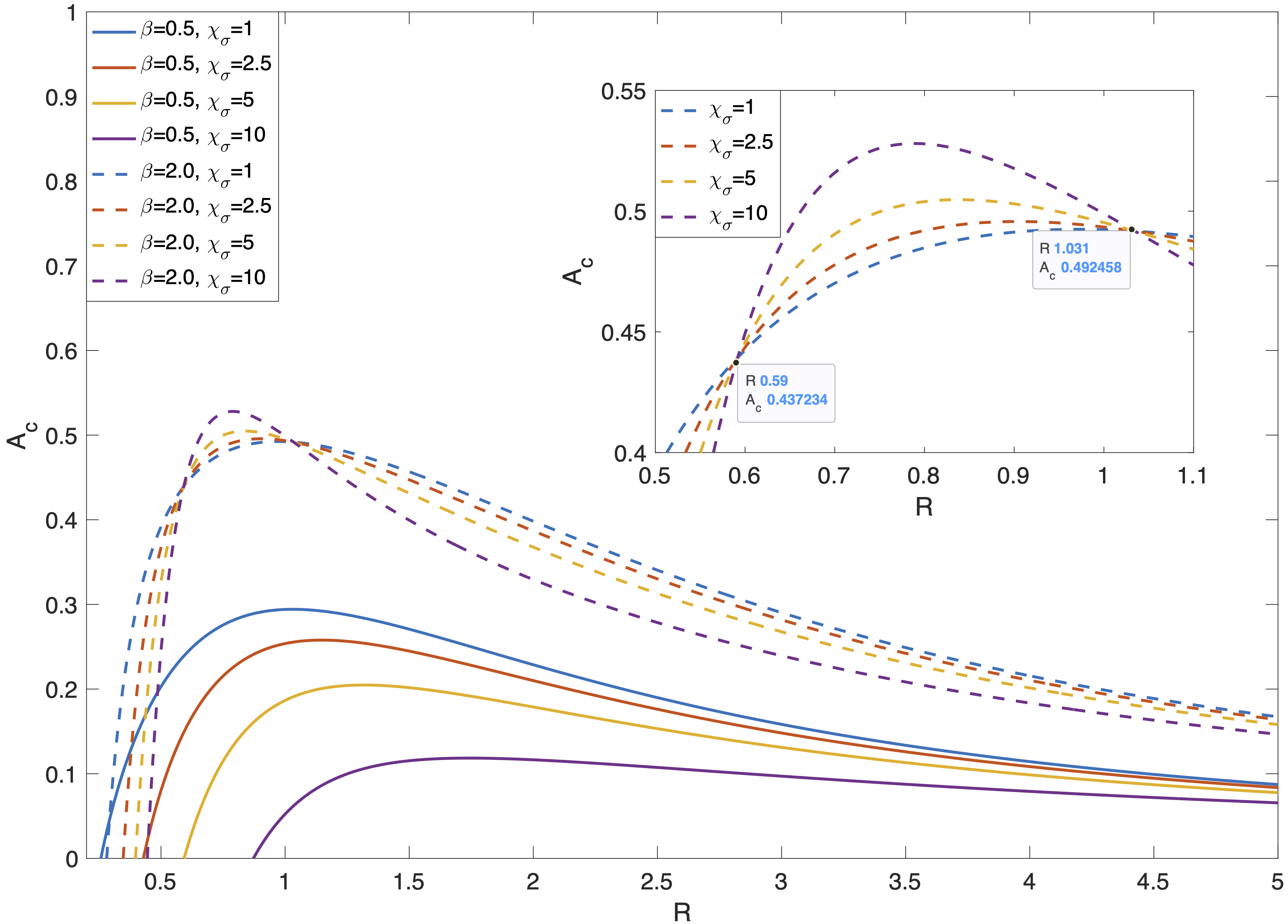}
\caption{Critical apoptosis parameter $\mathcal{A}_c$ as a function of unperturbed radius $R$ from equation \eqref{apoptc}, circular necrotic boundary $R_0= 0.1$, and $\chi_\sigma$ labeled in the legend. Solid: $\beta = 0.5$; Dashed: $\beta = 2$. See text for details.}\label{fig:apoptc}
\end{figure}

In Fig. \ref{fig:apoptc}, we plot $\mathcal{A}_c$ as a function of $R$, $\beta=0.5$ (solid), $\beta=2$ (dashed), and $\chi_\sigma$ as labeled in the legend. The figure reveals that the unstable regime expands in general with stronger taxis, with the inset as an exception that shows an opposite tendency for $\beta=2$ when the tumor radius $R$ is small and between 0.59 to 1.03. Moreover, all the dashed curves are pulled upward under a richer supply of nutrients from vasculature $(\beta=2)$. These two observations suggest that as the tumor is growing, \textit{chemotaxis} may enhance the morphological instability while \textit{angiogenesis} might inhibit it, which will be further investigated in Section \ref{sec:nonlinear}.

\section{Results}\label{sec:Result}
\subsection{Numerical Convergence in time and space}

\label{sec:4.1}
In this section, we test the convergence of our method.

\begin{figure}
\centering
\begin{minipage}{0.85\linewidth}
\centering
\includegraphics[width=\textwidth]{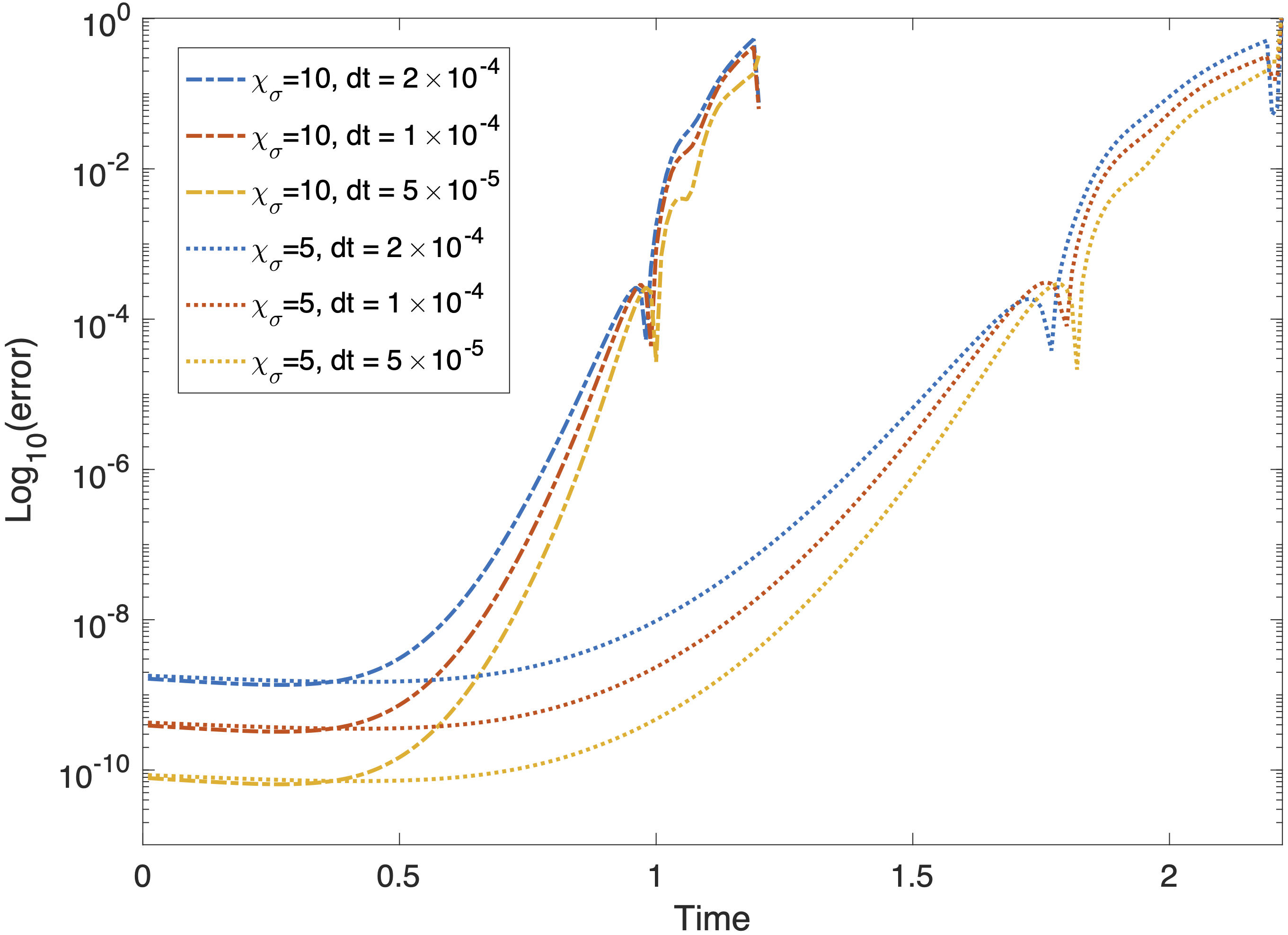}[a]
\end{minipage}
\begin{minipage}{0.85\linewidth}
\centering
\includegraphics[width=\textwidth]{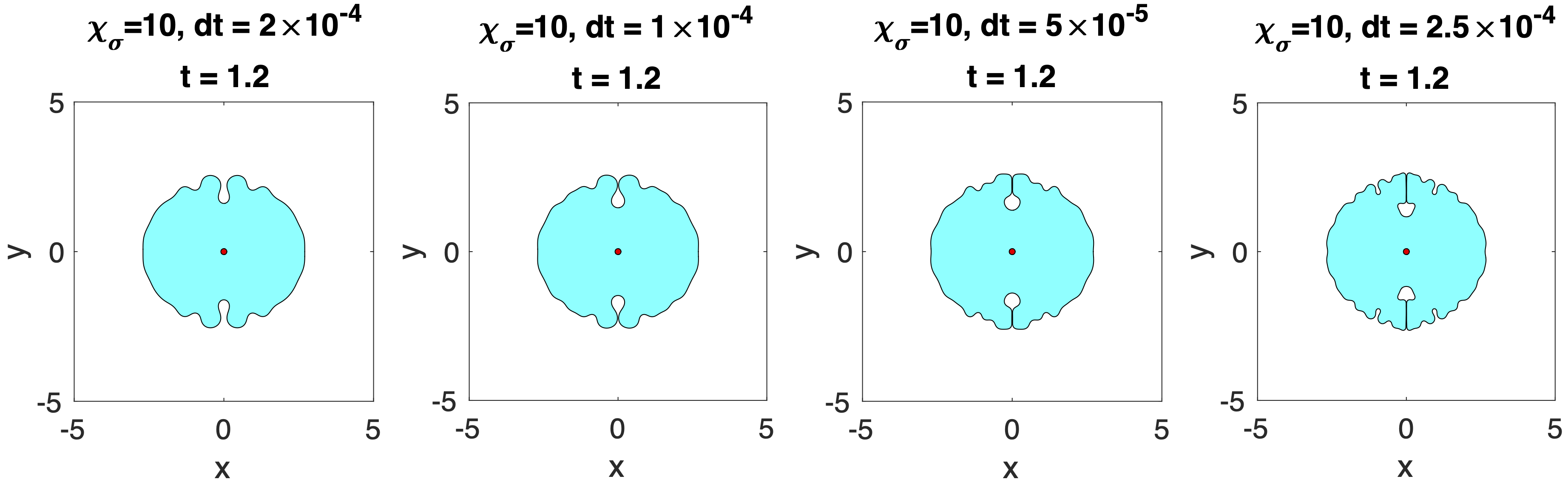}[b]
\end{minipage}
\begin{minipage}{0.85\linewidth}
\centering
\includegraphics[width=\textwidth]{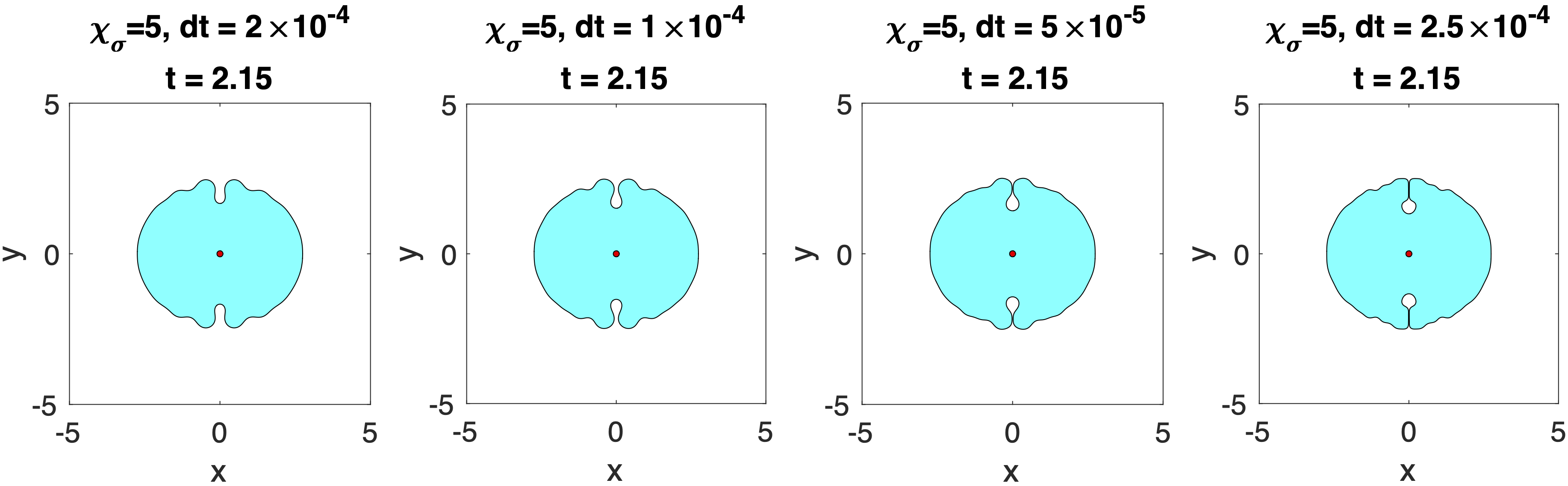}[c]
\end{minipage}

\caption{Temporal resolution studies for the case with parameters $\mathcal P=5$, $\mathcal A=0.25$, $\chi_\sigma=5$ $\text{ (dash-dot lines)}$, $\chi_\sigma=10\text{ (dot lines)}$, $\beta=0.5$,  $\underline{\sigma}=0.2$ and $\mathcal{G}^{-1}=0.001$. The necrotic boundary is circular with radius $R_0=0.1$, and initial tumor boundary is $r=2.5+ 0.1\cos (2 \theta)$. The errors shown are calculated as 
the differences of tumor area
between the solution 
with $\Delta t=2.5\times 10^{-5}$ and those with $\Delta t=2\times10^{-4},1\times10^{-4},5\times10^{-5}$. Dash-dot lines are the errors for $\chi_{\sigma}=10$ and dot lines are the errors for $\chi_{\sigma}=5$ in [a]. The corresponding tumor morphologies before the computations terminate are shown in [b] with $\chi_\sigma=10$ and in [c] with $\chi_\sigma=5$. In all cases $N=512$.
 }
\label{fig:time}
\end{figure}

First, we present a temporal resolution study in Fig. \ref{fig:time} [a]. 
The errors calculated by differences of tumor area 
between the simulation with $\Delta t=2.5\times 10^{-5}$ and $\Delta t=2\times10^{-4},1\times10^{-4},5\times10^{-5}$ , respectively,  are plotted versus time. In all cases, the number of spatial collocation points is $N = 512$. We first examine the case plotted in dash-dot lines with $\chi_{\sigma}=10$, the necrotic radius $R_0=0.1$ and other parameters $P=5$, $\mathcal A=0.25$, $\beta=0.5$, $\underline{\sigma}=0.2$ and $\mathcal{G}^{-1}=0.001$. 
We find a factor of 4 is observed when $\Delta t$ is halved, indicating a second order convergence rate. This is expected since the time stepping scheme is second order accurate (\ref{Appendix:Time stepping}.)
Similar observations can also be found for the case plotted in dot lines with $\chi_{\sigma}=5$. 
The corresponding tumor morphologies are shown in Fig. \ref{fig:time} [b] with $\chi_\sigma=10$ and [c] with $\chi_\sigma=5$. 

\begin{figure}
    \centering
    \includegraphics[width=\textwidth]{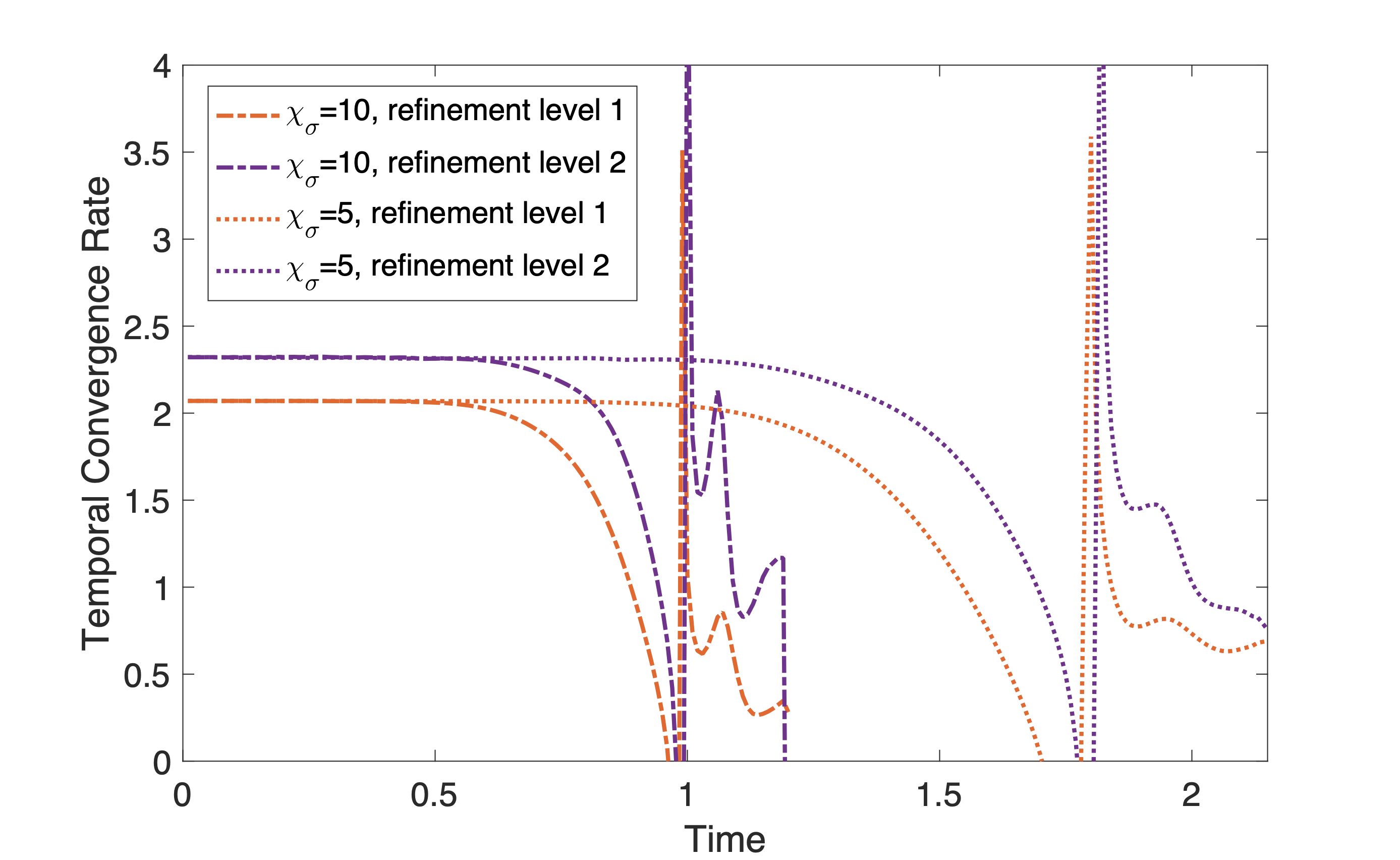}
    \caption{The temporal convergence rate $C_{n}=\frac{\ln (e_n/e_{n+1})}{\ln 2}$, where $e_n=|A_n-A_0|$ denotes the error of the $n-$th level, $A_0$ is the area of the finest case ($\Delta t = 2.5\times 10^{-5}$), the first refinement level $A_1$ is that of the roughest case ($\Delta t=2\times 10^{-4}$) and $A_{n+1}$ is finer than $A_{n}$. Refinement level 1 uses $e_1,e_2$ and refinment level 2 uses $e_2,e_3$. Dot lines correspond to $\chi_\sigma=10$ and dash-dot lines correspond to $\chi_\sigma=5$.}
    \label{fig:CR}
\end{figure}

To be more precise, define temporal convergence rate $C_{n}=\frac{\ln (e_n/e_{n+1})}{\ln 2}$, where $e_n=|A_n-A_0|$ denotes the error of the $n-$th level, $A_0$ is the area of the finest case ($\Delta t = 2.5\times 10^{-5}$), the first refinement level $A_1$ is that of the roughest case ($\Delta t=2\times 10^{-4}$) and $A_{n+1}$ is finer than $A_{n}$. Thus, what we have plotted in Fig. \ref{fig:time} [a] are what we denoted as $e_i,i=1,2,3$ here for the two cases $\chi_\sigma=5,10$. In Fig. \ref{fig:CR}, the temporal convergence rate of refinement level 1 (using $e_1,e_2$) and 2 (using $e_2,e_3$) are plotted versus time for the two cases $\chi_\sigma=5,10$. We see the convergence rate is of second order, even better for higher refinement level, 
before $t\approx 0.8$ for $\chi_\sigma=10$ and $t\approx 1.5$ for $\chi_\sigma\approx 1.5$. 
We remark that at $t\approx 1\ (\chi_\sigma =10)$ and $t\approx 1.8\  (\chi_\sigma=5)$ the convergence rate deteriorates when there are regions of the tumor boundary starting to touch one another (which can be found in Fig. \ref{fig:time} [b] and [c] for the case with finest $\Delta t$). This is caused by the approaching of the topological singularity of the tumor boundary, which we will discuss later.

\begin{figure}
\centering
\begin{minipage}{\linewidth}
\centering
\includegraphics[width=0.85\textwidth]{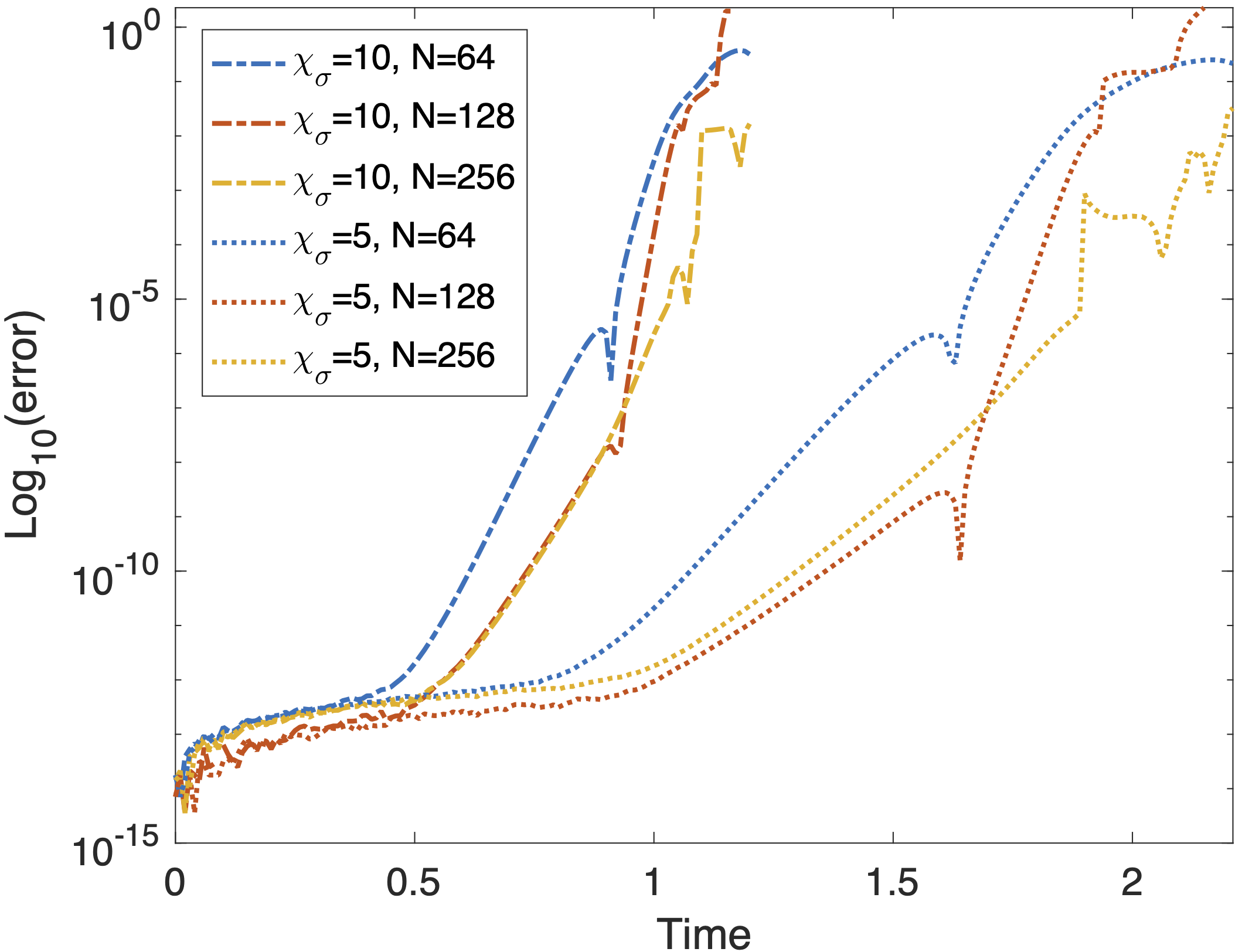}[a]
\end{minipage}
\begin{minipage}{\linewidth}
\centering
\includegraphics[width=0.85\textwidth]{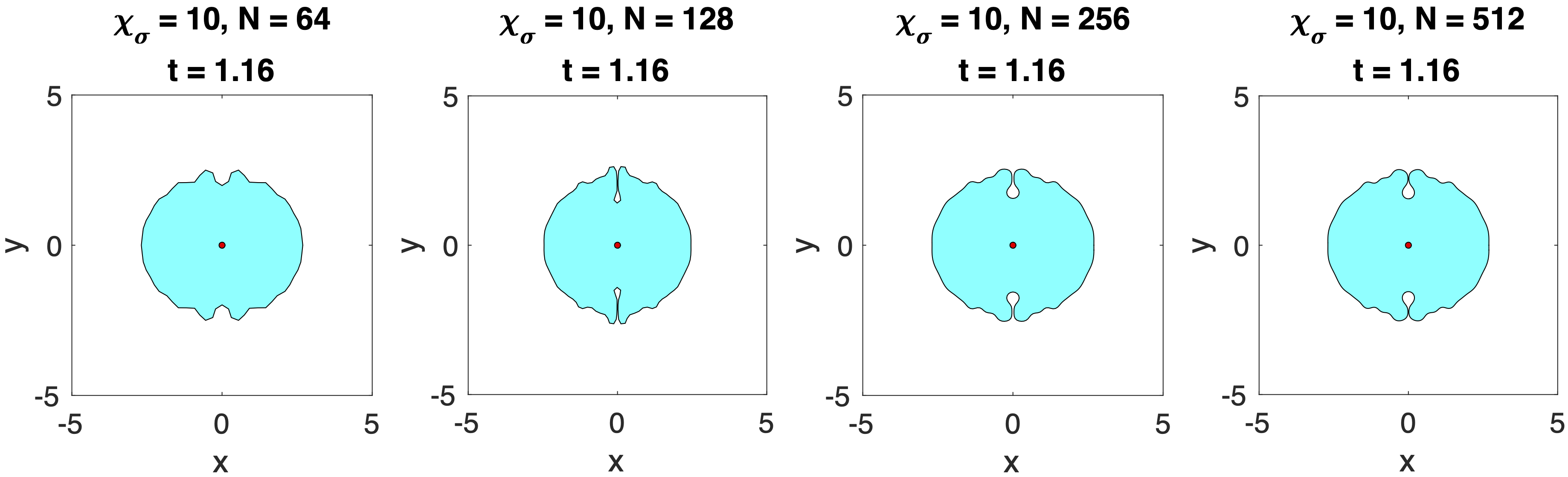}[b]
\end{minipage}
\begin{minipage}{\linewidth}
\centering
\includegraphics[width=0.85\textwidth]{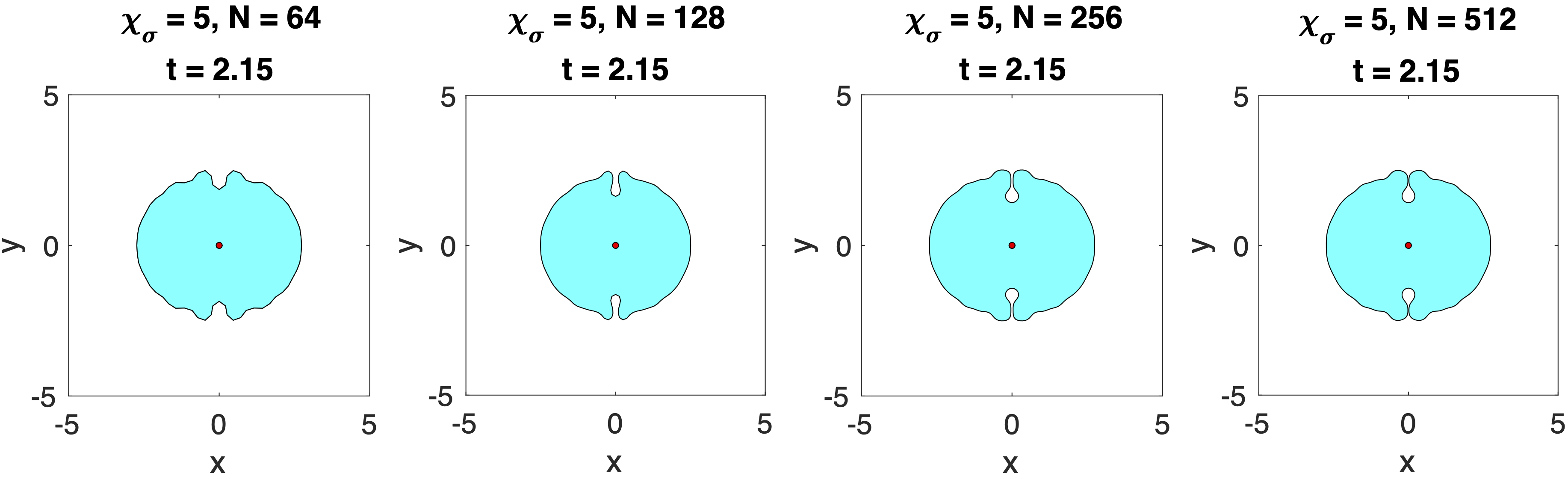}[c]
\end{minipage}

\caption{Spatial resolution studies for the case with parameters $\mathcal P=5$, $\mathcal A=0.25$, $\chi_\sigma=5$ $\text{ (dash-dot lines)}$, $\chi_\sigma=10\text{ (dot lines)}$, $\beta=0.5$,  $\underline{\sigma}=0.2$ and $\mathcal{G}^{-1}=0.001$. The necrotic boundary is circular with radius $R_0=0.1$, and initial tumor boundary is $r=2.5+ 0.1\cos (2 \theta)$. The errors are calculated as the 
differences of tumor area between the solution with $N=512$
and those with $N=64, 128, 256$.
The corresponding tumor morphologies are shown in [b] with $\chi_\sigma=10$ and in [c] with $\chi_\sigma=5$. In all cases $\Delta t=5\times 10^{-5}$.
 }
\label{fig:space}
\end{figure}

In space, the accuracy of our simulation is established by a resolution study of the simulation shown in Fig. \ref{fig:space}, with all biophysical parameters the same as the case in Fig. \ref{fig:time}. The spatial error is investigated by varying the number $N$ of spatial collocation points representing the tumor boundary $\Gamma(t)$. The errors calculated by 
tumor area of the solution between the simulation with $N=512$ and those with $N=64,128,256$ respectively are plotted versus time in Fig. \ref{fig:space} [a]. In all cases, the time step is $\Delta t=5\times 10^{-5}.$ At early times, the error is dominated by the tolerance for solving the integral equations $(1\times 10^{-10}).$ This is consistent with the spectral accuracy of our method. 
We remark that such error control lasts longer by refining the time step size $\Delta t$.
The corresponding tumor morphologies are shown in Fig. \ref{fig:space} [b] with $\chi_\sigma=10$ and [c] with $\chi_\sigma=5$.

\begin{figure}
\centering
\begin{minipage}{\linewidth}
\centering
\includegraphics[width=0.85\textwidth]{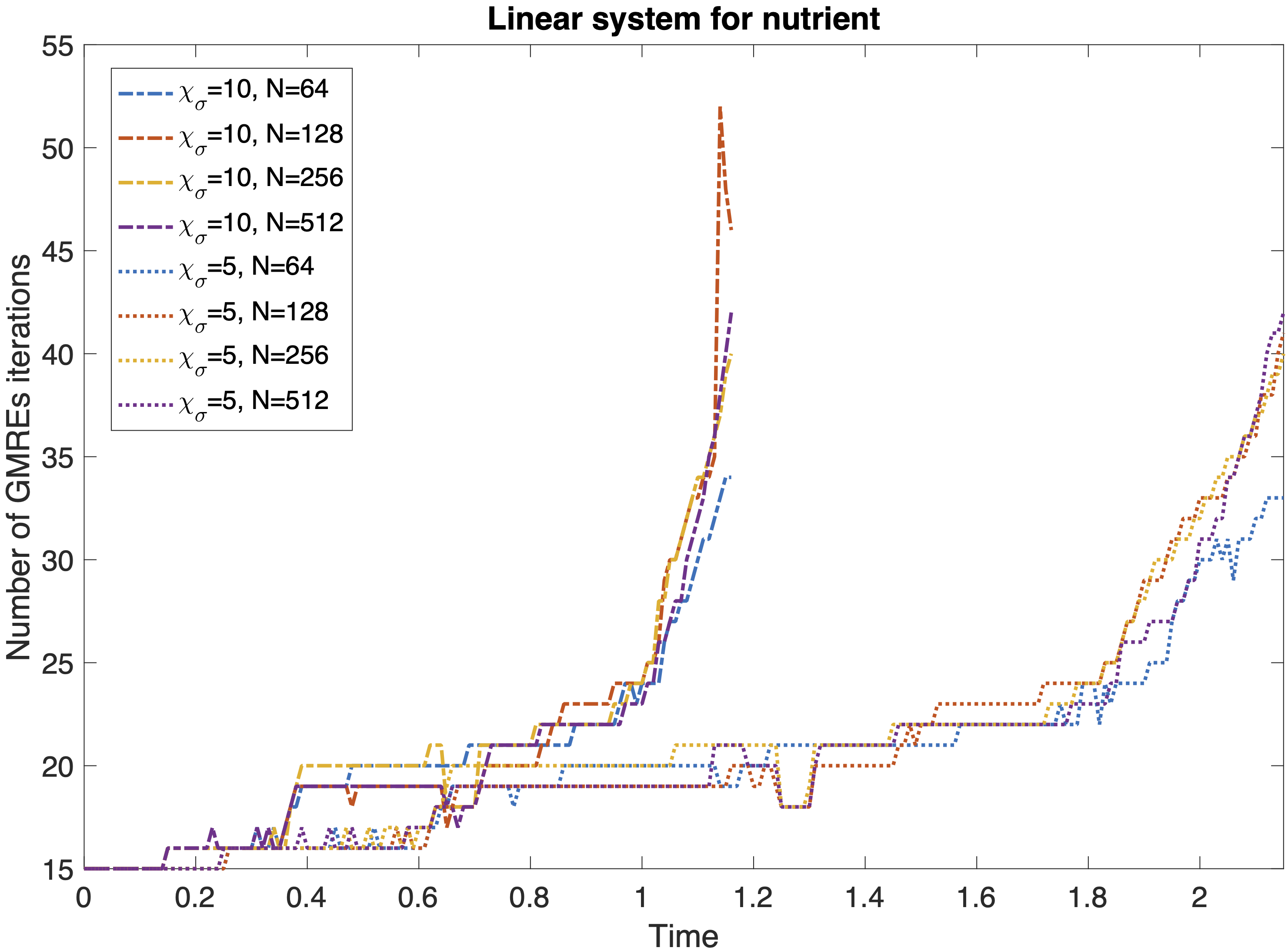}[a]
\end{minipage}
\begin{minipage}{\linewidth}
\centering
\includegraphics[width=0.85\textwidth]{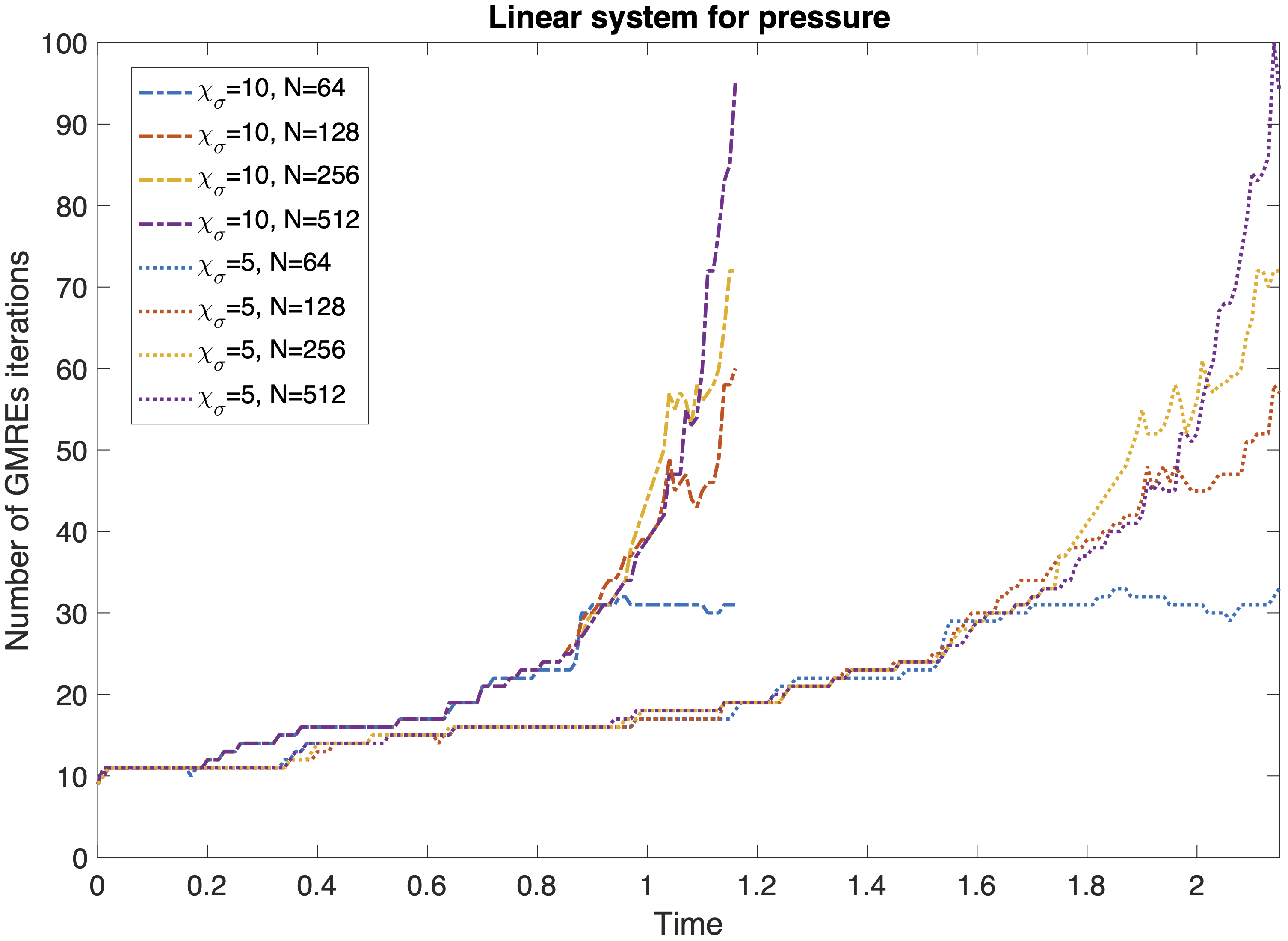}[b]
\end{minipage}

\caption{GMRES iteration numbers of linear system for nutrient in [a] and for pressure in [b] for the case with parameters $\mathcal P=5$, $\mathcal A=0.25$, $\chi_\sigma=5$ $\text{ (dash-dot lines)}$, $\chi_\sigma=10\text{ (dot lines)}$,  $\beta=0.5$,  $\underline{\sigma}=0.2$ and $\mathcal{G}^{-1}=0.001$. The necrotic boundary is circular with radius $R_0=0.1$, and initial tumor boundary is $r=2.5+ 0.1\cos (2 \theta)$. In all cases $\Delta t=5\times 10^{-5}$.}
\label{fig:gmres}
\end{figure}

In Fig. \ref{fig:gmres}, we show the GMRES iteration number of the linear system in Eq. \eqref{bimnut} for nutrient in [a] and of the one in Eq.  \eqref{bimpressure} for pressure in [b] with all the parameters the same as those used in Fig. \ref{fig:space}. We observe that the iteration number increases as time evolves.
An abrupt increase was observed before the computation stops. To understand this phenomena, we need to mention that as commented in \cite{jou1997} for the application of  boundary integral methods to elastic media problems in the end of their \textit{Sec. 3.1}, if the precipitates merge, which corresponds to that neighboring tumor fingers are getting really close (almost merging) to one another as shown in Fig. $\ref{fig:space}$ [b] and [c], then a topological change occurs and the kernels of the boundary integral systems become singular. In computations, this difficulty is reflected by a rapid increase in condition number and GMRES iteration count when the boundary become very close. Thus, it is related to the geometry of $\Gamma(t)$ and the coefficient matrix of the boundary integral system will become ill-conditioned when any two columns of it (calculated by the distance-dependent Green functions as described in $\textit{Sec. 4}$) become nearly the same due to the geometry change. Although this difficulty exists, we observe that the numerical approach still works for our purpose to study the evolution of the tumor and all the biophysical quantities required on both interfaces before the topological change.

\subsection{Comparison with linear analysis when the necrotic core is circular}\label{sec:comparison non vs lin}

\begin{figure}
\begin{minipage}{1.0\linewidth}
\centering
\includegraphics[width=\textwidth]{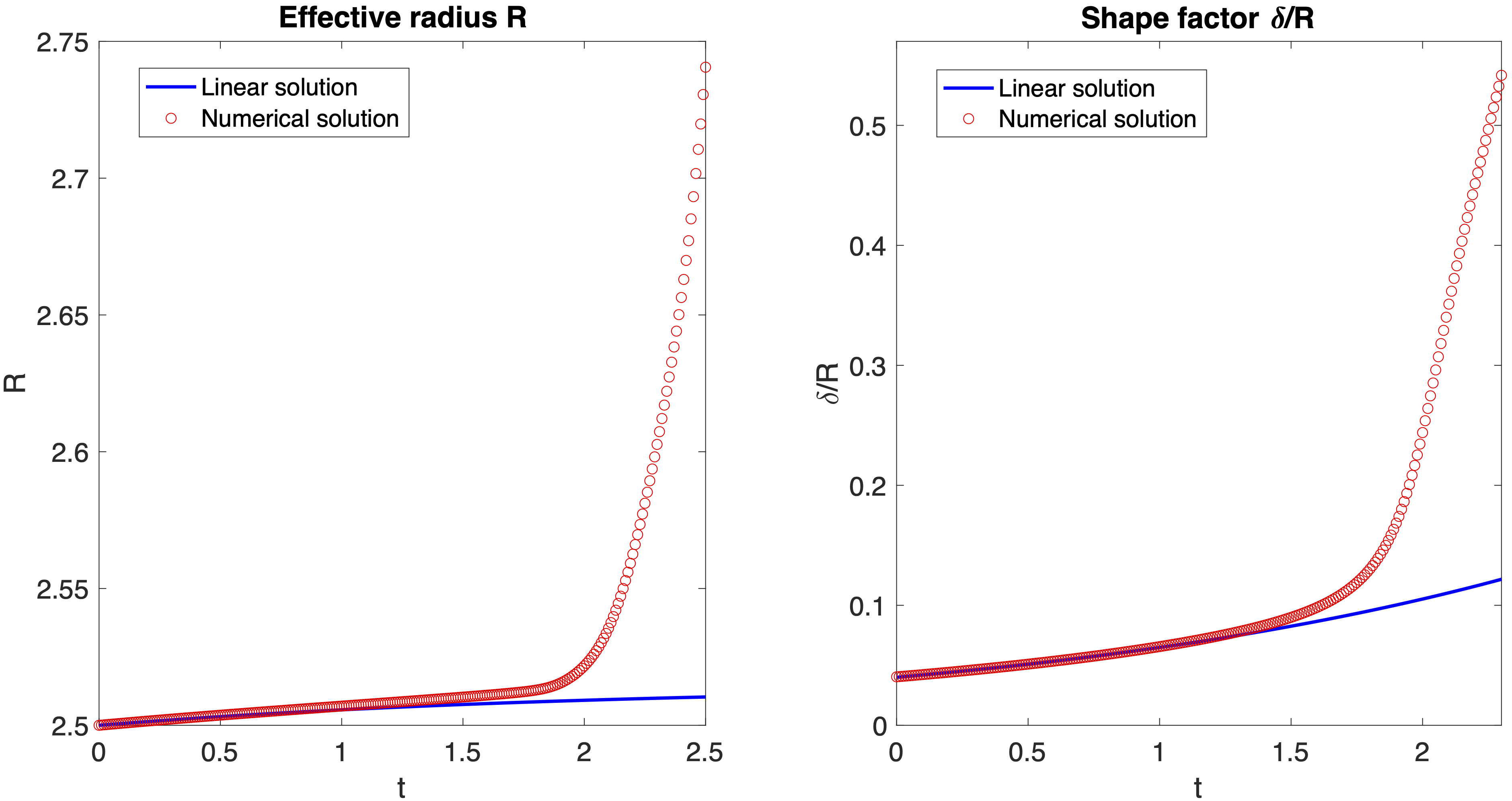}[a]
\end{minipage}
\\
\begin{minipage}{1.0\linewidth}
\centering
\includegraphics[width=\textwidth]{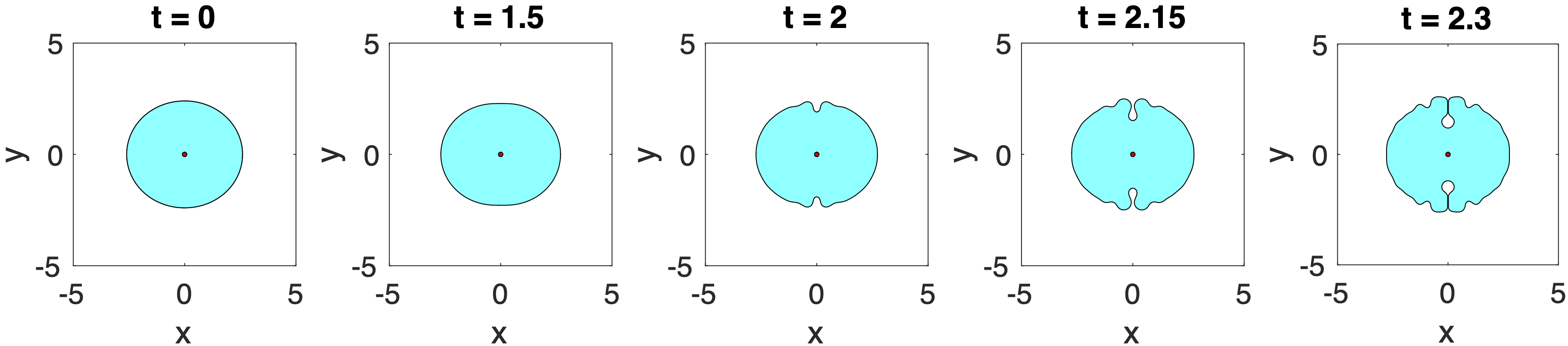}[b]
\end{minipage}
\caption{In [a]: A comparison between linear theory (blue curves) and the nonlinear simulations (red circles) for the effective radius $R$ (left) and shape factor $\frac{\delta}{R}$ (right). Red circles: nonlinear simulations; Blue lines: linear solutions. In [b]: The nonlinear tumor morphologies. The parameters are $\mathcal P=5$, $\mathcal{A}=0.25$, $\beta=0.5$, $\chi_{\sigma}=5$, $\underline{\sigma}=0.2$ and $\mathcal{G}^{-1}=0.001$. The necrotic boundary is circular with $R_0=0.1$, and initial tumor boundary is $r=2.5+ 0.1\cos (2 \theta)$. Here, $N=512$ and $\Delta t=1\times 10^{-4}$.}
\label{fig:comparison}
\end{figure}

We next compare the nonlinear simulation with linear theory. The results are shown in Fig. \ref{fig:comparison} [a] where we consider the case with the parameters $\mathcal P=5$, $\mathcal{A}=0.25$, $\beta=0.5$, $\chi_{\sigma}=5$, $\underline{\sigma}=0.2$ and $\mathcal{G}^{-1}=0.001$, circular necrotic boundary $R_0=0.1$, and initial tumor boundary $r=2.5+ 0.1\cos (2 \theta)$. While there is good agreement between the linear and nonlinear results at early times, both the effective radius and shape factors are under-predicted by linear theory at later times. The nonlinear tumors, shown in Fig. \ref{fig:comparison} [b], show the development of the tumor tissue encapsulating the surrounding tissue. In the next section, we will use the nonlinear simulation to investigate those factors that influence tumor progression.

\subsection{Nonlinear simulation with a circular necrotic core}
\label{sec:nonlinear}

In this section we study the factors of tumor growth with a fixed circular necrotic core from three different aspects: \textit{angiogenesis}, \textit{chemotaxis} and \textit{necrosis}. In the first aspect, we take $\beta=0.5,1,2$ to show the effect of \textit{angiogenesis}. In the second, we take $R_0=0.1,1,1.5$ to study the effect of \textit{necrosis}. In the third, we take $\chi_{\sigma}=5,10$ to demonstrate the effect of \textit{chemotaxis}. 

\paragraph{\textbf{Angiogenesis}}
\begin{figure}
\centering
\includegraphics[width=\textwidth]{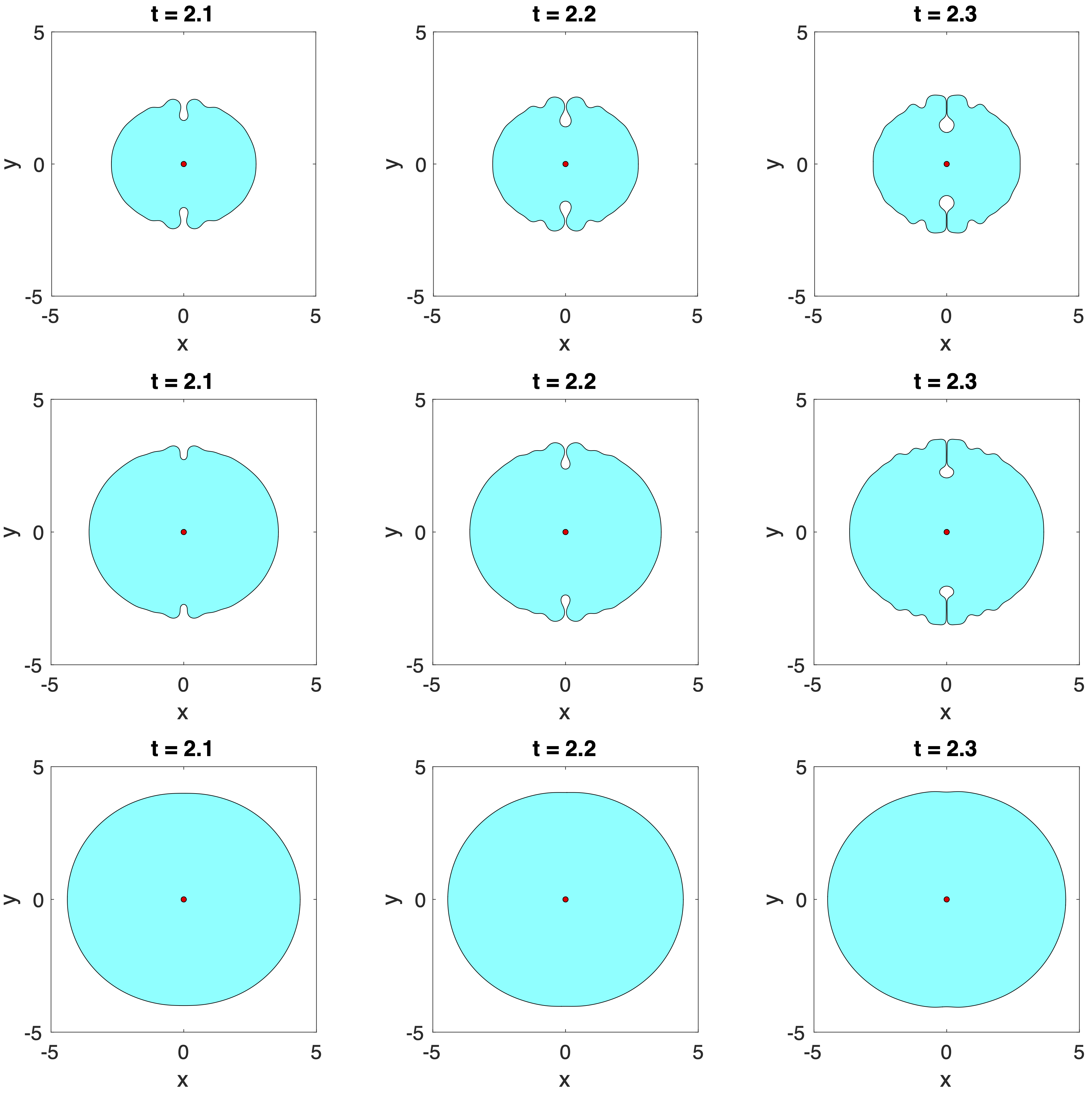}
\caption{The tumor morphologies under different values of the angiogenesis factor $\beta=0.5$ (first row), $\beta=1$ (second row) and $\beta=2$ (third row). The remaining parameters are $\mathcal P=5, \mathcal{A}=0.25, \chi_{\sigma}=5, R_0=0.1,\underline{\sigma}=0.2$ and $\mathcal{G}^{-1}=0.001$. The initial tumor boundary is $r=2.5+0.1\cos(2\alpha)$. Here, $N=512$ and $\Delta t=1\times 10^{-4}$.}\label{fig:angiogenesis}
\end{figure}

We present the evolution of tumor in Fig. \ref{fig:angiogenesis} where the angiogenesis factor $\beta=0.5,1,2$ corresponds to the rows from top to bottom, and the columns from left to right correspond to different time $t=2.1,2.2,2.3$. Hence, by comparing the evolution among the rows, we can see the effect of \textit{angiogenesis}. 

At $t=2.1$ (first column), for example, as $\beta$ increases, a larger tumor size is observed, which indicates that tumor vascularization will enhance the growth rate of the tumor. For $\beta=2$ (third row), $t=2.3$, we can see the tumor eventually evolves into a compact spheroid 
which suggests that $\textit{angiogenesis}$ may inhibit the instability and is consistent with our linear stability results.
We also remark that in our case with finite \textit{angiogenesis} factor $\beta$ (Robin boundary condition), the tumor interface has a much more unstable shape than the case with ``infinite'' $\beta$ (Dirichlet boundary condition) in \cite{cristini2003nonlinear}.

\paragraph{\textbf{Necrosis and Chemotaxis}}

\begin{figure}
\includegraphics[width=\textwidth]{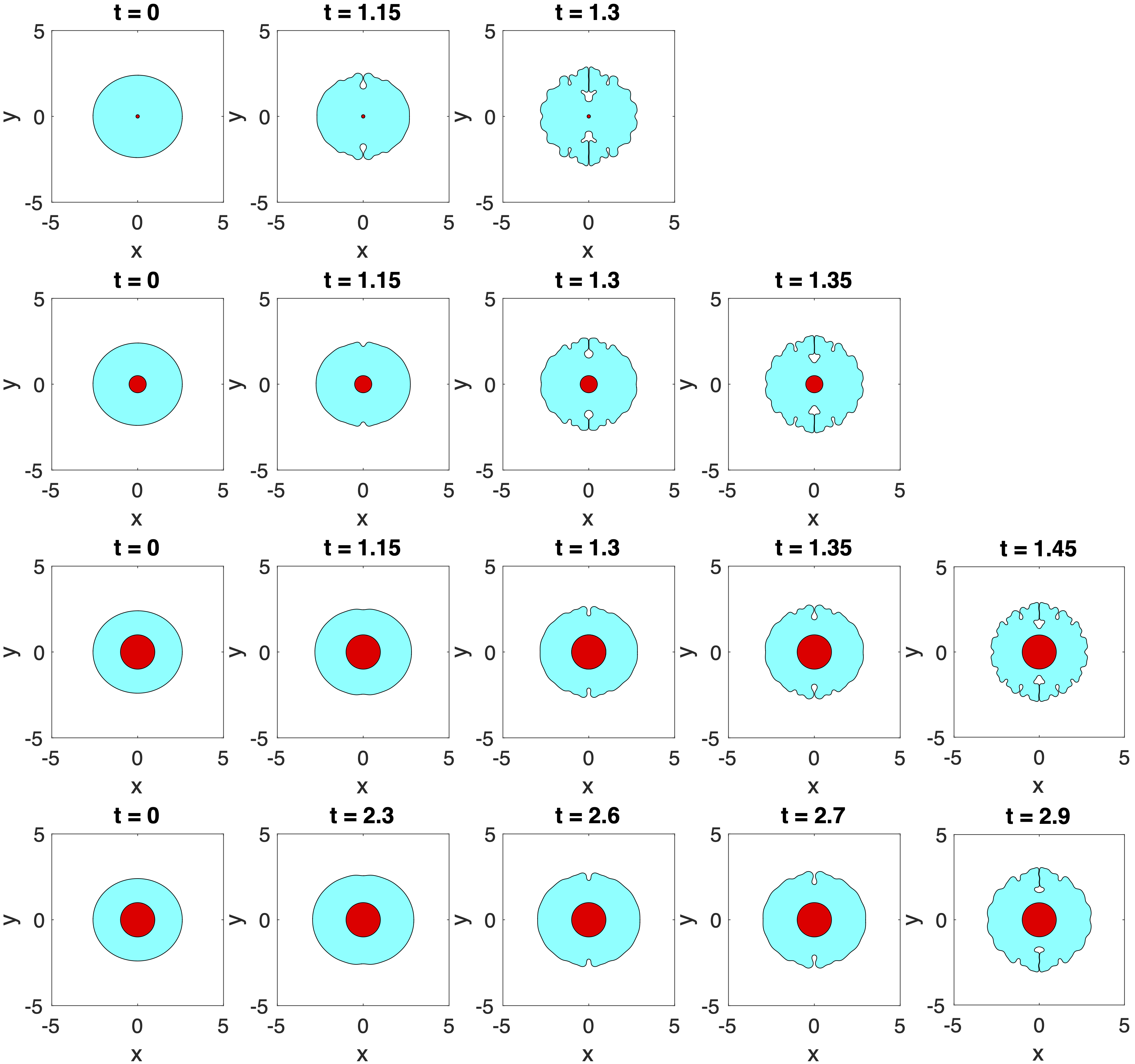}
\caption{The tumor morphologies under different sizes of the necrotic core $R_0=0.1$ (first row), $0.5$ (second row) and $1$ (third and fourth row). The remaining parameters are $\mathcal P=5, \mathcal{A}=0.25, \beta=0.5, \chi_\sigma=10$ (first 3 rows), $5$ (fourth row) $,\underline{\sigma}=0.2$ and $\mathcal{G}^{-1}=0.001$. The initial tumor boundary is $r=2.5+0.1\cos(2\alpha)$. Here, $N=512$ and $\Delta t=5\times 10^{-5}$.}\label{fig:chemonecrosis}
\end{figure}

We present the tumor evolution in the first three rows in Fig. \ref{fig:chemonecrosis} where the radius of the (fixed) necrotic core $R_0=0.1, 1, 1.5$ corresponds to the rows from the first to the third row, and the columns from left to right correspond to different time $t$ as labeled. In the first column at $t=0$, all tumors are with the same initial shape $r=2.5+0.1\cos(2\alpha)$ but different necrotic radii. By comparing the first three rows from left to right, we see the tumor develop cavities along $y$-axis and such development is slower on tumors with a larger necrotic core. We also observe tumor further develops crown-like spikes protruding from their surfaces and such development tends to induce more spikes on tumors with a larger necrotic core.

In the fourth row, we change the chemotaxis constant from $\chi_\sigma=10$ to $\chi_\sigma=5$, with the remaining parameters the same as in the third row. Here we plot the time $t$ at the double amount of time to compare the development of morphology since smaller taxis decelerates tumor growth. By comparing the third and the fourth row, we see tumor with larger $\chi_\sigma$ eventually generates more aggressive spikes, which indicates the effect of \textit{chemotaxis} in triggering the development of spikes.


\subsection{Nonlinear simulation with a non-circular necrotic core}
In this section, instead of fixing the necrotic core with circular shape, we fix the necrotic core with a non-circular  (3-fold) shape and present the evolution of tumor morphology from the perspective of \textit{Proliferation} and \textit{Apoptosis}.

\paragraph{\textbf{Proliferation and Apoptosis}}
\begin{figure}
\includegraphics[width=\textwidth]{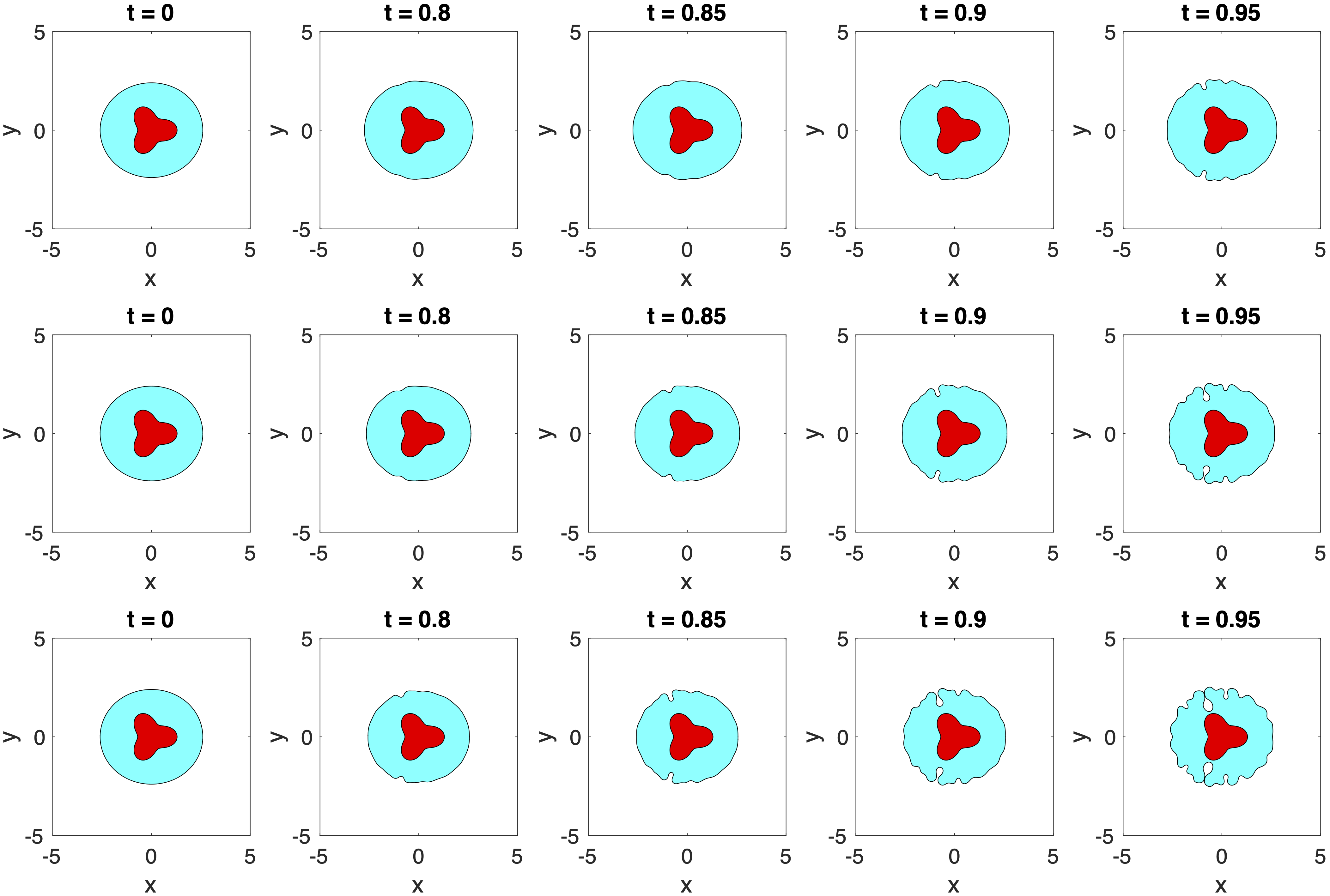}
\caption{The tumor morphologies with non-circular necrotic core $R_0=1+0.3\cos(3\alpha)$ for $\mathcal P=5,\mathcal A=0.25$ (first row), $\mathcal P=1,\mathcal A=0.25$ (second row) and $\mathcal P=1,\mathcal A=0.35$ (third row). The remaining parameters are $\beta=0.5, \chi_\sigma=10,\underline{\sigma}=0.2$ and $\mathcal{G}^{-1}=0.001$. The initial tumor boundary is $r=2.5+0.1\cos(2\alpha)$. Here, $N=512$ and $\Delta t=5\times 10^{-5}$.}\label{fig:proapop3}
\end{figure}

In Fig. \ref{fig:proapop3}, we demonstrate the dependence of tumor growth on \textit{proliferation} and \textit{apoptosis} with non-circular necrotic core $R_0=2+0.3\cos(3\alpha)$. The first cases in each row start with the same initial tumor boundary $r=2.5+0.1\cos(2\alpha)$ and necrotic boundary $R_0=1+0.3\cos(3\alpha)$. In the first row, we take $\mathcal P=5,\mathcal A=0.25$ and observe that the evolution of tumor morphology is influenced by the morphology of its necrotic core
and eventually develop two cavities on the left part of tumor.
In the second row, we decrease the proliferation rate $\mathcal P$ from $5$ to $1$ and find that at the same time $t$, the tumor attains a similar pattern with additional minor cavities
, which shows that the proliferation rate $\mathcal P$ 
stabilizes tumor morphology
In the third row, comparing to the second, we increase the apoptosis rate $\mathcal A$ from $0.25$ to $0.35$ and observe that the size of tumor shrinks and that when tumor cells are removed through apoptosis, more space is released for the aggressive patterns like protruding fingers to develop.



\subsection{The concentrations and fluxes for the control of necrotic region}
In this section, we study the control of a fixed necrotic core by observing the evolution of nutrient concentrations and fluxes on the boundaries. Here we consider two prototypical shapes of the necrotic region: circular and non-circular (3-fold). The parameters are the same as those used to produce the third row in Fig. \ref{fig:proapop3}.

\paragraph{\textbf{Control of a circular necrotic region}}

\begin{figure}

\begin{minipage}{1.0\linewidth}
\centering
\includegraphics[width=\textwidth]{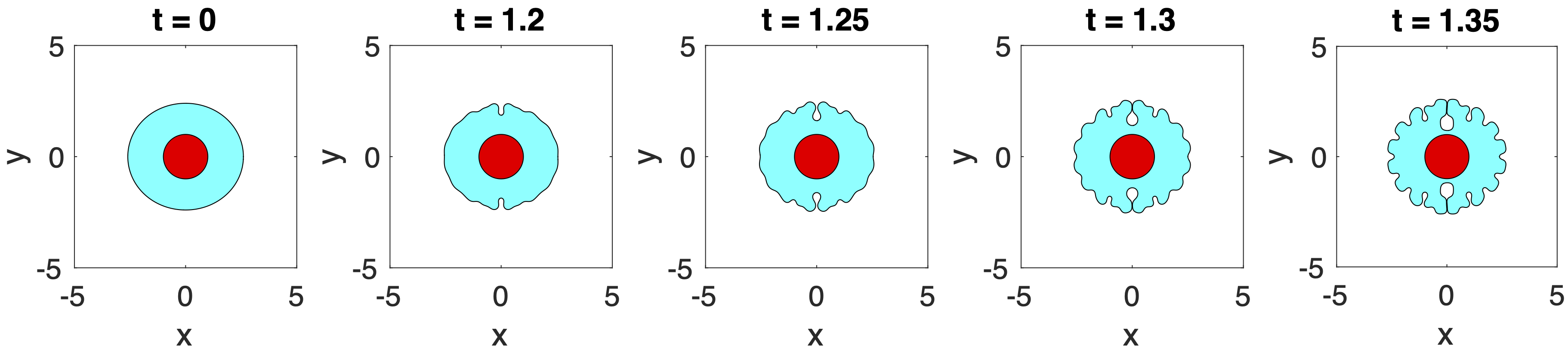}[a]
\end{minipage}
\begin{minipage}{1.0\linewidth}
\centering
\includegraphics[width=\textwidth]{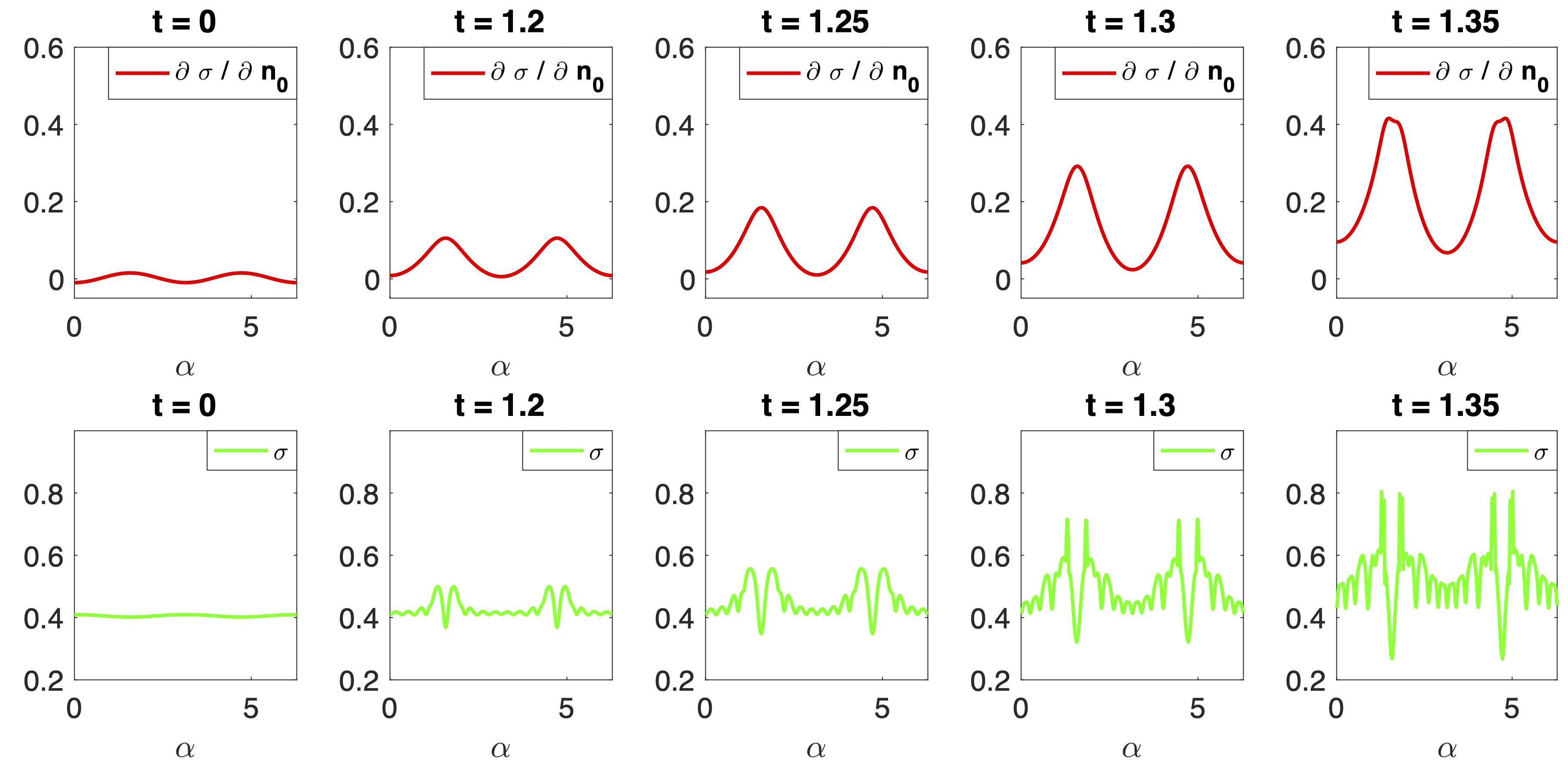}[b]
\end{minipage}
\begin{minipage}{1.0\linewidth}
\centering
\includegraphics[width=\textwidth]{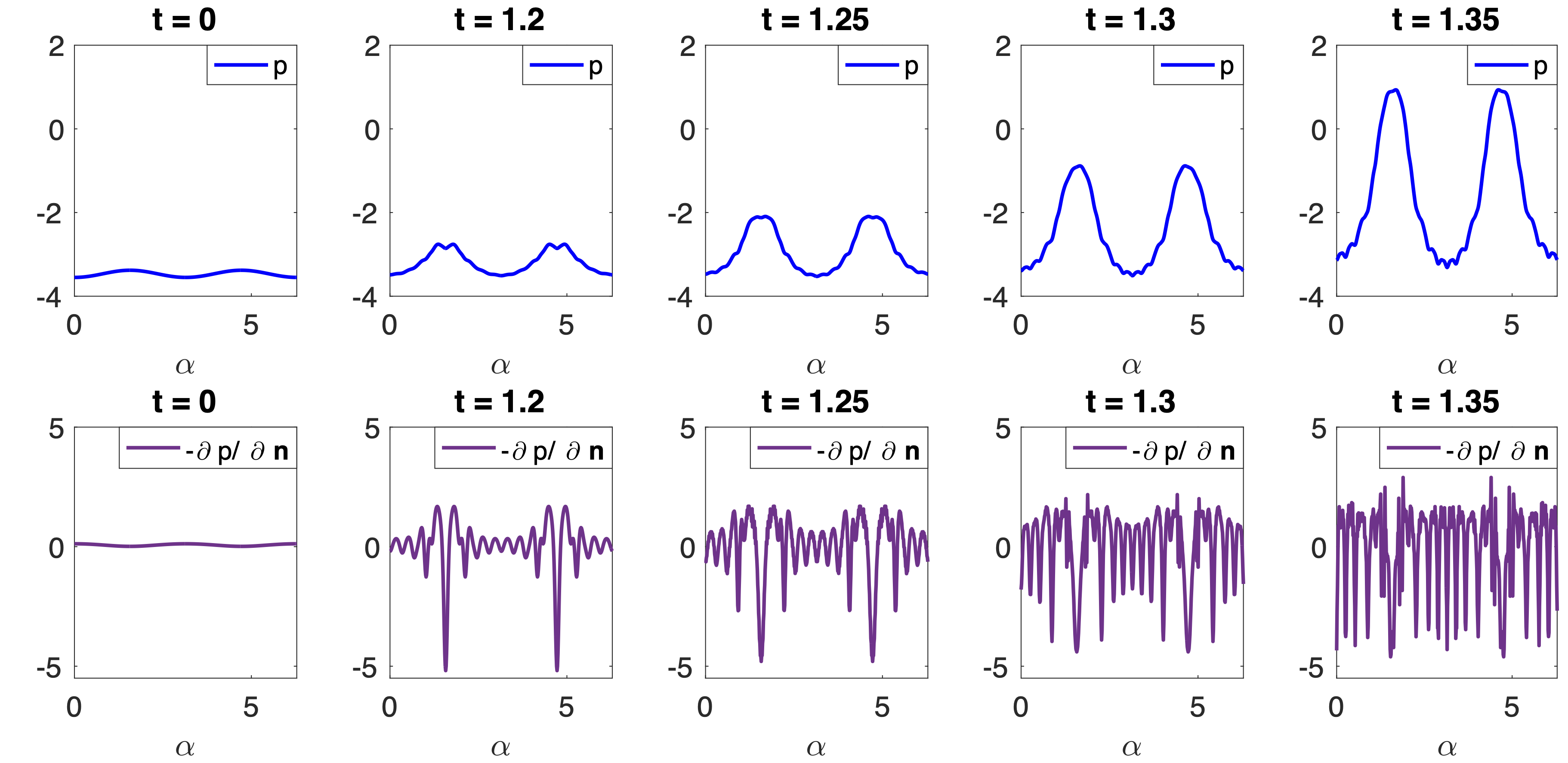}[c]
\end{minipage}

\caption{In [a]: The tumor morphologies with a circular necrotic boundary $R_0=1$. In [b]: the nutrient fluxes $\frac{\partial\sigma}{\partial\mathbf{n}_0}$ at necrotic boundary $\Gamma_0$ (first row), and the nutrient concentrations $\sigma$ at tumor boundary $\Gamma$ (second row). In [c]: the hydrostatic pressure $p$ at necrotic boundary $\Gamma_0$ (first row), and pressure fluxes $-\frac{\partial p}{\partial\mathbf{n}}$ at tumor boundary $\Gamma$ (second row). The remaining parameters are $\mathcal P=1,\mathcal A=0.35$, $\beta=0.5, \chi_\sigma=10,\underline{\sigma}=0.2$ and $\mathcal{G}^{-1}=0.001$. The initial tumor boundary is $r=2.5+0.1\cos(2\alpha)$. Here, $N=512$ and $5\times 10^{-5}$.}\label{fig:k0}
\end{figure}

In Fig. \ref{fig:k0}, we present the the control of a fixed circular necrotic core, where [a] shows the evolution of tumor morphologies with a circular necrotic boundary $R_0=1$, [b] shows the corresponding nutrient fluxes $\left.\frac{\partial \sigma}{\partial \mathbf{n}_0}\right|_{\Gamma_0}$ at the fixed necrotic boundary (first row) and the nutrient concentrations $\left.\sigma\right|_{\Gamma(t)}$ at the evolving tumor boundary (second row), [c] shows the corresponding hydrostatic pressure $\left.p\right|_{\Gamma_0}$  at the fixed necrotic boundary (first row) and the pressure fluxes $\left.\frac{\partial p}{\partial \mathbf{n}}\right|_{\Gamma(t)}$ at the evolving tumor boundary (second row). We remark that these four quantities are essentially the unknowns we solved at each time step using the boundary integral method (see Section \ref{sec:bimreformulation}; the modified pressure can be recovered back to the hydrostatic pressure by Eq. \eqref{eq:modpressure}.)

At $t=0$, we see the nutrient fluxes on the necrotic boundary reach two local maximums at the places where the two boundaries are closest $(\alpha = \pi/2,\ 3\pi/2)$, while the nutrient concentrations on the tumor boundary reach local minimums at those directions. At $t=1,2$ and $1.25$, the two local maximums of the nutrient fluxes on the necrotic boundary grow higher, while the two local minimums of the nutrient concentrations on the tumor boundary become deeper and then start to oscillate in response to the unstable tumor morphology. For $t=1.2$ and $1.35$, we see in [a] the developing of crown-like protrusions in tumor morphology. Interestingly, we observe that the envelope of the oscillating curve of nutrient concentrations except for the region close to local minima has a similar profile with the fluxes.

For the hydrostatic pressure, by comparing the first rows in \ref{fig:k0} [b] and [c], we see that the pressure level has a similar but more oscillating profile than that of the level of the nutrient fluxes at the necrotic boundary, which reveals a balance at the necrotic boundary for the size of the necrotic core to be maintained fixed. In the second row of \ref{fig:k0} [c], we see that the signs of the pressure fluxes reflects the region of the expansion $\left(-\frac{\partial p}{\partial \mathbf{n}}>0\right)$ or the shrinking $\left(-\frac{\partial p}{\partial \mathbf{n}}<0\right)$ of tumor. For example, as $t=1.2$, the negative sign of the two local minimums of $-\frac{\partial p}{\partial \mathbf{n}}$ corresponds to the cavities developed along the y-axis $(\alpha = \pi/2,\ 3\pi/2)$.


\paragraph{\textbf{Control of a non-circular region}}
\begin{figure}
\begin{minipage}{1.0\linewidth}
\centering
\includegraphics[width=\textwidth]{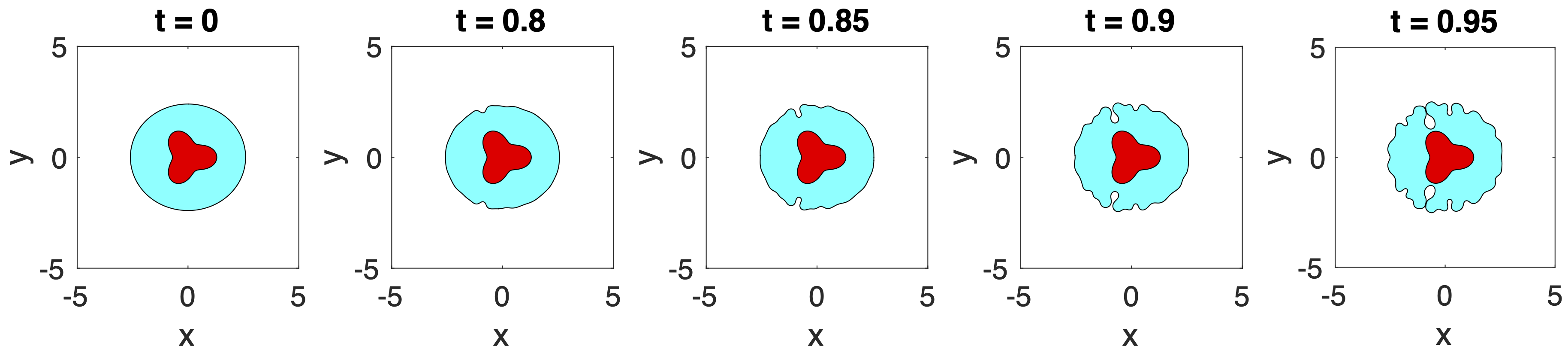}[a]
\end{minipage}
\begin{minipage}{1.0\linewidth}
\centering
\includegraphics[width=\textwidth]{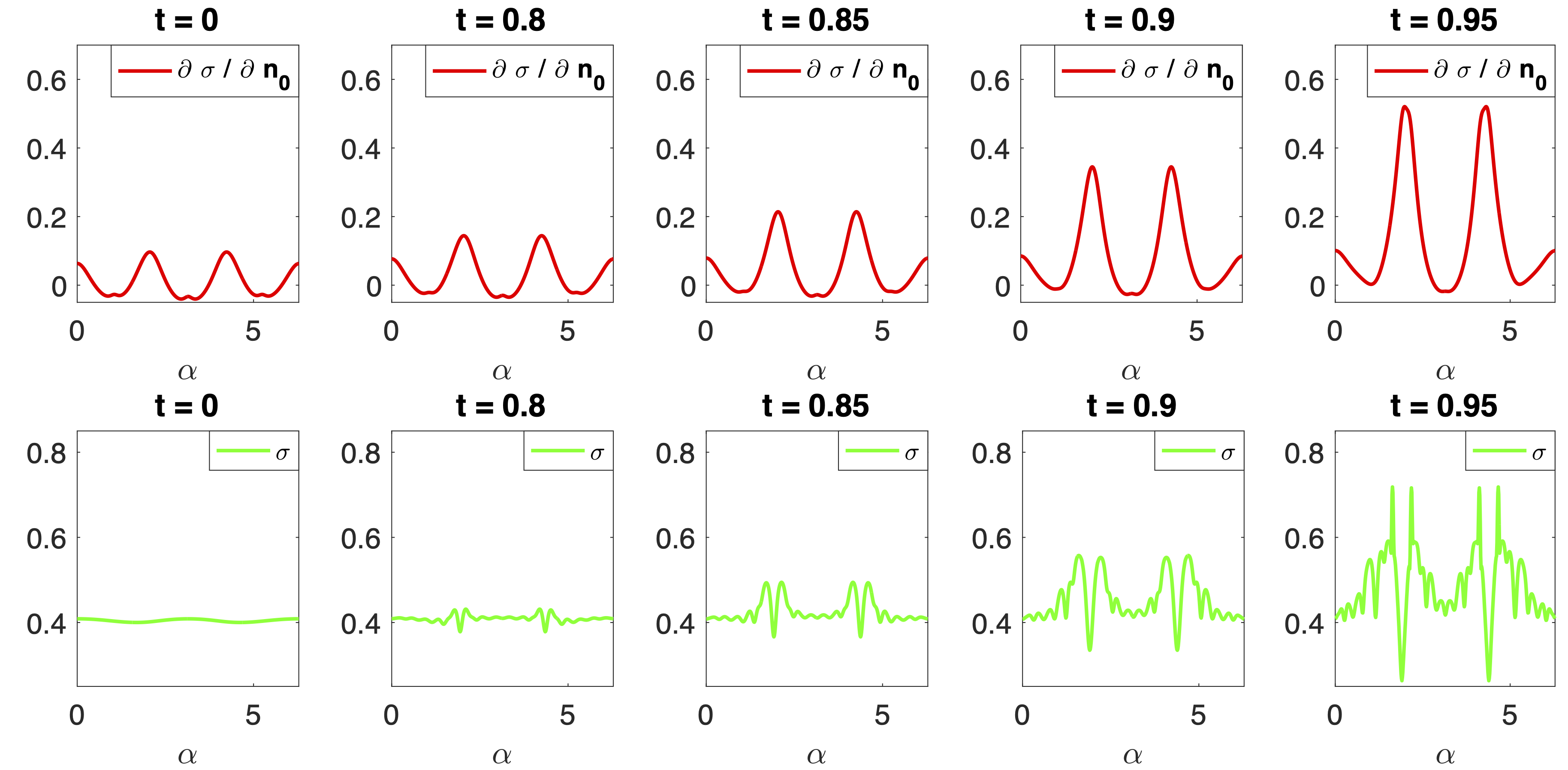}[b]
\end{minipage}
\begin{minipage}{1.0\linewidth}
\centering
\includegraphics[width=\textwidth]{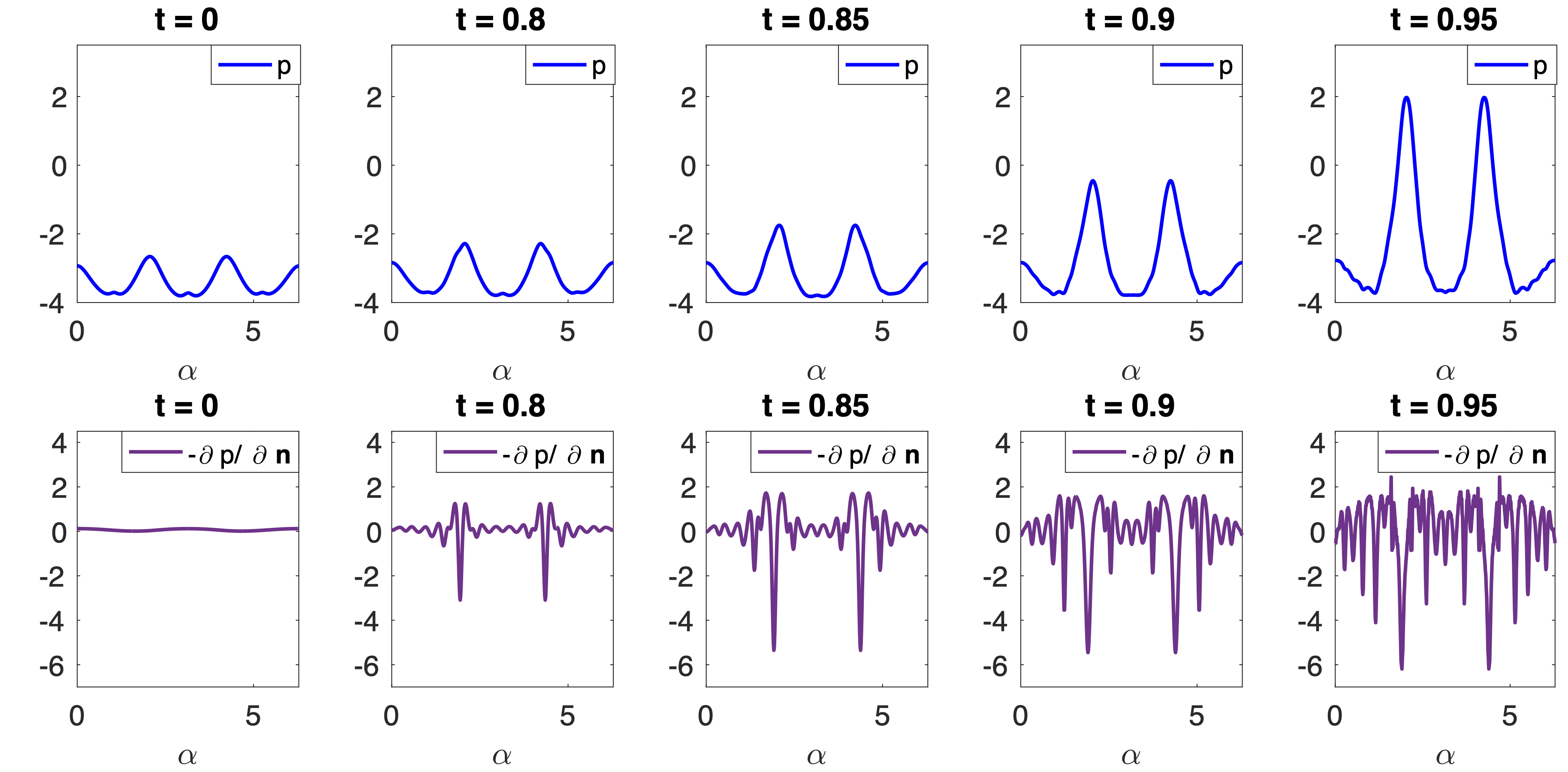}[c]
\end{minipage}
\caption{In [a]: The tumor morphologies with a non-circular (3-fold) necrotic boundary $R_0=1+0.3\cos(3\alpha)$. In [b]: the nutrient fluxes $\frac{\partial\sigma}{\partial\mathbf{n}_0}$ at necrotic boundary $\Gamma_0$ (first row), and the nutrient concentrations $\sigma$ at tumor boundary $\Gamma$ (second row). In [c]: the hydrostatic pressure $p$ at necrotic boundary $\Gamma_0$ (first row), and pressure fluxes $-\frac{\partial p}{\partial\mathbf{n}}$ at tumor boundary $\Gamma$ (second row). The remaining parameters are $\mathcal P=1,\mathcal A=0.35$, $\beta=0.5, \chi_\sigma=10,\underline{\sigma}=0.2$ and $\mathcal{G}^{-1}=0.001$. The initial tumor boundary is $r=2.5+0.1\cos(2\alpha)$. Here, $N=512$ and $5\times 10^{-5}$.}\label{fig:k3}
\end{figure}

\begin{figure}
\begin{minipage}{1.0\linewidth}
\centering
\includegraphics[width=\textwidth]{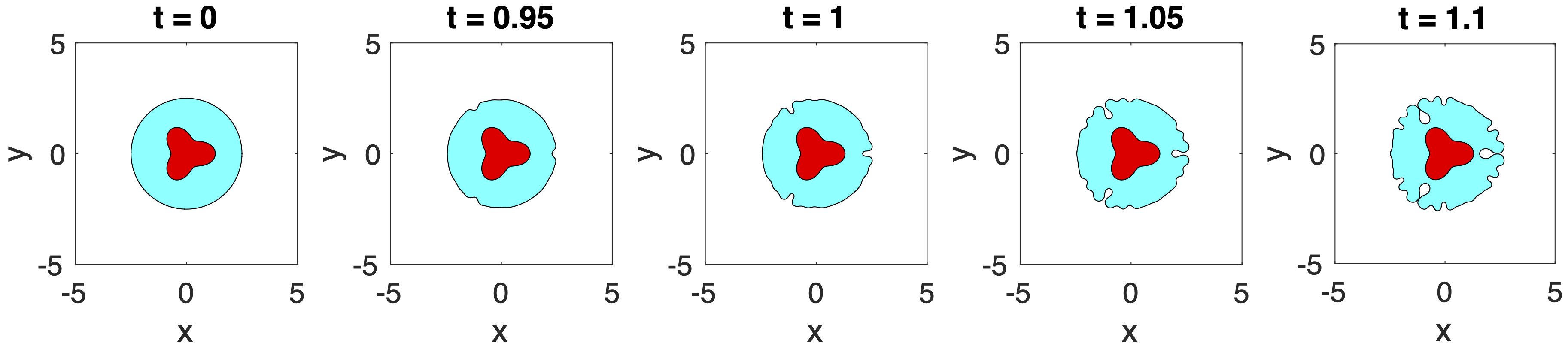}[a]
\end{minipage}
\begin{minipage}{1.0\linewidth}
\centering
\includegraphics[width=\textwidth]{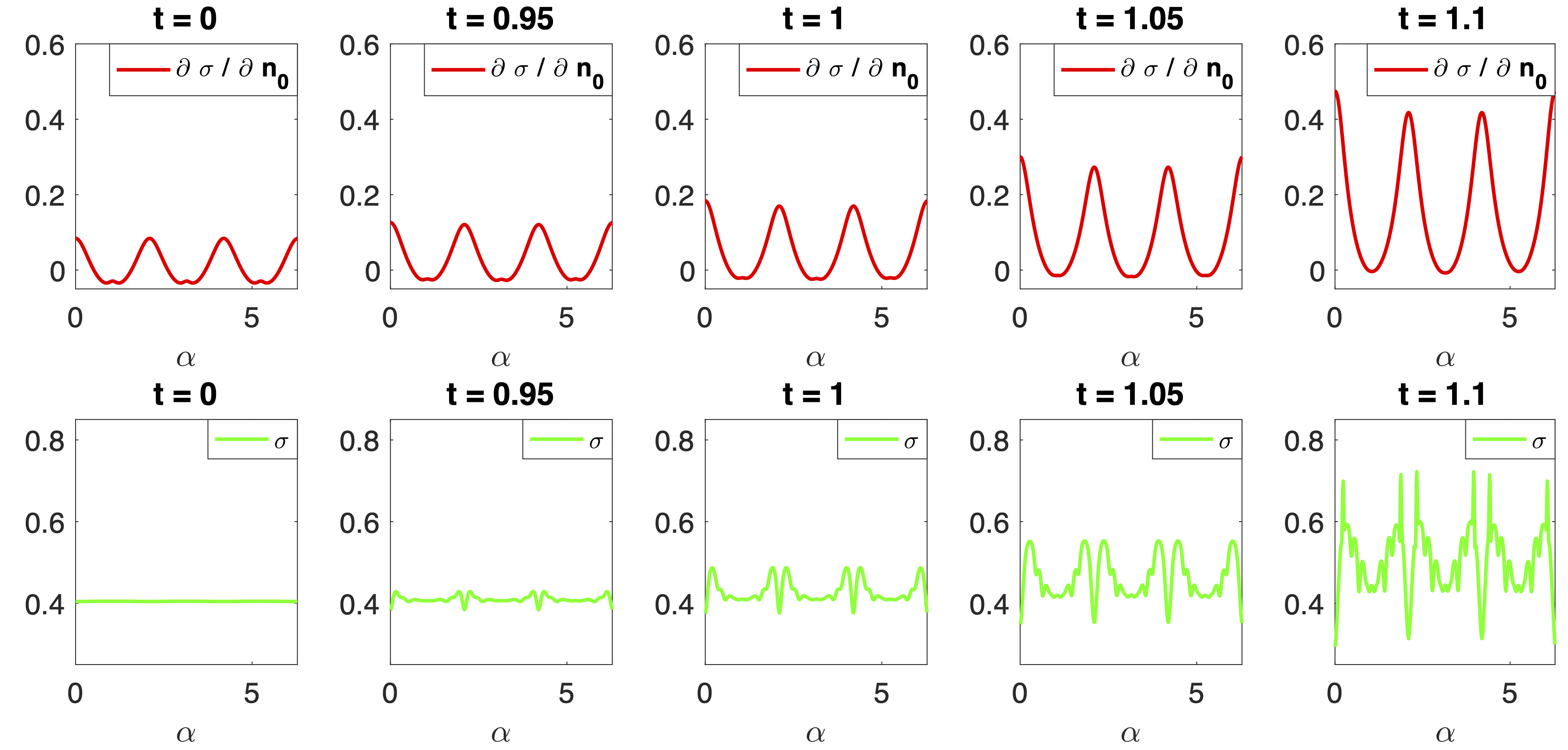}[b]
\end{minipage}
\begin{minipage}{1.0\linewidth}
\centering
\includegraphics[width=\textwidth]{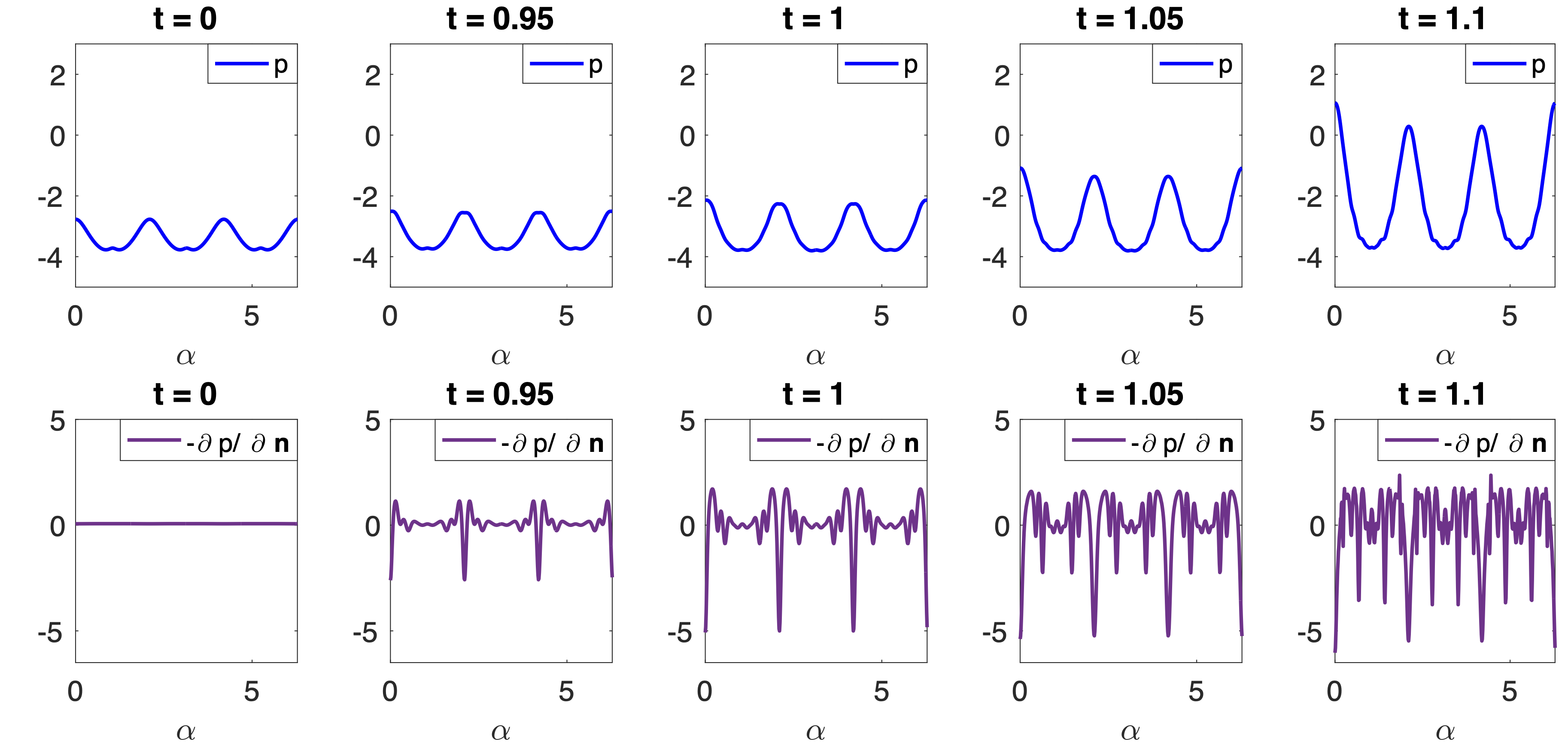}[c]
\end{minipage}
\caption{In [a]: The tumor morphologies with a non-circular (3-fold) necrotic boundary $R_0=1+0.3\cos(3\alpha)$. In [b]: the nutrient fluxes $\frac{\partial\sigma}{\partial\mathbf{n}_0}$ at necrotic boundary $\Gamma_0$ (first row), and the nutrient concentrations $\sigma$ at tumor boundary $\Gamma$ (second row). In [c]: the hydrostatic pressure $p$ at necrotic boundary $\Gamma_0$ (first row), and pressure fluxes $-\frac{\partial p}{\partial\mathbf{n}}$ at tumor boundary $\Gamma$ (second row). The remaining parameters are $\mathcal P=1,\mathcal A=0.35$, $\beta=0.5, \chi_\sigma=10,\underline{\sigma}=0.2$ and $\mathcal{G}^{-1}=0.001$. The initial tumor boundary is $r=2.5$. Here, $N=512$ and $5\times 10^{-5}$.}\label{fig:k3d0}
\end{figure}

Next, we present the the control of a fixed non-circular (3-fold) necrotic core in Fig. \ref{fig:k3}, where [a] is the evolution of tumor morphologies with a non-circular (3-fold) necrotic boundary $R_0=1+0.3\cos(3\alpha)$ (same as the third row in Fig. \ref{fig:proapop3}), [b] and [c] are in the same arrangement as in Fig. \ref{fig:k0}.

At $t=0$, we see the nutrient fluxes on the necrotic boundary reach 3 local maximums at the places where the two boundaries are closest $(\alpha = 0,\ 2\pi/3,\ 4\pi/3)$, while the nutrient concentrations look flat throughout the tumor boundary. As for $t=0.8$ and $0.85$, similar to Fig. \ref{fig:k0}, the two local maximums of the nutrient fluxes on the necrotic boundary grow higher and two local minimums of the nutrient concentrations on the tumor boundary appear, become deeper, and start to oscillate in response to the unstable tumor morphology. For $t=0.85$ and $0.95$, we see the development of splitting fingers in tumor morphology. Again, we observe that the envelope of the oscillating curve of nutrient concentrations except for the region close to local minima still has a similar profile with the fluxes as in Fig. \ref{fig:k0}. We also have similar observations for the hydrostatic pressure as in Fig. \ref{fig:k0}. However, by comparing the values in Fig. \ref{fig:k0} and \ref{fig:k3}, we can see the extremes in \ref{fig:k3} possess larger absolute magnitudes than those in \ref{fig:k0}, which indicates that the perturbation of the control shape of the necrotic core may cause the microenvironment of tumor to become more heterogeneous.

In Fig. \ref{fig:k3d0}, we change the initial tumor shape to be a pure circle with radius $r=2.5$, while other settings remain the same as in Fig. \ref{fig:k3}. Interestingly, the tumor morphology in [a] still evolves into a three-fold unstable pattern. In addition, we can see the values of the nutrient and the pressure at the boundaries in [b] and [c] all possess more evenly distributed and symmetric profiles than those in Fig. \ref{fig:k3}, which results from the more symmetric pattern of the tumor morphology in Fig. \ref{fig:k3d0} [a].



\section{Conclusions} 
\label{sec:conclusions}
In this paper, we have developed, analyzed, and solved numerically a tumor growth model that investigates the intratumoral structure using a controlled necrotic core and the extratumoral nutrient supply from vasculature, which is modeled by a Robin boundary condition at the tumor boundary. The model incorporates cell proliferation, death, \textit{angiogenesis, necrosis} and \textit{chemotaxis} up gradients of nutrients that are transported diffusionally from the vascularized tumor boundary and uptaken by tumor cells. Linear analysis, though limited to a simple geometry (performed here for a circular necrotic boundary), reveals the presence of rich pattern formation mechanisms via unstable tumor growth.

To gain insight into the nonlinear solutions, we developed a novel boundary integral method that naturally incorporates the Robin boundary condition to accurately and efficiently simulate the system. Direct layer potential representations were used for both pressure and nutrient fields, which enables us to obtain the value of the nutrient concentration with its fluxes and the hydrostatic pressure with its gradients on the interfaces accurately by solving two systems of integral equations. The tumor interface was evolved using a semi-implicit time-stepping method developed previously (\textit{e.g.}, \cite{hou1994removing,cristini2003nonlinear}). The method is spectrally accurate in space and second-order accurate in time.

With the advantage of boundary integral methods in addressing the complex boundary geometries and naturally incorporating the Robin boundary condition for nutrient field, our nonlinear simulations explored various unstable morphologies caused by \textit{angiogenesis, chemotaxis, necrosis} and cancer cell \textit{proliferation} and \textit{apoptosis}. When the tumor is growing with a fixed circular necrotic core, we investigate the effect of \textit{angiogenesis} in inhibiting the morphological instability and the effect of \textit{chemotaxis and necrosis} in destabilizing the tumor morphology, which is also observed in \cite{macklin2007nonlinear} and \cite{lu2020complex}. When the tumor is growing with a fixed non-circular necrotic core, we show the effect of \textit{proliferation} in accelerating tumor growth and stabilizing tumor morphology and \textit{apoptosis} acting in the opposite way. Finally, our numerical approach provides us with accurate physical quantities required to maintain the control of the shapes of the necrotic region for both nutrient and pressure fields.

In future work, we can investigate the motion of the necrotic boundary (\textit{e.g.}, two moving boundaries in the system) as considered in \cite{kohlmann2012necrotic}. Note that multiple moving interfaces have been studied recently in the context of Hele-Shaw flow \cite{Pedro1,Pedro2}. In addition, we can consider the secretion of the Tumor Angiogenesis Factor (TAF) from the moving necrotic boundary or the diffusion of the inhibitors (\textit{e.g.} drugs) from the moving tumor boundary. Furthermore, the nutrient concentration on the necrotic boundary need not be constant and the cell division and uptake rates need not be uniform, as assumed here. Another immediate extension is to use the Stokes equations (see for example \cite{minjhe2019}) to study the fluid properties of tumor cell and extracellular matrix mixtures through their different viscosities. The regulation of cell fates and motility, proliferation and apoptosis rates by mechanical and thermal stresses can also be incorporated. Finally, while we presented the results in two dimensions, similar behaviors are expected to hold qualitatively in three dimensions, and we plan to perform full 3D simulations to confirm this.

\section*{Acknowledgments}
ML acknowledges partial support from NSF-Simons Center for Multiscale Cell Fate Research through Interdisciplinary Opportunity Award IOA \#1901. M. L. is also grateful to Yifan Wang and Yuchi Qiu for stimulating discussions. WH is supported by  the National Science Foundation (NSF) grant DMS-2052685. SL acknowledges the support from the NSF, Division of Mathematical Sciences grant DMS-1720420 and ECCS-1307625. CL is partially supported by the NSF, Division of Mathematical Sciences grant DMS-1950868. JL acknowledges partial support from the NSF through grants DMS-1953410, DMS-1719960, and DMS-1763272 and the Simons Foundation (594598QN) for a NSF-Simons Center for Multiscale Cell Fate Research. JL also thanks the National Institutes of Health for partial support through grants 1U54CA217378-01A1 for a National Center in Cancer Systems Biology at UC Irvine and P30CA062203 for the Chao Family Comprehensive Cancer Center at UC Irvine.

\appendix

\section{Linear Stability Analysis}
\label{appendix:linear analysis}

The governing equations are

\begin{equation}\label{nutrientfieldap}
\left\{
\begin{array}{ccc}
\begin{aligned}
{\Delta}{\sigma}&={\sigma}&\text{in }&\Omega(t),\\
\left.\sigma\right|_{\Gamma_0}&=\underline\sigma &\text{on }&\Gamma_0, \\
\left.\frac{\partial\sigma}{\partial\mathbf{n}}\right|_{\Gamma(t)}&=\beta(1-\left.\sigma\right|_{\Gamma(t)}) &\text{on }&\Gamma(t).
\end{aligned}
\end{array}
\right.
\end{equation}

\begin{equation}\label{pressurefieldap}
\left\{
\begin{array}{ccc}
\begin{aligned}
\Delta p&=0
&\text{in }&\Omega(t),\\
\left.\frac{\partial p}{\partial \mathbf{n}_0}\right|_{\Gamma_0}&=\mathcal{P}\left.\frac{\partial\sigma}{\partial\mathbf{n}_0}\right|_{\Gamma_0} -\mathcal{P}\mathcal{A}\left.\frac{\mathbf{n}_0\cdot\mathbf{x}}{d}\right|_{\Gamma_0}
&\text{on }&\Gamma_0,\\
\left.p\right|_{\Gamma(t)}&={\mathcal{G}}^{-1}\left.\kappa\right|_{\Gamma(t)}
+(\mathcal{P}-{\chi_{\sigma}})\left.\sigma\right|_{\Gamma(t)}
-\mathcal{P}\mathcal{A}\left.\frac{\mathbf{x}\cdot \mathbf{x}}{2 d}\right|_{\Gamma(t)}
&\text{on }&\Gamma(t).
\end{aligned}
\end{array}
\right.
\end{equation}

\begin{equation}
V=-\left.\frac{\partial p}{\partial \mathbf{n}}\right|_{\Gamma(t)}
-\mathcal{P}
\left(\mathcal{A} \left.\frac{\mathbf{n}\cdot \mathbf{x}}{d}\right|_{\Gamma(t)}
-\beta\left(1-\left.\sigma\right|_{\Gamma(t)}\right)\right) \quad\text{on }\Gamma(t).
\end{equation}

Consider a perturbed tumor interface $\Gamma(t)$:

\begin{equation}
r(t)=R(t)+\delta(t) e^{i l \theta}\label{perturbedcircle}.
\end{equation}

In cylindrical coordinates modified Helmholtz equation satisfies

\begin{equation}
r^{-1}\left(r \sigma_{r}\right)_{r}+r^{-2} \sigma_{\theta \theta}+r^{-2} \sigma_{z z}- \sigma=0 \text { in } \Omega(t).
\end{equation}

Assume axial symmetry, $\textit{i.e}$, $\sigma=\sigma(r, \theta)$ is independent of $z$, then

\begin{equation}
r^{-1}\left(r \sigma_{r}\right)_{r}+r^{-2} \sigma_{\theta \theta}- \sigma=0\label{MHHeq}.
\end{equation}

\subsection{Radial solutions}
We first consider the radial solution, \textit{i.e.}, $\sigma=\sigma(r)$, then $\eqref{MHHeq}$ reduces to modified Bessel differential equation:
\begin{equation}
\left(r^{2} \frac{d^{2}}{d r^{2}}+r \frac{d}{d r}- r^{2}\right) \sigma(r)=0 \text { in } \Omega(t).
\end{equation}
Recall the general form of modified Bessel differential equation is
\begin{equation}
\left(x^{2} \frac{d^{2}}{d x^{2}}+x \frac{d}{d x}-\left( x^{2}+n^{2}\right)\right) y(x)=0.
\end{equation}

The general solutions are
\begin{eqnarray}{} 
y&=& a_{1} J_{n}(-i  x)+a_{2} Y_{n}(-i  x) \nonumber\\ 
&=&c_{1} I_{n}(x)+c_{2} K_{n}( x), 
\end{eqnarray}
where $J_{n}(x)$ is $a$ Bessel function of the first kind, $Y_{n}(x)$ is $a$ Bessel function of the second kind, $I_{n}(x)$ is a modified Bessel function of the first kind and $K_{n}(x)$ is a modified Bessel function of the second kind.

The following recurrence relation is useful in the linear analysis
\begin{eqnarray}
I_{n}^{\prime}(x)&=&\frac{1}{2}\left(I_{n-1}(x)+I_{n+1}(x)\right),\nonumber\\
I_{n}^{\prime}(x)&=&I_{n-1}(x)-\frac{n}{x} I_{n}(x)=\frac{n}{x} I_{n}(x)+I_{n+1}(x),\nonumber\\
I_{0}^{\prime}(x)&=&I_{1}(x).\label{besseliid}\\
K_{n}^{\prime}(x)&=&-\frac{1}{2}\left(K_{n-1}(x)+K_{n+1}(x)\right),\nonumber\\
K_{n}^{\prime}(x)&=&-K_{n-1}(x)-\frac{n}{x} K_{n}(x)=\frac{n}{x} K_{n}(x)-K_{n+1}(x),\nonumber\\
K_{0}^{\prime}(x)&=&-K_{1}(x).
\end{eqnarray}

We obtain $\sigma=A_{1} I_{0}\left(r\right)+A_{2} K_{0}\left( r\right) \operatorname{in}\  \Omega(t).$ Applying the boundary conditions on the circles $r=R_0,R$, we have
\begin{align}
A_{1} I_{0}\left(R_{0}\right)+A_{2} K_{0}\left(R_{0}\right)&=\underline{\sigma},\\
A_{1} I_{1}(R)-A_{2} K_{1}(R)&=\beta\left(1-A_{1} I_{0}(R)-A_{2} K_{0}(R)\right).\label{rblin}
\end{align}

Solving for $A_{1}, A_{2}$, we obtain
\begin{eqnarray}{}
A_1&=&
\frac{\underline{\sigma}\left(K_{1}(R)-\beta K_{0}(R)\right)+\beta K_{0}\left(R_{0}\right)}{K_{0}\left(R_{0}\right)\left(\beta I_{0}(R)+I_{1}(R)\right)+I_{0}\left(R_{0}\right)\left(K_{1}(R)-\beta K_{0}(R)\right)},\label{AA1}\\
A_2&=&
\frac{\underline{\sigma}\left(\beta I_{0}(R)+I_{1}(R)\right)-\beta I_{0}\left(R_{0}\right)}{K_{0}\left(R_{0}\right)\left(\beta I_{0}(R)+I_{1}(R)\right)+I_{0}\left(R_{0}\right)\left(K_{1}(R)-\beta K_{0}(R)\right)}\label{AA2}.
\end{eqnarray}

As $R_0\rightarrow 0$ we have
\begin{align}
A_1&\rightarrow \frac{1}{I_0(R)+\frac{I_1(R)}{\beta}},\\
A_2&\rightarrow 0.
\end{align}

\subsection{Perturbation of radial solutions}
Now we seek a solution of the modified Helmholtz' s equation on the perturbed circle given by $\eqref{perturbedcircle}$. Since $\delta$ is the perturbation size, following \cite{mullins1963morphological} we consider the Fourier expansion of the solution to the $1^{st}$ order in $\delta$:

\begin{equation}
\sigma(r, \theta)=\sigma_{0}(r)+\delta e^{i l \theta} \sigma_{1}(r) \text { in } \Omega(t).
\end{equation}

Note here that $r, \theta$ and $\delta$ are all functions of time $t,$ \textit{i.e.} $r=r(t), \theta=\theta(t), \delta=\delta(t)$.
Multiplying Eq. $\eqref{MHHeq}$ by $r^2$, we obtain

\begin{align}
\left(r^{2} \partial_{r}^{2}+r \partial_{r}+\partial_{\theta}^{2}- r^{2}\right)\left(\sigma_{0}(r)+\delta e^{i l \theta} \sigma_{1}(r)\right)&=0 &\text { in } \Omega(t),\\
\left(r^{2} \frac{d^{2}}{d r^{2}}+r \frac{d}{d r}- r^{2}\right) \sigma_{0}(r)&=0 &\text { in } \Omega(t),\\
\left(r^{2} \frac{d^{2}}{d r^{2}}+r \frac{d}{d r}-\left( r^{2}+l^{2}\right)\right) \sigma_{1}(r)&=0 &\text { in } \Omega(t).
\end{align}

Therefore it is sufficient to consider the expression :

\begin{equation}
\sigma=A_{1} I_{0}(r)+A_{2} K_{0}(r)+\delta e^{i l \theta}\left(B_{1} I_{l}(r)+B_{2} K_{l}(r)\right) \text{ in } \Omega(t).
\end{equation}

Apply nutrient boundary conditions on the interface $r=R+\delta e^{i l \theta}$ with $\delta\ll1$. 
(Orders higher than $O(\delta)$ are all discarded in the following calculations.) At $O(1)$, the equations are the same as the radial solution.

The equations at $O(\delta)$ determine the coefficients $B_{1}, B_{2}:$
\begin{align}
B_{1} I_{l}\left( R_{0}\right)+B_{2} K_{l}\left( R_{0}\right)&=0,\\
B_{1}\left(I_{l-1}(R)-\frac{l}{R} I_{l}(R)\right)-B_{2}\left(K_{l-1}(R)+\frac{l}{R} K_{l}(R)\right)\notag\\
+A_{1}\left(I_{0}(R)-\frac{1}{R} I_{1}(R)\right)+A_{2}\left(K_{0}(R)+\frac{1}{R} K_{1}(R)\right)
&=-\beta\left(B_{1} I_{l}(R)+B_{2} K_{l}(R)+A_{1} I_{1}(R)-A_{2} K_{1}(R)\right).\label{neclin}
\end{align}

Solving for $B_{1}, B_{2}$, we have
\begin{align}
B_1&=-\frac{K_{l}\left(R_{0}\right)\left(A_{1}\left((\beta R-1) I_{1}(R)+R I_{0}(R)\right)+A_{2}\left((1-\beta R) K_{1}(R)+R K_{0}(R)\right)\right)}{I_{l}\left(R_{0}\right)\left((l-\beta R) K_{l}(R)+R K_{l-1}(R)\right)+K_{l}\left(R_{0}\right)\left((\beta R-l) I_{l}(R)+R I_{l-1}(R)\right)},\\
B_2&=\frac{I_{l}\left(R_{0}\right)\left(A_{1}\left((\beta R-1) I_{1}(R)+R I_{0}(R)\right)+A_{2}\left((1-\beta R) K_{1}(R)+R K_{0}(R)\right)\right)}{I_{l}\left(R_{0}\right)\left((l-\beta R) K_{l}(R)+R K_{l-1}(R)\right)+K_{l}\left(R_{0}\right)\left((\beta R-l) I_{l}(R)+R I_{l-1}(R)\right)}.
\end{align}

Applying the definition for $A_1,A_2$ in Eq. $\eqref{AA1}, \eqref{AA2}$, we have
\scriptsize
\begin{align}
B_1&=-\frac{K_{l}\left(R_{0}\right)\left(\underline{\sigma}\left(\beta\left(\frac{1}{R}-\beta\right)+1\right)+\beta K_{0}\left(R_{0}\right)\left((\beta R-1) I_{1}(R)+R I_{0}(R)\right)+\beta I_{0}\left(R_{0}\right)\left((\beta R-1) K_{1}(R)-R K_{0}(R)\right)\right.}
{\left(K_{0}(R_{0})\left(\beta I_{0}(R)+I_{1}(R)\right)+I_{0}(R_{0})\left(K_{1}(R)-\beta K_{0}(R)\right)\right)
\left(I_{l}(R_{0})\left((l-\beta R) K_{l}(R)+R K_{l-1}(R)\right)+K_{l}(R_{0})\left((\beta R-l) I_{l}(R)+R I_{l-1}(R)\right)\right)},\\
B_2&=
-\frac{I_{l}(R_{0})\left(\underline{\sigma}\left(\beta^{2}-\frac{\beta}{R}-1\right)+\beta K_{0}(R_{0})\left(I_{1}(R)-R\left(\beta I_{1}(R)+I_{0}(R)\right)\right)+\beta I_{0}\left(R_{0}\right)\left((1-\beta R) K_{1}(R)+R K_{0}(R)\right)\right)}
{\left(K_{0}(R_{0})\left(\beta I_{0}(R)+I_{1}(R)\right)+I_{0}(R_{0})\left(K_{1}(R)-\beta K_{0}(R)\right)\right)
\left(I_{l}(R_{0})\left((l-\beta R) K_{l}(R)+R K_{l-1}(R)\right)+K_{l}(R_{0})\left((\beta R-l) I_{l}(R)+R I_{l-1}(R)\right)\right)}.
\end{align}

\normalsize 
As $\beta\rightarrow\infty$
\begin{align}
B_1&\rightarrow -\frac{\left(R K_{l}\left(R_{0}\right)\left(I_{1}(R) K_{0}\left(R_{0}\right)+I_{0}\left(R_{0}\right) K_{1}(R)\right)-\sigma K_{l}\left(R_{0}\right)\right)}{R\left(I_{0}\left(R_{0}\right) K_{0}(R)-I_{0}(R) K_{0}\left(R_{0}\right)\right)\left(I_{l}\left(R_{0}\right) K_{l}(R)-I_{l}(R) K_{l}\left(R_{0}\right)\right)},\\
B_2&\rightarrow -\frac{\left(\underline\sigma I_{l}\left(R_{0}\right)-R I_{l}\left(R_{0}\right)\left(I_{1}(R) K_{0}\left(R_{0}\right)+I_{0}\left(R_{0}\right) K_{1}(R)\right)\right)}{R\left(I_{0}\left(R_{0}\right) K_{0}(R)-I_{0}(R) K_{0}\left(R_{0}\right)\right)\left(I_{l}\left(R_{0}\right) K_{l}(R)-I_{l}(R) K_{l}\left(R_{0}\right)\right)}.
\end{align}

The nutrient $\sigma$ on $\Gamma$ is given by
\begin{eqnarray}{}
(\sigma)_{\Gamma}&=&\left(A_{1} I_{0}\left(r\right)+A_2K_0(r)+ (B_{1} I_{l}\left(r\right)+B_2 K_l(r)) \delta e^{i l \theta}\right)_{\Gamma}\nonumber\\
&=&A_{1} I_{0}\left(R\right)+A_2K_0(R)+\left(A_{1} I_{1}\left(R\right)-A_2K_1(R)+B_{1} I_{l}(R)+B_2 K_l(R)\right) \delta e^{i l \theta}.\label{nutatbdry}
\end{eqnarray}

The normal derivative of $\sigma$ on $\Gamma$ is given by

\begin{eqnarray} \left(\frac{\partial \sigma}{\partial n}\right)_{\Gamma} 
&=&\left(\frac{\partial \sigma}{\partial r}\right)_{\Gamma} \nonumber\\
&=&\left(\left(A_{1} I_{0}(r)+A_2K_0(r)+(B_{1} I_{l}(r)+B_2K_l(r)) \delta e^{i l \theta}\right)_{r}\right)_\Gamma\nonumber \\
&=&\left(A_{1} I_{1}\left( r\right)-A_2K_1(r)+\left(B_{1}\left(I_{l-1}\left(r\right)-\frac{l}{r} I_{l}\left( r\right)\right)-B_2\left(K_{l-1}\left(r\right)+\frac{l}{r} K_{l}\left( r\right)\right)\right) \delta e^{i l \theta}\right)_{\Gamma}\nonumber \\
&=& A_{1} I_{1}\left(R\right)-A_2K_1(R)\nonumber\\
&&+\left(A_{1}\left(I_{0}(R)- \frac{I_{1}( R)}{R}\right)+A_{2}\left(K_{0}(R)+ \frac{K_{1}( R)}{R}\right)\right.\nonumber\\
&&\left.+B_{1}\left(I_{l-1}(R)-l \frac{I_{l}( R)}{R}\right)-B_2\left(K_{l-1}(R)+l\frac{ K_{l}(R)}{R}\right)\right) \delta e^{i l\theta}.\label{nutfluxatbdry}
\end{eqnarray}

Similarly we seek a solution of Laplace equation on the perturbed circle given by Eq. $\eqref{perturbedcircle}$.
It is sufficient to consider the expression
\begin{equation}
p=C_1+C_2 \ln r+\delta e^{i l \theta} \left(D_1 r^{l}+\frac{D_2}{r^{l}}\right).
\end{equation}

For the perturbed circle defined by Eq.  $\eqref{perturbedcircle}, \kappa$ is given by

\begin{equation}
\kappa=\frac{1}{R}\left(1+\frac{l^{2}-1}{R} \delta e^{i l \theta}\right).
\end{equation}

On the interface we obtain
\begin{eqnarray}{}
(p)_{\Gamma}
&=&C_{1}+C_{2} \ln R+\delta e^{i l \theta}\left(\frac{C_{2}}{R}+D_{1} R^{l}+\frac{D_{2}}{R^{l}}\right) \nonumber\\
&=&\mathcal{G}^{-1}(\kappa)_{\Gamma}+\left(\mathcal{P}-\chi_{\sigma}\right)(\sigma)_{\Gamma}-\mathcal{P} \mathcal{A} \frac{(x \cdot x)_{\Gamma}}{4}\nonumber\\
&=&\scriptsize\mathcal{G}^{-1} \frac{1}{R}-\frac{\mathcal{P} \mathcal{A}}{4} R^{2}+\left(\mathcal{P}-\chi_{\sigma}\right)\left(A_{1} I_{0}(R)+A_{2} K_{0}(R)\right)\nonumber\\
&&+\left(\mathcal{G}^{-1} \frac{l^{2}-1}{R^{2}}-\frac{\mathcal{P} \mathcal{A}}{2} R+\left(\mathcal{P}-\chi_{\sigma}\right)\left(A_{1} I_{1}(R)-A_{2} K_{1}(R)+B_{1} I_{l}(R)+B_{2} K_{l}(R)\right)\right) \delta e^{i l \theta}.\nonumber
\end{eqnarray}

It is straightforward to derive that
\begin{align}{}
C_{1}+C_{2} \ln R&=\mathcal{G}^{-1} \frac{1}{R}-\frac{\mathcal{P} \mathcal{A}}{4} R^{2}+\left(\mathcal{P}-\chi_{\sigma}\right)\left(A_{1} I_{0}(R)+A_{2} K_{0}(R)\right),\label{C1}\\
\frac{C_{2}}{R}+D_{1} R^{l}+\frac{D_{2}}{R^{l}}&=\mathcal{G}^{-1} \frac{l^{2}-1}{R^{2}}-\frac{\mathcal{P} \mathcal{A}}{2} R+\left(\mathcal{P}-\chi_{\sigma}\right)\left(A_{1} I_{1}(R)-A_{2} K_{1}(R)+B_{1} I_{l}(R)+B_{2} K_{l}(R)\right).\label{D1}
\end{align}

The normal derivative of $p$ is given by
\footnotesize
\begin{eqnarray}{} \left(\frac{\partial p}{\partial n}\right)_{\Gamma_0} &=&\left(\frac{\partial p}{\partial r}\right)_{\Gamma_0}\nonumber\\
&=&\frac{C_{2}}{R_{0}}+\delta e^{i l \theta} l\left(D_{1} R_{0}^{l-1}-\frac{D_{2}}{R_{0}^{l+1}}\right) \notag\\ &=&\mathcal{P}\left(\frac{\partial \sigma}{\partial n_{0}}\right)_{\Gamma_{0}}-\mathcal{P} \mathcal{A} \frac{\left(n_{0} \cdot x\right)_{\Gamma_{0}}}{2}\notag\\
&=&\scriptsize\mathcal{P}\left(A_{1} I_{1}\left(R_{0}\right)-A_{2} K_{1}\left(R_{0}\right)+\left(B_{1}\left(I_{l-1}\left(R_{0}\right)-\frac{l}{R} I_{l}\left(R_{0}\right)\right)-B_{2}\left(K_{l-1}\left(R_{0}\right)+\frac{l}{R} K_{l}\left(R_{0}\right)\right)\right) \delta e^{i l \theta}\right)-\frac{\mathcal{P} \mathcal{A}}{2} R_{0}\notag\\
&=&\mathcal{P}\left(A_{1} I_{1}\left(R_{0}\right)-A_{2} K_{1}\left(R_{0}\right)+\left(B_{1} I_{l-1}\left(R_{0}\right)-B_{2} K_{l-1}\left(R_{0}\right)\right) \delta e^{i l \theta}\right)-\frac{\mathcal{P} \mathcal{A}}{2} R_{0}
\label{pressureatbdry},\nonumber
\end{eqnarray}

\normalsize 
where we have used Eq. $\eqref{neclin}$.

Then it is straightforward to derive that
\begin{align}
\frac{C_{2}}{R_{0}}&=\mathcal{P}\left(A_{1} I_{1}\left(R_{0}\right)-A_{2} K_{1}\left(R_{0}\right)\right)-\frac{\mathcal{P} \mathcal{A}}{2} R_{0}, \label{C2}\\
l\left(D_{1} R_{0}^{l-1}-\frac{D_{2}}{R_{0}^{l+1}}\right)&=\mathcal{P}\left(B_{1} I_{l-1}\left(R_{0}\right)-B_{2} K_{l-1}\left(R_{0}\right)\right).\label{D2}
\end{align}

Now solving $C_1,C_2$ by Eqs. $\eqref{C1},\eqref{C2}$, we have
\scriptsize
$$
\begin{aligned}{}
C_{1}&=\mathcal{P}\left(A_{1}\left(I_{0}(R)-R_{0} \ln (R) I_{1}\left(R_{0}\right)\right)+A_{2}\left(K_{0}(R)+R_{0} \ln (R) K_{1}\left(R_{0}\right)\right)\right)
-\chi_{\sigma} \left(A_{1} I_{0}(R)+A_{2} K_{0}(R)\right)
+\frac{\mathcal{P} \mathcal{A}}{2}\left(R_{0}^{2} \ln \left(R\right)-\frac{R^{2}}{2}\right)+\mathcal{G}^{-1} \frac{1}{R}, \nonumber \\
C_{2}&=\mathcal{P}\left(A_{1} I_{1}\left(R_{0}\right)-A_{2} K_{1}\left(R_{0}\right)\right) R_{0}-\frac{\mathcal{P} \mathcal{A}}{2} R_{0}^{2}.\nonumber
\end{aligned}
$$
\normalsize 
Next, solving $D_1,D_2$ by Eqs. $\eqref{D1},\eqref{D2}$, we have
\scriptsize
\begin{align} \notag 
D_{1}=&\frac{\mathcal P}{R^{2 l}+R_{0}^{2 l}}\left(R^{l}\left(A_{1} I_{1}(R)-A_{2} K_{1}(R)+B_{1} I_{l}(R)+B_{2} K_{l}(R)\right)-R_{0} R^{l-1}\left(A_{1} I_{1}\left(R_{0}\right)-A_{2} K_{1}\left(R_{0}\right)\right)+\frac{R_{0}^{l+1}}{l}\left(B_{1} I_{l-1}\left(R_{0}\right)-B_{2} K_{l-1}\left(R_{0}\right)\right)\right) \\
-&\frac{\chi_{\sigma}}{R^{2 l}+R_{0}^{2l}}\left(R^l \left( A_{1} I_{1}(R)-A_{2} K_{1}(R)+B_{1} I_{l}(R)+B_{2} K_{l}(R)\right)-R_{0} R^{l-1}\left(A_{1} I_{1}\left(R_{0}\right)-A_{2} K_{1}\left(R_{0}\right)\right)\right)+\frac{R^{l}}{R^{2 l}+R_{0}^{2 l}}\left(\frac{\mathcal P \mathcal{A} R_{0}^{2}}{2 R}-\frac{\mathcal P \mathcal{A} R}{2}+\mathcal{G}^{-1} \frac{l^{2}-1}{R^{2}}\right),
\\
D_{2}=&\frac{\mathcal P}{R^{2 l}+R_{0}^{2 l}}\left(R^{l}R_0^{2l}\left(A_{1} I_{1}(R)-A_{2} K_{1}(R)+B_{1} I_{l}(R)+B_{2} K_{l}(R)-\frac{R_{0}}{R} \left(A_{1} I_{1}\left(R_{0}\right)-A_{2} K_{1}\left(R_{0}\right)\right)\right)-\frac{R^{2l}R_{0}^{l+1}}{l}\left(B_{1} I_{l-1}\left(R_{0}\right)-B_{2} K_{l-1}\left(R_{0}\right)\right)\right)\\
-&\frac{\chi_{\sigma}}{R^{2 l}+R_{0}^{2l}}\left(R^l R_0^{2l} \left( A_{1} I_{1}(R)-A_{2} K_{1}(R)+B_{1} I_{l}(R)+B_{2} K_{l}(R)-\frac{R_{0}}{R} \left(A_{1} I_{1}\left(R_{0}\right)-A_{2} K_{1}\left(R_{0}\right)\right)\right)\right)+\frac{R^{l}R_0^{2l}}{R^{2 l}+R_{0}^{2 l}}\left(\frac{\mathcal P \mathcal{A} R_{0}^{2}}{2 R}-\frac{\mathcal P \mathcal{A} R}{2}+\mathcal{G}^{-1} \frac{l^{2}-1}{R^{2}}\right).
\end{align}
\normalsize 
Thus
\small
$$
\begin{aligned}
p=\ &\mathcal{P}\left(A_{1}\left(I_{0}(R)+I_{1}\left(R_{0}\right) R_{0} \ln\left(\frac{r}{R}\right)\right)+A_{2}\left(K_{0}(R)-K_{1}\left(R_{0}\right) R_{0} \ln \left(\frac{r}{R}\right)\right)\right)\notag\\
&-\chi_{\sigma}\left(A_{1} I_{0}(R)+A_{2} K_{0}(R)\right) 
-\frac{\mathcal{P} \mathcal{A}}{2}\left(R_{0}^{2} \ln \left(\frac{r}{R}\right)+\frac{R^{2}}{2}\right)+\frac{\mathcal{G}^{-1}}{R}\notag\\
&+\delta e^{i l \theta}\left(\mathcal{P}\left(\left(A_{1} I_{1}(R)-A_{2} K_{1}(R)+B_{1} I_{l}(R)+B_{2} K_{l}(R)-\left(A_{1} I_{1}\left(R_{0}\right)-A_{2} K_{1}\left(R_{0}\right)\right) \frac{R_{0}}{R}\right) \frac{(R r)^{l}+\left(R_{0} R\right)^{l}\left(\frac{R_{0}}{r}\right)^{l}}{R^{2l}+R_{0}^{2 l}}\right.\right. \notag\\
&+\left.\frac{R_{0}}{l}\left(B_{1}I_{l-1}\left(R_{0}\right)-B_{2}K_{l-1}\left(R_{0}\right)\right)\frac{\left(R_{0} r\right)^{l}-R^{2 l}\left(\frac{R_{0}}{r}\right)^{l}}{R^{2 l}+R_{0}^{2 l}}\right)\notag\\
&-\chi_\sigma\left(A_{1} I_{1}(R)-A_{2} K_{1}(R)+B_{1} I_{l}(R)+B_{2} K_{l}(R)-\left(A_{1} I_{1}(R_{0})-A_{2} K_{1}(R_{0})\right) \frac{R_{0}}{R}\right) \frac{(R_{0} r)^{l}-R^{2 l}\left(\frac{R_{0}}{r}\right)^{l}}{R^{2 l}+R_{0}^{2 l}}
\notag\\
&+\left.\left(\frac{\mathcal{P} \mathcal{A} R_{0}^{2}}{2R}-\frac{\mathcal{P} \mathcal{A} R}{2}+\mathcal{G}^{-1} \frac{l^{2}-1}{R^{2}}\right) \frac{(R r)^{l}+(R_{0} R)^{l}\left(\frac{R_{0}}{r}\right)^{l}}{R^{2l}+R_{0}^{2l}}\right).\notag
\label{plin}
\end{aligned}
$$
\normalsize 
And
\small
$$
\begin{aligned}
\left(\frac{\partial p}{\partial n}\right)_{\Gamma}=&\left(\mathcal{P}\left(A_{1} I_{1}\left(R_{0}\right)-A_{2} K_{1}\left(R_{0}\right)\right)-\frac{\mathcal{P} \mathcal{A} R_{0}}{2}\right) \frac{R_{0}}{R}\notag\\
&+\delta e^{i l \theta}\left({\mathcal{P}}\left(-{\left(A_{1} I_{1}\left(R_{0}\right)-A_{2} K_{1}\left(R_{0}\right)\right)} \frac{R_{0}}{R^{2}}\right.\right.\notag\\
&+\left(A_{1} I_{1}(R)-A_{2} K_{1}(R)+B_{1} I_{l}(R)+B_{2} K_{l}(R)-\left(A_{1} I_{1}\left(R_{0}\right)-A_{2} K_{1}\left(R_{0}\right)\right) \frac{R_{0}}{R}\right) \frac{R^{2 l}-R_{0}^{2 l}}{R^{2 l}+R_{0}^{2 l}} \frac{l}{R}\notag\\
&+\left.2\left(B_{1} I_{l-1}\left(R_{0}\right)-B_{2} K_{l-1}\left(R_{0}\right)\right) \frac{R^{l} R_{0}^{l}}{R^{2 l}+R_{0}^{2l} } \frac{R_{0}}{R}\right)\notag\\
&-\chi_{\sigma}\left(A_{1} I_{1}(R)-A_{2} K_{1}(R)+B_{1} I_{l}(R)+B_{2} K_{l}(R)-\left(A_{1} I_{1}\left(R_{0}\right)-A_{2} K_{1}\left(R_{0}\right)\right) \frac{R_{0}}{R}\right) \frac{R^{2 l}-R_{0}^{2 l}}{R^{2 l}+R_{0}^{2 l}} \frac{l}{R}\notag\\
&\left.+\frac{\mathcal{P} \mathcal{A}}{2}\left(\frac{R_{0}}{R}\right)^{2}
+\left(\frac{\mathcal{P} \mathcal{A} R_{0}^{2}}{2 R}-\frac{\mathcal{P} \mathcal{A} R}{2}+\mathcal{G}^{-1} \frac{l^{2}-1}{R^{2}}\right) \frac{R^{2 l}-R_{0}^{2 l}}{R^{2 l}+R_{0}^{2 l}} \frac{l}{R}\right).\label{pplin}
\end{aligned}
$$
\normalsize
Note that
\begin{equation}
n \cdot(\mathbf{x})_{\Gamma}=r+O\left(\delta^{2}\right)=R+\delta e^{i l \theta}+O\left(\delta^{2}\right).\label{normalxatbdry}
\end{equation}

Combining Eqs. $\eqref{nutatbdry},\eqref{pplin}$ and $\eqref{normalxatbdry}$,  we obtain

\small
\begin{align}{} 
V=&\frac{d{R}}{dt}+\frac{d\delta}{dt} e^{i l \theta} \nonumber\\
=&-\left(\frac{\partial p}{\partial n}\right)_{\Gamma}-\mathcal{P}\left(\mathcal{A}\frac{\mathbf{n}\cdot(\mathbf{x})_\Gamma}{2}-\beta(1-(\sigma)_\Gamma)\right)\notag\\
=&\ \mathcal{P}\left(\beta\left(1-A_{1} I_{0}(R)-A_{2} K_{0}(R)\right)-\frac{\mathcal{A}}{2} \frac{R^{2}-R_{0}^{2}}{R}\right)-\mathcal{P}\left(A_{1} I_{1}\left(R_{0}\right)-A_{2} K_{1}\left(R_{0}\right)\right) \frac{R_{0}}{R} \notag\\
&+\delta e^{i l \theta}\left(\frac{\mathcal{P} \mathcal{A}}{2}\left(\frac{R^{2}-R_{0}^{2}}{R^{2}} \frac{R^{2 l}-R_{0}^{2 l}}{R^{2 l}+R_{0}^{2 l}} l-\frac{R^{2}+R_{0}^{2}}{R^{2}}\right)-\mathcal{G}^{-1} \frac{l\left(l^{2}-1\right)}{R^{3}} \frac{R^{2 l}-R_{0}^{2 l}}{R^{2 l}+R_{0}^{2 l}}\right.\notag\\
&-\mathcal{P}{\beta}\left(A_{1} I_{1}(R)-A_{2} K_{1}(R)+B_{1} I_{l}(R)+B_{2} K_{l}(R)\right)+\mathcal{P}\left(\left(A_{1} I_{1}\left(R_{0}\right)-A_{2} K_{1}\left(R_{0}\right)\right) \frac{R_{0}}{R^{2}}\right. \notag\\
&-\left(A_{1} I_{1}(R)-A_{2} K_{1}(R)+B_{1} I_{l}(R)+B_{2} K_{l}(R)-\left(A_{1} I_{1}\left(R_{0}\right)-A_{2} K_{1}\left(R_{0}\right)\right) \frac{R_{0}}{R}\right) \frac{R^{2l}-R_{0}^{2 l}}{R^{2 l}+R_{0}^{2 l}} \frac{l}{R}\notag\\
&\left.-2\left(B_{1} I_{l-1}\left(R_{0}\right)-B_{2} K_{l-1}\left(R_{0}\right)\right) \frac{R^{l} R_{0}^{l}}{R^{2 l}+R_{0}^{2 l}} \frac{R_{0}}{R}\right)\notag\\
&\left.-\chi_{\sigma}\left(A_{1} I_{1}(R)-A_{2} K_{1}(R)+B_{1} I_{l}(R)+B_{2} K_{l}(R)-\left(A_{1} I_{1}\left(R_{0}\right)-A_{2} K_{1}\left(R_{0}\right)\right) \frac{R_{0}}{R}\right) \frac{R^{2l}-R_{0}^{2 l}}{R^{2 l}+R_{0}^{2 l}} \frac{l}{R}\right)\notag
.
\end{align}

\normalsize
Equating coefficients of like harmonics, we obtain
\begin{equation}
\frac{dR}{dt}=\mathcal{P} \left(\beta\left(1-A_1 I_0(R)-A_2 K_0(R)\right)-\frac{\mathcal{A}}{2}\frac{R^2-R_0^2}{R}\right)-\mathcal{P}\left(A_1 I_1\left(R_0\right)-A_2 K_1\left(R_0\right)\right)\frac{R_0}{R},
\end{equation}

or equivalently, simplifying by Eq. $\eqref{rblin}$, we have
\begin{equation}
\frac{dR}{dt}=\underbrace{
\mathcal{P} 
\left(A_1 I_1(R)-A_2 K_1(R)
-\frac{R_0 }{R}\left(A_1 I_1(R_0)-A_2 K_1(R_0)\right)\right)}_\text{Proliferation}
-
\underbrace{\frac{\mathcal{P}\mathcal{A}}{2}
\frac{R^2-R_0^2
}{R}}_\text{Apoptosis}
,
\end{equation}

\begin{equation}
R^{-1}\frac{dR}{dt}=\frac{\mathcal{P}}{R}\left(A_{1} I_{1}(R)-A_{2} K_{1}(R)\right)-\frac{\mathcal{P} \mathcal{A}}{2}\left(1-\left(\frac{R_{0}}{R}\right)^{2}\right)-\mathcal{P}\left(A_{1} I_{1}\left(R_{0}\right)-A_{2} K_{1}\left(R_{0}\right)\right) \frac{R_{0}}{R^{2}}.
\end{equation}

$O(\delta)$ terms 
\small
\begin{align}
\delta^{-1}\frac{d\delta}{dt}=&
\frac{\mathcal{P} \mathcal{A}}{2}\left(\frac{R^{2}-R_{0}^{2}}{R^{2}} \frac{R^{2 l}-R_{0}^{2 l}}{R^{2 l}+R_{0}^{2 l}} l-\frac{R^{2}+R_{0}^{2}}{R^{2}}\right)-\mathcal{G}^{-1} \frac{l\left(l^{2}-1\right)}{R^{3}} \frac{R^{2 l}-R_{0}^{2 l}}{R^{2 l}+R_{0}^{2 l}}\notag\\
&-\mathcal{P} \beta\left(A_{1} I_{1}(R)-A_{2} K_{1}(R)+B_{1} I_{l}(R)+B_{2} K_{l}(R)\right)+\mathcal{P}\left(\left(A_{1} I_{1}\left(R_{0}\right)-A_{2} K_{1}\left(R_{0}\right)\right) \frac{R_{0}}{R^{2}}\right.\notag\\
&-\left(\left(A_{1} I_{1}(R)-A_{2} K_{1}(R)+B_{1} I_{l}(R)+B_{2} K_{l}(R)\right)-\left(A_{1} I_{1}\left(R_{0}\right)-A_{2} K_{1}\left(R_{0}\right)\right) \frac{R_{0}}{R}\right) \frac{R^{2 l}-R_{0}^{2 l}}{R^{2 l}+R_{0}^{2 l}} \frac{l}{R}\notag\\
&-\left.2\left(B_{1} I_{l-1}\left(R_{0}\right)-B_{2} K_{l-1}\left(R_{0}\right)\right) \frac{R^{l} R_{0}^{l}}{R^{2 l}+R_{0}^{2 l}} \frac{R_{0}}{R}\right)\notag\\
&-\chi_{\sigma}\left(\left(A_{1} I_{1}(R)-A_{2} K_{1}(R)+B_{1} I_{l}(R)+B_{2} K_{l}(R)\right)-\left(A_{1} I_{1}\left(R_{0}\right)-A_{2} K_{1}\left(R_{0}\right)\right) \frac{R_{0}}{R}\right) \frac{R^{2 l}-R_{0}^{2 l}}{R^{2 l}+R_{0}^{2 l}} \frac{l}{R}.
\end{align}

\normalsize
The equation of shape perturbation is given by
\small
\begin{eqnarray}
&&\left(\frac{\delta}{R}\right)^{-1}
\frac{d}{dt}{\left(\frac{\delta}{R}\right)}\notag\\
&=&\delta^{-1} \frac{d \delta}{d t}-R^{-1} \frac{d R}{d t}\nonumber\\
&=&\overbrace{\mathcal{P} \mathcal{A}\left(\left(1-\left(\frac{R_{0}}{R}\right)^{2}\right)\left(1-\frac{2 R_{0}^{2 l}}{R^{2 l}+R_{0}^{2 l}}\right) \frac{l}{2}-\left(\frac{R_{0}}{R}\right)^{2}\right)}^\textbf{Apoptosis}-\overbrace{\mathcal{G}^{-1} \frac{l\left(l^{2}-1\right)}{R^{3}}\left(1-\frac{2 R_{0}^{2 l}}{R^{2 l}+R_{0}^{2 l}}\right)}^\textbf{Cell-cell adhesion} \notag\\ 
&&
-\overbrace{\mathcal{P}{\beta}\left(\frac{1}{R}+A_{1}\left(I_{1}(R)-\frac{I_{0}(R)}{R}\right)-A_{2}\left(K_{1}(R)+\frac{K_{0}(R)}{R}\right)+B_{1} I_{l}(R)+B_{2} K_{l}(R)\right)}^\textbf{Angiogenesis}\notag\\
&&+\overbrace{\chi_{\sigma}\left(A_{1} I_{1}(R)-A_{2} K_{1}(R)+B_{1} I_{l}(R)+B_{2} K_{l}(R)-\left(A_{1} I_{1}\left(R_{0}\right)-A_{2} K_{1}\left(R_{0}\right)\right) \frac{R_{0}}{R}\right) \left(1-\frac{2 R_{0}^{2 l}}{R^{2 l}+R_{0}^{2l}}\right) \frac{l}{R}
}^\textbf{Chemotaxis}\notag\\
&&+\overbrace{\mathcal{P}\left(\left(A_{1} I_{1}\left(R_{0}\right)-A_{2} K_{1}\left(R_{0}\right)\right) \frac{R_{0}}{R^{2}}\left(2+l\left(1-\frac{2 R_{0}^{2 l}}{R^{2 l}+R_{0}^{2l}}\right)\right)-2\left(B_{1} I_{l-1}\left(R_{0}\right)-B_{2} K_{l-1}\left(R_{0}\right)\right) \frac{R^{l} R_{0}^{l}}{R^{2 l}+R_{0}^{2 l}} \frac{R_{0}}{R}\right)
}^\textbf{Proliferation}\notag\\
&&-\overbrace{\mathcal{P}\left(\left(A_{1} I_{1}(R)-A_{2} K_{1}(R)+B_{1} I_{l}(R)+B_{2} K_{l}(R) \right)\left(1-\frac{2 R_{0}^{2 l}}{R^{2 l}+R_{0}^{2l}}\right) \frac{l}{R}\right)
}^\textbf{Proliferation}.
\end{eqnarray}

\normalsize
The apoptosis parameter $\mathcal{A}_{c}$ as a function of $R$ such that $\frac{d}{dt}{\left(\frac{{\delta}}{R}\right)}=0$ is given by	

\small
\begin{eqnarray}
\mathcal{A}_{c}&=& \left(\overbrace{ \mathcal{G}^{-1} \frac{l\left(l^{2}-1\right)}{\mathcal{P} R^{3}}\left(1-\frac{2 R_{0}^{2 l}}{R^{2 l}+R_{0}^{2 l}}\right)}^\textbf{Cell-cell adhesion}\right.\notag\\
&&+\overbrace{\beta\left(\frac{1}{R}+A_{1}\left(I_{1}(R)-\frac{I_{0}(R)}{R}\right)-A_{2}\left(K_{1}(R)+\frac{K_{0}(R)}{R}\right)+B_{1} I_{l}(R)+B_{2} K_{l}(R)\right)}^{\textbf{Angiogenesis}}\notag\\
&&-\overbrace{\frac{\chi_{\sigma}}{\mathcal{P}}\left(A_{1} I_{1}(R)-A_{2} K_{1}(R)+B_{1} I_{l}(R)+B_{2} K_{l}(R)-\left(A_{1} I_{1}\left(R_{0}\right)-A_{2} K_{1}\left(R_{0}\right)\right) \frac{R_{0}}{R}\right) \left(1-\frac{2 R_{0}^{2 l}}{R^{2 l}+R_{0}^{2l}}\right) \frac{l}{R}
}^\textbf{Chemotaxis to Proliferation}\notag\\
&&-
\left(A_{1} I_{1}(R_{0})-A_{2} K_{1}(R_{0}\right) \frac{R_{0}}{R^{2}}
\left(2+l\left(1-\frac{2R_{0}^{2l}}{R^{2l}+R_{0}^{2l}}\right)\right)\notag\\
&&+\left(A_{1} I_{1}(R)-A_{2} K_{1}(R)+B_{1} I_{l}(R)+B_{2} K_{l}(R)\right)\left(1-\frac{2 R_{0}^{2 l}}{R^{2 l}+R_{0}^{2 l}}\right) \frac{l}{R}\notag\\
&&\left.+2\left(B_{1} I_{l-1}(R_{0})-B_{2} K_{l-1}(R_{0})\right) 
\frac{R^{l} R_{0}^{l}}{R^{2 l}+R_{0}^{2 l}} \frac{R_{0}}{R}\right)\notag\\
&&/\left(\left(1-\left(\frac{R_{0}}{R}\right)^{2}\right)\left(1-\frac{2 R_{0}^{2 l}}{R^{2 l}+R_{0}^{2 l}}\right) \frac{l}{2}-\left(\frac{R_{0}}{R}\right)^{2}\right).
\nonumber
\end{eqnarray}

\normalsize

\section{The Evaluation of the Boundary Integrals}
\label{Appendix:Evaluation of BIM}
With the integral formulation above, we assume interface curves $\Gamma$ and $\Gamma_{\infty}$ are analytic and given by $\big\{\mathbf{x}(\alpha,t)=(x(\alpha,t),y(\alpha,t): 0\leq \alpha \leq 2 \pi \big\}$, where $\mathbf{x}$ is $2 \pi$-periodic in the parametrization $\alpha$. The unit tangent and normal(outward) vectors can be calculated as $\mathbf{s}=(x_\alpha,y_\alpha)/s_\alpha$, $\mathbf{n}=(y_\alpha,-x_\alpha)/s_\alpha$, where the local variation of the arclength $s_\alpha=\sqrt{x_\alpha^2+y_\alpha^2}$. Subscripts refer to partial differentiation.
We track the interfaces $\Gamma$ and $\Gamma_{\infty}$ by introducing N marker points to discretize the planar curves, parametrized by $\alpha_j=jh$, $h=\frac{2\pi}{N}$, $N$ is a power of $2$. Here we focus on the numerical evaluation of integrals following \cite{jou1997,li2011boundary,minjhe2019}. A rigorous convergence and error analysis of the boundary integral method for a simplified tumor problem can be found in \cite{Wenrui2018}.

\paragraph{\textbf{Computation of the single-layer potential type integral}}
In Eqs. \eqref{bimnut} and \eqref{bimpressure}, the single-layer potential type integrals contain the Green functions with a logarithmic singularity at $r=0$. They can be rewritten in the following form under the parametrization $\alpha$
\begin{equation}\label{eq50}
    \int_{\Gamma}
    \Psi(\alpha,\alpha')\phi(\alpha')s_{\alpha}(\alpha')d\alpha',
\end{equation}where $\Psi$ are the Green functions, which can be $G$ or $\Phi$ from Eqs. \eqref{eq34} or \eqref{fundamentalLaplace} and $\Gamma$ may be either $\Gamma(t)$ or $\Gamma_{0}$.
We may decompose the Green functions as below
\begin{equation}\label{eq51}
    \Phi(\alpha,\alpha')=
    -\frac{1}{2 \pi} \ln r
    =-\frac{1}{2\pi}\left(
    \ln{2 \left| \sin{\frac{\alpha-\alpha'}{2}}\right|}
    +\left[\ln r-\ln{2 \left| \sin{\frac{\alpha-\alpha'}{2}}\right|}\right]
    \right),
\end{equation}

\begin{equation}\label{eq52}
    G(\alpha,\alpha')=
    \frac{1}{2 \pi} K_{0}( r)=
    -\frac{1}{2 \pi}
    \left(I_{0}(r)\ln{2 \left| \sin{\frac{\alpha-\alpha'}{2}}\right|}
    +\left[-K_{0}(r)
    -I_{0}(r)\ln{2 \left| \sin{\frac{\alpha-\alpha'}{2}}\right|}
    \right]\right),
\end{equation}
where $I_{0}$ is a modified Bessel function of the first kind,
$r=|\mathbf{x}(\alpha)-\mathbf{x}'(\alpha ')|$. The square brackets on the right-hand side of Eqs. \eqref{eq51}, \eqref{eq52}
have removable singularity at $\alpha=\alpha'$, since $r=
s_{\alpha}\left|\alpha-\alpha'\right|
\sqrt{1+\mathscr{O}(\alpha-\alpha')}
=s_{\alpha}\left|\alpha-\alpha'\right|
(1+\mathscr{O}(\alpha-\alpha'))$ for $\alpha \approx \alpha'$, where $\mathscr{O(\alpha-\alpha')}$ denotes a smooth function that vanishes as $\alpha\rightarrow\alpha'$, and since $K_0$ has the expansion
\begin{equation}
K_{0}(z)=-\left(\log \frac{z}{2}+C\right) I_{0}(z)+\Sigma_{n=1}^{\infty} \frac{\psi(n)}{(n !)^{2}}\left(\frac{z}{2}\right)^{2 n}.
\end{equation}
Thus, for an analytic and $2\pi$-periodic function $f(\alpha,\alpha')$, a standard trapezoidal rule or alternating point rule can be used to evaluate the integral
\begin{equation}\label{eq53}
    \int_{0}^{2\pi}
f(\alpha,\alpha')
\ln{\frac{r}{2 \left| \sin{\frac{\alpha-\alpha'}{2}}\right|}}
d\alpha'.
\end{equation}
The remaining terms on the right-hand side of Eqs. \eqref{eq51}, \eqref{eq52}
have logarithmic singularity and can be evaluated through the following spectrally accurate quadrature \cite{kress1995numerical}
\begin{equation}\label{eq54}
    \int_{0}^{2\pi}f(\alpha_i,\alpha')
\ln{2 \left| \sin{\frac{\alpha_i-\alpha'}{2}}\right|}
d\alpha'\approx
\Sigma_{j=0}^{2m-1}q_{\left|j-i\right|}f(\alpha_i,\alpha_j),
\end{equation}where $m=\frac{N}{2}$, $\alpha_i=\frac{\pi i}{m}$ for $i=0,1,...,2m-1$, and weight coefficients
\begin{equation}\label{eq55}
    q_j=-\frac{\pi}{m}\Sigma_{k=1}^{m-1}\frac{1}{k}\cos{\frac{kj\pi}{m}}-\frac{(-1)^j\pi}{2m^2} , \text{for } j=0,1,...,2m-1.
\end{equation}
\paragraph{\textbf{Computation of the double-layer potential-type integral}}
In Eqs. \eqref{bimnut} and \eqref{bimpressure}, the double-layer potential type integrals contain the Green functions with singularity at $r=0$. 
They can be rewritten as in the following form under the parametrization $\alpha$
\begin{equation}\label{eq56}
    \int_{\Gamma}
    \frac{\partial \Psi(\alpha,\alpha')}{\partial \mathbf{n}(\alpha')}\phi(\alpha')s_{\alpha}(\alpha')d\alpha',
\end{equation} where $\Psi$ are the Green functions $G$ or $\Phi$ from Eqs. \eqref{eq34} or \eqref{fundamentalLaplace} and $\Gamma$ may be either $\Gamma(t)$ or $\Gamma_{0}$. 
Further, in Eq. \eqref{bimpressure},
\begin{equation}\label{eq57}
    \frac{\partial \Phi(\alpha,\alpha')}{\partial \mathbf{n}(\alpha')}s_{\alpha}(\alpha')=
    h(\alpha,\alpha')\frac{1}{r},
\end{equation} where the auxiliary function $h(\alpha,\alpha')=\frac{(\mathbf{x(\alpha)}-\mathbf{x(\alpha')})\cdot\mathbf{n(\alpha')}s_\alpha(\alpha')}{2\pi r}$ with $r=\left|\mathbf{x(\alpha)}-\mathbf{x(\alpha')}\right|$. Note that $h(\alpha,\alpha')\sim\mathscr{O}(\alpha-\alpha')$.
Since $\frac{\partial \Phi}{\partial \mathbf{n}}$ has no logarithmic singularity, we may simply use the alternating point rule to evaluate it.
For $\frac{\partial G}{\partial \mathbf{n}}$ in Eq.  \eqref{bimnut}, we decompose it as below
\begin{equation}\label{eq58}
    \frac{\partial G(\alpha,\alpha')}{\partial \mathbf{n}(\alpha')}s_{\alpha}(\alpha')=
    h(\alpha,\alpha') K_{1}(r)=
    g_1(\alpha,\alpha')\ln{2 \left| \sin{\frac{\alpha-\alpha'}{2}}\right|}
    +g_2(\alpha,\alpha'),
\end{equation}where $g_1(\alpha,\alpha')$ and $g_2(\alpha,\alpha')$ are analytic and $2\pi$-periodic functions with
\begin{equation}\label{eq59}
    g_1(\alpha,\alpha')=h(\alpha,\alpha')I_1(r),
\end{equation}
\begin{equation}\label{eq60}
    g_2(\alpha,\alpha')=h(\alpha,\alpha')
    \left[K_{1}(r)
    -I_{1}(r)\ln{2 \left| \sin{\frac{\alpha-\alpha'}{2}}\right|}
    \right],
\end{equation}
where we have used the fact
\begin{equation}
\frac{d}{dr}K_0(r)=-K_1(r).
\end{equation}
Since $K_1$ has the expansion
\begin{equation}
K_{1}(z)=\frac{1}{z}+\left(\log \frac{z}{2}+C\right) I_{1}(z)-\frac{1}{2} \sum_{n=0}^{\infty} \frac{\psi(n+1)+\psi(n)}{n !(n+1) !}\left(\frac{z}{2}\right)^{2 n+1},
\end{equation}
the  square bracket on the right-hand side of Eq. $\eqref{eq60}$ also has removable singularity at $\alpha=\alpha'$ thus the integral involving $g_2(\alpha,\alpha')$ can be evaluated by a standard trapezoidal rule or alternating point rule. Note that
\begin{equation}\label{eq61}
    g_2(\alpha,\alpha)=\frac{h(\alpha,\alpha)}{ r}=\frac{1}{4\pi}
    \frac{x_\alpha y_{\alpha\alpha}-x_{\alpha\alpha}y_\alpha}{x_\alpha^2+y_\alpha^2}.
\end{equation}The first term on the right-hand side of Eq. $\eqref{eq58}$ is still singular and evaluated through the quadrature given in Eqs. \eqref{eq54} and \eqref{eq55}.\\
To summarize, using $\text{Nystr\"om}$ discretization with the Kress quadrature rule described above, we reduce the boundary integral Eqs. \eqref{bimnut} and \eqref{bimpressure} to two dense linear systems with the unknowns as the discretization of $p$, $\frac{\partial\sigma}{\partial\mathbf{n_0}}$ on $\Gamma_0$ and $\sigma$, $\frac{\partial p}{\partial \mathbf{n}}$ on $\Gamma(t)$, which can be solved using an iterative solver, \textit{e.g.}, GMRES \cite{saad1986gmres}.

\section{The Evolution of the Interface}
\label{Appendix:Evolution of the interface}
As indicated by \cite{hou1994removing}, the curvature-driven motion introduces high-order derivatives, both non-local and nonlinear, into the dynamics through the Laplace-Young condition at the interface. Explicit time integration methods thus suffer from severe stability constraints and implicit methods are difficult to apply since the stiffness enters nonlinearly. Hou et al. resolved these difficulties by adopting the $\theta-L$ formulation and the small-scale decomposition (SSD), which we apply here.
\paragraph{\textbf{$\theta-L$ formulation}}
This formulation helps to circumvent the problem of point clustering. Consider a point $\mathbf{x}(\alpha,t)=(x(\alpha,t),y(\alpha,t))\in \Gamma(t)$. Denote the unit tangent and normal (outward) vectors as $\hat{\mathbf{s}}=(x_\alpha,y_\alpha)/s_\alpha$ and $\hat{\mathbf{n}}=(y_\alpha,-x_\alpha)/s_\alpha$, the normal velocity and tangent velocity by $V(\alpha,t)=u\cdot\hat{\mathbf{n}}$ and $T(\alpha,t)=u\cdot\hat{\mathbf{s}}$, respectively, where $u=\mathbf{x}_t=V \hat{\mathbf{n}}+T \hat{\mathbf{s}}$ gives the motion of $\Gamma(t)$. The tangent angle that the planar curve $\Gamma(t)$ forms with the horizontal axis at $\mathbf{x}$, called $\theta$, satisfies $\theta=\tan^{-1}{\frac{y_\alpha}{x_\alpha}}$. The length of one period of the curve is $L(t)=\int_0^{2\pi}s_\alpha d\alpha$, where $s_\alpha$, the derivative of the arclength, satisfies $s_\alpha^2=x_\alpha^2+y_\alpha^2$. Differentiating these two equations in time, we obtain the following evolution equations:
\begin{equation}\label{eq62}
    \theta_t=\kappa T - V_s=\frac{1}{s_\alpha}(\theta_\alpha T- V_\alpha),
\end{equation}
\begin{equation}\label{eq63}
    s_{\alpha t}=(T_s+\kappa V)s_\alpha=T_\alpha+\theta_\alpha V.
\end{equation}
Instead of using the $(x,y)$ coordinates, $(L,\theta)$ becomes the dynamical variables. The unit tangent and normal vectors become $\Hat{\mathbf{s}}=(\cos{\theta},\sin{\theta})$, $\Hat{\mathbf{n}}=(\sin{\theta},-\cos{\theta})$.

 The normal velocity $V$ is calculated using Eq. \eqref{normal velocity}. The tangent velocity $T$ is chosen (independent of the morphology of the interface) such that the marker points are equally spaced in arclength to prevent point clustering:
\begin{equation}\label{tangent velocity}
    T(\alpha,t)= \frac{\alpha}{2\pi}\int_0^{2\pi}\theta_{\alpha'}V' d\alpha'-\int_0^\alpha \theta_{\alpha'}V' d\alpha'.
\end{equation}
It follows that $s_\alpha$ is independent of $\alpha$ thus is everywhere equal to its mean:
\begin{equation}\label{eq65}
    s_\alpha=\frac{1}{2\pi}\int_0^{2\pi}s_\alpha(\alpha,t)d\alpha=\frac{L(t)}{2\pi}.
\end{equation}
The procedure for obtaining the initial equal arclength parametrization is presented in ``Appendix B" of \cite{baker1990connection}. The
idea is to solve the nonlinear equation
\begin{equation}\label{eq66}
    \int_0^{\alpha_j} s_{\beta}d\beta=\frac{j}{N}L
\end{equation}for $\alpha_j$ using Newton's method and evaluate the equal arclength marker points $\mathbf{x}(\alpha_j)$by interpolation in Fourier space.
We may recover the interface by simply integrating:
\begin{equation}\label{eq67}
    \mathbf{x}_\alpha=\mathbf{x}_s s_\alpha=\frac{L(t)}{2\pi}(\cos{\theta(\alpha,t)},\sin{\theta(\alpha,t)}).
\end{equation}
\paragraph{\textbf{Small scale decomposition (SSD)}}
The idea of the small scale decomposition (SSD) is to extract the dominant part of the equations at small spatial scales \cite{hou1994removing}. To remove the stiffness, we use SSD in our problem and develop an explicit, non-stiff time integration algorithm.  In Eqs. \eqref{bimnut} and \eqref{bimpressure}, based on  the analysis of the single-layer- and double-layer- type terms, the only singularity in the integrands comes from the logarithmic kernel. Following \cite{hou1994removing} and noticing the curvature term in Eq. \eqref{pressurefield}, one can show that at small spatial scales,
\begin{equation}\label{eq68}
V(\alpha,t) \sim \frac{1}{s_\alpha^2} \mathcal{H}[\theta_{\alpha\alpha}],
\end{equation}
where $\mathcal{H}(\xi)=\frac{1}{2\pi}\int_0^{2\pi}\xi'\cot{\frac{\alpha-\alpha'}{2}}d\alpha'$ is the Hilbert transform for a $2\pi$-periodic function $\xi$.\\ We rewrite Eq. \eqref{eq62},
\begin{equation}\label{eq69}
\theta_t=\frac{1}{s_\alpha^3} \mathcal{H}[\theta_{\alpha\alpha\alpha}]+N(\alpha,t),
\end{equation}
where the Hilbert transform term is the dominating high-order term at small spatial scales, and $ \displaystyle N= (\kappa T-V_s)-\frac{1}{s_\alpha^3} \mathcal{H}[\theta_{\alpha\alpha\alpha}]$ contains other lower-order terms in the evolution. This demonstrates that an explicit time-stepping method has the high-order constraint $\displaystyle \Delta t \le \left ( \frac{h}{s_\alpha} \right)^3$ where $\Delta t$ and $h$ are the time-step and spatial grid size, respectively. This has been demonstrated numerically in the seminal work \cite{hou1994removing} for a Hele-Shaw problem.  For the tumor growth problem, the semi-implicit time-stepping scheme (see Eq. \eqref{eq69}) requires $\Delta t = O(h)$ instead of explicit schemes which would require $\Delta t = O(h^3)$.  

\section{Semi-implicit Time-Stepping Scheme}
\label{Appendix:Time stepping}
Taking the Fourier transform of Eq. \eqref{eq69}, we get
\begin{equation}\label{tangent angle ode}
{\hat \theta}_t=-\frac{|k|^3}{s_\alpha^3} {\hat \theta}(k,t) +{\hat N}(k,t).
\end{equation}
We solve Eq. \eqref{tangent angle ode} using the second order accurate linear propagator method in the Adams-Bashforth form \cite{hou1994removing} in Fourier space and apply the inverse Fourier transform to recover $\theta$.  Specifically,  we discretize Eq. \eqref{tangent angle ode} as
\begin{equation}\label{linear propagator method}
\quad\quad\;{\hat \theta}^{n+1}(k)=e_k(t_n,t_{n+1}){\hat \theta}^{n}(k)+\frac{\Delta t}{2}(3e_k(t_n,t_{n+1}){\hat N}^{n}(k)-e_k(t_{n-1},t_{n+1}){\hat N}^{n-1}(k),
\end{equation}
where  the superscript $n$ denotes the numerical solutions at $t=t_n$ and the integrating factor
\begin{equation}\label{eq72}
e_k(t_1,t_2)=\exp\left (-{|k|^3}\int_{t_1}^{t_2}\frac{dt}{s_\alpha^3(t)}\right ).
\end{equation}
Note that by setting the integrating factors in Eq. \eqref{linear propagator method} to $1$, we recover the Adams-Bashforth explicit time-stepping method.
The integrating factors in Eq. \eqref{linear propagator method} can be evaluated simply using the trapezoidal rule,
\begin{eqnarray}\label{eq73}
\int_{t_n}^{t_{n+1}}\frac{dt}{s_\alpha^3(t)} &\approx& \frac{\Delta t}{2} \left (\frac{1}{(s_\alpha^n)^3}+\frac{1}{(s_\alpha^{n+1})^3} \right ) \nonumber, \\
\int_{t_{n-1}}^{t_{n+1}}\frac{dt}{s_\alpha^3(t)} &\approx& {\Delta t} \left (\frac{1}{2(s_\alpha^{n-1})^3}+\frac{1}{(s_\alpha^{n})^3}+\frac{1}{2(s_\alpha^{n+1})^3} \right ).
\end{eqnarray}
To compute the arclength $s_\alpha$, Eq. \eqref{eq63} is discretized using the explicit second-order Adams-Bashforth method \cite{hou1994removing},
\begin{equation}\label{linear propagator method 2}
s_\alpha^{n+1}=s_\alpha^n+\frac{\Delta t}{2}(3M^n-M^{n-1}),
\end{equation}
where $M$ is calculated using
\begin{equation}\label{eq75}
M=\frac{1}{2\pi}\int_0^{2\pi}V(\alpha,t)\theta_\alpha d\alpha.
\end{equation}

 Note that the second order linear propagator and Adams-Bashforth methods are multi-step method and require two previous time steps. The first time step is realized using an explicit Euler method for $s_\alpha^1$ and a first-order linear propagator of a similar form for $\hat{\theta}^1$.

 To reconstruct the tumor-host interface $(x(\alpha,t_{n+1}),y(\alpha,t_{n+1}))$ from the updated $\theta^{n+1}(\alpha)$ and $s_\alpha^{n+1}$, we first update a reference point $(x(0,t_{n+1}),y(0,t_{n+1})$ using a second-order explicit Adams-Bashforth method to discretize the equation of motion $\mathbf{x}_t=V\hat{\mathbf{n}}$ with the tangential part dropped since it does not change the morphology:
 \begin{equation}\label{eq76}
     (x(0,t_{n+1}),y(0,t_{n+1}))=(x(0,t_{n}),y(0,t_{n}))+\frac{\Delta t}{2} \left( 3V(0,t_n)\hat{\mathbf{n}}(0,t_{n})-V(0,t_{n-1})\hat{\mathbf{n}}(0,t_{n-1}) \right).
 \end{equation}
 Once we update the reference point, we obtain the configuration of the interface from the $\theta^{n+1}(\alpha)$ and $s_\alpha^{n+1}$ by integrating Eq. \eqref{eq67} following \cite{hou1994removing}:
 \begin{equation}
     \begin{aligned}\label{eq77}
     x(\alpha,t_{n+1})&=&x(0,t_{n+1})+s_\alpha^{n+1}\left( \int_0^{\alpha} \cos(\theta^{n+1}(\alpha'))d\alpha'
     -\frac{\alpha}{2\pi}\int_0^{2\pi}\cos(\theta^{n+1}(\alpha'))d\alpha'\right),\nonumber\\
     y(\alpha,t_{n+1})&=&y(0,t_{n+1})+s_\alpha^{n+1}\left( \int_0^{\alpha} \sin(\theta^{n+1}(\alpha'))d\alpha'
     -\frac{\alpha}{2\pi}\int_0^{2\pi}\sin(\theta^{n+1}(\alpha'))d\alpha'\right),
 \end{aligned}
 \end{equation}
 where the indefinite integration is performed using the discrete Fourier transform.

 We use a 25th order Fourier filter to damp the highest nonphysical mode and suppress the  aliasing error \cite{hou1994removing}. We also use Krasny filtering \cite{krasny1986study} to prevent the accumulation of round-off errors during the computation.

We solve first the nutrient field $\sigma$ then the pressure field $p$. Next we compute the normal velocity $V$ and update the interface $\Gamma(t)$ and repeat this procedure.

 \bibliographystyle{elsarticle-num} 
 \bibliography{cas-refs}

\begin{thebibliography}{10}
\expandafter\ifx\csname url\endcsname\relax
  \def\url#1{\texttt{#1}}\fi
\expandafter\ifx\csname urlprefix\endcsname\relax\def\urlprefix{URL }\fi
\expandafter\ifx\csname href\endcsname\relax
  \def\href#1#2{#2} \def\path#1{#1}\fi

\bibitem{cristini2005morphologic}
V.~Cristini, H.~B. Frieboes, R.~Gatenby, S.~Caserta, M.~Ferrari, J.~Sinek,
  Morphologic instability and cancer invasion, Clinical Cancer Research 11~(19)
  (2005) 6772--6779.

\bibitem{araujo2004history}
R.~P. Araujo, D.~S. McElwain, A history of the study of solid tumour growth:
  the contribution of mathematical modelling, Bulletin of mathematical biology
  66~(5) (2004) 1039--1091.

\bibitem{fasano2006mathematical}
A.~Fasano, A.~Bertuzzi, A.~Gandolfi, Mathematical modelling of tumour growth
  and treatment, in: Complex systems in biomedicine, Springer, 2006, pp.
  71--108.

\bibitem{roose2007mathematical}
T.~Roose, S.~J. Chapman, P.~K. Maini, Mathematical models of avascular tumor
  growth, SIAM review 49~(2) (2007) 179--208.

\bibitem{bellomo2008selected}
N.~Bellomo, E.~de~Angelis, Selected topics in cancer modeling: genesis,
  evolution, immune competition, and therapy, Springer Science \& Business
  Media, 2008.

\bibitem{lowengrub2009nonlinear}
J.~S. Lowengrub, H.~B. Frieboes, F.~Jin, Y.-L. Chuang, X.~Li, P.~Macklin, S.~M.
  Wise, V.~Cristini, Nonlinear modelling of cancer: bridging the gap between
  cells and tumors, Nonlinearity 23~(1) (2009) R1.

\bibitem{byrne2010dissecting}
H.~M. Byrne, Dissecting cancer through mathematics: from the cell to the animal
  model, Nature Reviews Cancer 10~(3) (2010) 221.

\bibitem{byrne2012mathematical}
H.~M. Byrne, Mathematical biomedicine and modeling avascular tumor growth
  (2012).

\bibitem{KimOthmer2015}
Y.~Kim, H.~Othmer, Hybrid models of cell and tissue dynamics in tumor growth,
  Math. Biosci. Eng. 12 (2015) 1141--1156.

\bibitem{alfonso2017biology}
J.~Alfonso, K.~Talkenberger, M.~Seifert, B.~Klink, A.~Hawkins-Daarud,
  K.~Swanson, H.~Hatzikirou, A.~Deutsch, The biology and mathematical modelling
  of glioma invasion: a review, Journal of the Royal Society Interface 14~(136)
  (2017) 20170490.

\bibitem{Yankeelov2018}
A.~Jarrett, E.~Lima, D.~n. Hormuth, M.~McKenna, F.~X., E.~D.A., A.~Resende,
  B.~A., T.~Yankeelov, Mathematical models of tumor cell proliferation: A
  review of the literature, Expert Rev. Anticancer Ther. 18 (2018) 1271--1286.

\bibitem{cristini2010multiscale}
V.~Cristini, J.~Lowengrub, Multiscale modeling of cancer: an integrated
  experimental and mathematical modeling approach, Cambridge University Press,
  2010.

\bibitem{cristini2017introduction}
V.~Cristini, E.~Koay, Z.~Wang, An Introduction to Physical Oncology: How
  Mechanistic Mathematical Modeling Can Improve Cancer Therapy Outcomes, CRC
  Press, 2017.

\bibitem{greenspan1976growth}
H.~Greenspan, On the growth and stability of cell cultures and solid tumors,
  Journal of theoretical biology 56~(1) (1976) 229--242.

\bibitem{friedman2001existence}
A.~Friedman, F.~Reitich, On the existence of spatially patterned dormant
  malignancies in a model for the growth of non-necrotic vascular tumors,
  Mathematical Models and Methods in Applied Sciences 11~(04) (2001) 601--625.

\bibitem{friedman2001symmetry}
A.~Friedman, F.~Reitich, Symmetry-breaking bifurcation of analytic solutions to
  free boundary problems: an application to a model of tumor growth,
  Transactions of the American Mathematical Society 353~(4) (2001) 1587--1634.

\bibitem{friedman2006bifurcation}
A.~Friedman, B.~Hu, Bifurcation from stability to instability for a free
  boundary problem arising in a tumor model, Archive for rational mechanics and
  analysis 180~(2) (2006) 293--330.

\bibitem{Hu20071}
A.~Friedman, B.~Hu, Bifurcation for a free boundary problem modeling tumor
  growth by stokes equation, SIAM Journal on Mathematical Analysis 39~(1)
  (2007) 174--194.
\newblock \href {https://doi.org/10.1137/060656292}
  {\path{doi:10.1137/060656292}}.

\bibitem{Hu20072}
A.~Friedman, B.~Hu, Bifurcation from stability to instability for a free
  boundary problem modeling tumor growth by stokes equation, Journal of
  mathematical analysis and applications 327~(1) (2007) 643--664.

\bibitem{friedman2008stability}
A.~Friedman, B.~Hu, Stability and instability of liapunov-schmidt and hopf
  bifurcation for a free boundary problem arising in a tumor model,
  Transactions of the American Mathematical Society 360~(10) (2008) 5291--5342.

\bibitem{zhao2020symmetry}
X.~E. Zhao, B.~Hu, Symmetry-breaking bifurcation for a free-boundary tumor
  model with time delay, Journal of Differential Equations 269~(3) (2020)
  1829--1862.

\bibitem{zhao2020impact}
X.~E. Zhao, B.~Hu, The impact of time delay in a tumor model, Nonlinear
  Analysis: Real World Applications 51 (2020) 103015.

\bibitem{cristini2003nonlinear}
V.~Cristini, J.~Lowengrub, Q.~Nie, Nonlinear simulation of tumor growth,
  Journal of mathematical biology 46~(3) (2003) 191--224.

\bibitem{cristini2009}
V.~Cristini, X.~Li, L.~J.S., S.~Wise, Nonlinear simulations of solid tumor
  growth using a mixture model: Invasion and branching, J. Math. Biol. 4-5~(06)
  (2009) 723--763.

\bibitem{fritz2019local}
M.~Fritz, E.~A. Lima, V.~Nikoli{\'c}, J.~T. Oden, B.~Wohlmuth, Local and
  nonlocal phase-field models of tumor growth and invasion due to ecm
  degradation, Mathematical Models and Methods in Applied Sciences 29~(13)
  (2019) 2433--2468.

\bibitem{mcdougall2002mathematical}
S.~R. McDougall, A.~Anderson, M.~Chaplain, J.~Sherratt, Mathematical modelling
  of flow through vascular networks: implications for tumour-induced
  angiogenesis and chemotherapy strategies, Bulletin of mathematical biology
  64~(4) (2002) 673--702.

\bibitem{mcdougall2006mathematical}
S.~R. McDougall, A.~R. Anderson, M.~A. Chaplain, Mathematical modelling of
  dynamic adaptive tumour-induced angiogenesis: clinical implications and
  therapeutic targeting strategies, Journal of theoretical biology 241~(3)
  (2006) 564--589.

\bibitem{foo2011stochastic}
J.~Foo, K.~Leder, F.~Michor, Stochastic dynamics of cancer initiation, Physical
  biology 8~(1) (2011) 015002.

\bibitem{hillen2013tumor}
T.~Hillen, H.~Enderling, P.~Hahnfeldt, The tumor growth paradox and immune
  system-mediated selection for cancer stem cells, Bulletin of mathematical
  biology 75~(1) (2013) 161--184.

\bibitem{pham2018nonlinear}
K.~Pham, E.~Turian, K.~Liu, S.~Li, J.~Lowengrub, Nonlinear studies of tumor
  morphological stability using a two-fluid flow model, Journal of mathematical
  biology (2018) 1--39.

\bibitem{cui2001analysis}
S.~Cui, A.~Friedman, Analysis of a mathematical model of the growth of necrotic
  tumors, Journal of Mathematical Analysis and Applications 255~(2) (2001)
  636--677.

\bibitem{hao2012bifurcation}
W.~Hao, J.~D. Hauenstein, B.~Hu, Y.~Liu, A.~J. Sommese, Y.-T. Zhang,
  Bifurcation for a free boundary problem modeling the growth of a tumor with a
  necrotic core, Nonlinear Analysis: Real World Applications 13~(2) (2012)
  694--709.

\bibitem{kohlmann2012necrotic}
M.~Kohlmann, Necrotic tumor growth: An analytic approach, Acta biotheoretica
  60~(3) (2012) 273--287.

\bibitem{hao2012continuation}
W.~Hao, J.~D. Hauenstein, B.~Hu, Y.~Liu, A.~J. Sommese, Y.-T. Zhang,
  Continuation along bifurcation branches for a tumor model with a necrotic
  core, Journal of Scientific Computing 53~(2) (2012) 395--413.

\bibitem{wu2019bifurcation}
J.~Wu, Bifurcation for a free boundary problem modeling the growth of necrotic
  multilayered tumors, Discrete \& Continuous Dynamical Systems-A 39~(6) (2019)
  3399.

\bibitem{zhuang2018analysis}
Y.~Zhuang, S.~Cui, Analysis of a free boundary problem modeling the growth of
  multicell spheroids with angiogenesis, Journal of Differential Equations
  265~(2) (2018) 620--644.

\bibitem{song2021stationary}
H.~Song, B.~Hu, Z.~Wang, Stationary solutions of a free boundary problem
  modeling the growth of vascular tumors with a necrotic core, Discrete \&
  Continuous Dynamical Systems-B 26~(1) (2021) 667.

\bibitem{lu2020complex}
M.-J. Lu, C.~Liu, J.~Lowengrub, S.~Li, Complex far-field geometries determine
  the stability of solid tumor growth with chemotaxis, Bulletin of mathematical
  biology 82~(3) (2020) 1--41.

\bibitem{macklin2007nonlinear}
P.~Macklin, J.~Lowengrub, Nonlinear simulation of the effect of
  microenvironment on tumor growth, Journal of theoretical biology 245~(4)
  (2007) 677--704.

\bibitem{friedman2006cancer}
A.~Friedman, Cancer models and their mathematical analysis, in: Tutorials in
  Mathematical Biosciences III, Springer, 2006, pp. 223--246.

\bibitem{grimes2014method}
D.~R. Grimes, C.~Kelly, K.~Bloch, M.~Partridge, A method for estimating the
  oxygen consumption rate in multicellular tumour spheroids, Journal of The
  Royal Society Interface 11~(92) (2014) 20131124.

\bibitem{roussos2011chemotaxis}
E.~T. Roussos, J.~S. Condeelis, A.~Patsialou, Chemotaxis in cancer, Nature
  Reviews Cancer 11~(8) (2011) 573--587.

\bibitem{veerapaneni2016integral}
S.~Veerapaneni, Integral equation methods for vesicle electrohydrodynamics in
  three dimensions, Journal of Computational Physics 326 (2016) 278--289.

\bibitem{jou1997}
H.~Jou, P.~H. Leo, J.~Lowengrub, Microstructural evolution in inhomogeneous
  elastic media, Journal of Computational Physics 131~(1) (1997) 109--148.

\bibitem{hou1994removing}
T.~Y. Hou, J.~S. Lowengrub, M.~J. Shelley, Removing the stiffness from
  interfacial flows with surface tension, Journal of Computational Physics
  114~(2) (1994) 312--338.

\bibitem{Pedro1}
P.~H.~A. Anjos, S.~Li, Weakly nonlinear analysis of the saffman-taylor problem
  in a radially spreading fluid annulus, Phys. Rev. Fluids 5 (2020) 054002.
\newblock \href {https://doi.org/10.1103/PhysRevFluids.5.054002}
  {\path{doi:10.1103/PhysRevFluids.5.054002}}.

\bibitem{Pedro2}
M.~Zhao, P.~H.~A. Anjos, J.~Lowengrub, S.~Li, Pattern formation of the
  three-layer saffman-taylor problem in a radial hele-shaw cell, Phys. Rev.
  Fluids 5 (2020) 124005.
\newblock \href {https://doi.org/10.1103/PhysRevFluids.5.124005}
  {\path{doi:10.1103/PhysRevFluids.5.124005}}.

\bibitem{minjhe2019}
M.-J. Lu, C.~Liu, S.~Li, Nonlinear simulation of an elastic tumor-host
  interface, Computational and Mathematical Biophyics. 7~(1) (2019) 25--47.

\bibitem{mullins1963morphological}
W.~W. Mullins, R.~F. Sekerka, Morphological stability of a particle growing by
  diffusion or heat flow, Journal of applied physics 34~(2) (1963) 323--329.

\bibitem{li2011boundary}
S.~Li, X.~Li, A boundary integral method for computing the dynamics of an
  epitaxial island, SIAM Journal on Scientific Computing 33~(6) (2011)
  3282--3302.

\bibitem{Wenrui2018}
W.~Hao, B.~Hu, S.~Li, L.~Song, Convergence of boundary integral method for a
  free boundary system, Journal of Computational and Applied Mathematics 334
  (2018) 128 -- 157.
\newblock \href {https://doi.org/https://doi.org/10.1016/j.cam.2017.11.016}
  {\path{doi:https://doi.org/10.1016/j.cam.2017.11.016}}.

\bibitem{kress1995numerical}
R.~Kress, On the numerical solution of a hypersingular integral equation in
  scattering theory, Journal of computational and applied mathematics 61~(3)
  (1995) 345--360.

\bibitem{saad1986gmres}
Y.~Saad, M.~H. Schultz, Gmres: A generalized minimal residual algorithm for
  solving nonsymmetric linear systems, SIAM Journal on scientific and
  statistical computing 7~(3) (1986) 856--869.

\bibitem{baker1990connection}
G.~Baker, M.~Shelley, On the connection between thin vortex layers and vortex
  sheets, Journal of Fluid Mechanics 215 (1990) 161--194.

\bibitem{krasny1986study}
R.~Krasny, A study of singularity formation in a vortex sheet by the
  point-vortex approximation, Journal of Fluid Mechanics 167 (1986) 65--93.

\end{thebibliography}





\end{document}